\input amstex
\documentstyle{amsppt}
\magnification=1200
\overfullrule=0pt
\hsize = 5.5 true in\hoffset=.5in
\vsize = 7 true in\voffset=.5in
\loadeusm

\def\fn{{\frak n}}
\def\fnl{{\fn_L}}
\def\fnm{{\fn_M}}
\def\fnn{{\fn_N}}
\def\fnnp{{\fn_N'}}

\def\fS*{{\frak S_*}}
\def\fs{{\frak s}}
\def\fsp{{\frak s'}}

\def\fW{{\frak W}}

\def\sB{{\eusm B}}

\def\sC{{\eusm C}}

\def\sL{{\eusm L}}

\def\sM{{\eusm M}}

\def\tM{{\widetilde {\sM}}}

\def\sU{{\eusm U}}

\def\sR{{\eusm R}}

\def\C{{\Bbb C}}

\def\R{{\Bbb R}}
\def\T{{\Bbb T}}
\def\Z{{\Bbb Z}}

\def\a{{\alpha}}
\def\da{{\dot\a}}
\def\b{{\beta}}
\def\db{{\dot\b}}
\def\f{{\varphi}}
\def\p{{\psi}}
\def\z{{\zeta}}
\def\Ad{{\text{\rm Ad}}}
\def\tp{{\tilde p}}
\def\tq{{\tilde  q}}
\def\tr{{\tilde r}}

\def\txti{{\text{\rm i}}}
\def\txm{{\text{\rm m}}}

\def\txmm{{\txm_M}}
\def\txmn{{\txm_N}}
\def\txn{{\text{\rm n}}}

\def\txnn{{\txn_N}}
\def\txnsps{{\txn_{\fsp, \fs}}}

\def\txs{{\text{\rm s}}}

\def\wtQ{{\widetilde Q}}

\def\wtG{{\widetilde G}}
\def\wtH{{\widetilde H}}

\def\th{{\theta}}
\def\tht{{\theta_T}}
\def\tht'{{\th_{T'}}}
\def\la{{\lambda}}

\def\AFD{{\text{\rm AFD}}}
\def\Aut{{\text{\rm Aut}}}
\def\Autf'{{\Aut_\f'}}
\def\Cnt{{\text{\rm Cnt}}}
\def\Ob{{\text{\rm Ob}}}
\def\Obm{{\Ob_{\text{\rm m}}}}
\def\cnt{{\Cnt}}
\def\cntr{{\Cnt_{\text{\rm r}}}}
\def\Hom{{\text{\rm Hom}}}

\def\homr{{\Hom_\R}}

\def\Int{{\text{\rm Int}}}

\def\Ker{{\text{\rm Ker}}}

\def\mod{{\text{\rm mod}}}

\def\Out{{\text{\rm Out}}}
\def\out{{\text{\rm out}}}

\def\Outt{{\Out_{\tau, \th}}}

\def\Autt{{\Aut_{\tau, \th}}}

\def\Autf{{\Aut_\f}}

\def\Autf'{{\Aut_\f'}}

\def\two{{\rm I\!I}}
\def\twoone{{\rm I\!I$_1$}}
\def\threee{{\text{\rm I\!I\!I}}}
\def\threeone{{\rm I\!I\!I$_1$}}
\def\threel{{\rm I\!I\!I$_{\lambda}$}}
\def\three0{{\rm I\!I\!I$_0$}}
\def\threez{{\three0}}
\def\twoinf{{\rm I\!I$_{\infty}$}}
\def\oneinf{{\rm I$_{\infty}$}}

\def\tsU{{\widetilde {\sU}}}

\def\tAd{{\widetilde {\text{\rm Ad}}}}
\def\ta{{\widetilde {\a}}}

\def\d{{\delta}}

\def\tj{{\tilde j}}

\def\tmu{{\tilde \mu}}

\def\La{{\Lambda}}

\def\part{{\partial}}
\def\dpart{{\dot\part}}
\def\partg{{\part_G}}
\def\parth{{\part_H}}
\def\partq{{\part_Q}}

\def\partwtq{{\part_{\wtQ}}}

\def\partth{{\part_\th}}
\def\dpartth{{\dpart_\th}}

\def\oct{{\text{\rm Oct}}}

\def\act{{\text{\rm Act}}}
\def\sig{{\sigma}}
\def\tsig{{\tilde \sig}}

\def\sigft{{\sigma_t^{\f}}}

\def\sigsps{{\sig_{\fsp, \fs}}}
\def\wt{{semi-finite normal weight}}
\def\fwt{{faithful \wt}}
\def\botimes{{\bar \otimes}}
\def\id{{\text{\rm id}}}

\def\r0{{\sR_0}}

\def\r01{{\sR_{0,1}}}

\def\Map{{\text{\rm Map}}}

\def\botimes{{\overline \otimes}}
\def\wt{{semi-finite normal weight}}
\def\wts{{semi-finite normal weights}}
\def\fwt{{faithful \wt}}
\def\fwts{{faithful \wts}}

\def\vna{{von Neumann algebra}}
\def\vnas{{ von Neumann algebras}}

\def\tpart{{\widetilde\partial}}

\def\tB{{\text{\rm B}}}

\def\tC{{\text{\rm C}}}

\def\tH{{\text{\rm H}}}

\def\tZ{{\text{\rm Z}}}

\def\QED{{\hfill$\heartsuit$}}

\def\dprime{{^{\prime\prime}}}
\def\cSQN#1#2#3#4#5#6#7#8#9{{
\CD
@. 1 @. 1 @. 1 @.\\
@. @VVV @VVV @VVV @.\\
1 @>>> #1 @>>> #2 @>\partial>> #3 @>>> 1 \\ 
@.@VVV @VVV @VVV @. \\   
1 @>>>#4
@>>>#5 @>\widetilde
\partial>>#6 @>>>1\\  
@. @VV{\Ad}V @VV{\widetilde{\text{\rm Ad}}}V @VVV
@.\\   
1 @>>> #7 @>>> #8 @>\dot \part>> #9@>>> 1\\   
@.@VVV @VVV @VVV @. \\  
@. 1 @. 1 @. 1 @.
\endCD
}}

\def\inv{{^{-1}}}
\document

\def\X{{\text{\rm Xext}}}

\def\res{{\text{\rm res}}}
\def\Res{{\text{\rm Res}}}
\def\rsn{{\text{\rm Rsn}}}
\def\inf{{\text{\rm inf}}}

\def\dfs{{\dot\fs}}

\def\fsg{{\fs_G}}
\def\sdpart{{\fs_{\dot\part}}}
\def\sh{{\fs\!_{\scriptscriptstyle H}}}

\def\sp{{\fs_\pi}}

\def\sj{{\fs_j}}

\def\se{{\fs_{\!\scriptscriptstyle E}}}

\def\sts{{\fs_{\tilde\sig}}}

\def\pig{{\pi\!_{\scriptscriptstyle G}}}

\def\piq{{\pi_Q}}
\def\piz{{\pi\!_{\scriptscriptstyle\tZ}}}
\def\dpi{{\dot \pi}}
\def\sz{{\fs_{\!\scriptscriptstyle\tZ}}}

\def\etat{{{\eta_{T}}}}

\def\etat'{{\eta_{T'}}}

\def\tpi{{\tilde \pi}}

\def\piq{{\pi_Q}}

\def\_#1{{_{#1}}}
\def\^#1{{^{#1}}}

\def\ainv#1{{\a_{#1}^{-1}}}

\def\(#1){{^{({#1})}}}

\def\scirc{{\lower-.3ex\hbox{{$\scriptscriptstyle\circ$}}}}

\def\bracett'#1{{\{#1\}_{T'}}}
\def\brackett#1{{\left[{#1}\right]_{T}}}
\def\brackett'#1{{\left[{#1}\right]_{T'}}}

\def\rt'z{{\R/T'\Z}}
\def\coron{{\!:\ }}

\def\cdag{{\text{countable discrete amenable group}}}

\def\hjr{{{\text{\rm HJR}-exact sequence}}}
\def\tbatw{{\tB_\a^2}}
\def\tcatw{{\tC_\a^2}}
\def\tzatw{{\tZ_\a^2}}
\def\thatw{{\tH_\a^2}}
\def\tbath{{\tB_\a^3}}

\def\tzath{{\tZ_\a^3}}
\def\tzasth{{\tZ_{\a, \text{\rm s}}^3}}
\def\tbasth{{\tB_{\a, \text{\rm s}}^3}}
\def\thasth{{\tH_{\a, \text{\rm s}}^3}}
\def\thath{{\tH_\a^3}}
\def\tbth1{{\tB_\th^1}}

\def\tbasout{{\tB_{\a, \fs}^\out}}
\def\tctw{{\tC^2}}
\def\tzth1{{\tZ_\th^1}}
\def\thth1{{\tH_\th^1}}
\def\ththr{{\tH^3}}
\def\ththra{{\tH_\a^3}}

\def\thaout{{\tH_\a^\out}}

\def\thasout{{\tH_{\a, \fs}^\out}}

\def\tzthr{{\tZ^3}}
\def\tzthra{{\tZ_\a^3}}

\def\tbthr{{\tB^3}}

\def\thtw{{\tH^2}}
\def\thatw{{\tH_\a^2}}
\def\tztw{{\tZ^2}}
\def\tzatw{{\tZ_\a^2}}

\def\tzasout{{\tZ_{\a, \fs}^\out}}

\def\chim{{\chi_{\text{\rm m}}}}

\def\pr{{\text{\rm pr}}}

\def\pijj'{{\pi_{J, J'}}}
\def\pijj's{{\pi_{J, J'}^*}}

\def\tmu{{\tilde \mu}}

\def\linfrt'z{{L^\infty(\R/T'\Z)}}

\def\rt'z{{\R/T'\Z}}

\def\muh{{\mu\!_{\scriptscriptstyle H}}}
\def\dhjr{{\d_{\scriptscriptstyle{\text{\rm HJR}}}}}
\def\Inf{{\text{\rm Inf}}}

\nologo
\TagsOnRight
\loadeusm

\topmatter
\title
{OUTER ACTIONS OF A DISCRETE AMENABLE GROUP 
ON APPROXIMATELY FINITE
DIMENSIONAL FACTORS I,\\
General Theory}
\endtitle
\rightheadtext{Outer Conjugacy}
\leftheadtext{Outer Conjugacy}
\author
{Yoshikazu Katayama and Masamichi Takesaki}
\endauthor
\address{Department of Mathematics,
Osaka Ky$\hat\text{o}$iku University Osaka, Japan}
\endaddress
\address{Department of Mathematics,
University of California,
Los Angeles, California 90095-1555}
\endaddress
\thanks{This research is support in part by the NSF Grant: DMS-9801324 and 
DMS-0100883  and
also by the Ministry of Education, Science, Sports and Culture, Grant-in-Aid for Scientific 
Research (C), 14540206, 2002.}
\endthanks
\abstract
{To each factor $\sM$, we associate an invariant $\Obm(\sM)$  to  be called the
{\it intrinsic modular obstruction} as a cohomological invariant which lives in
the ``third" cohomology group:
$$\aligned
&\thasout(\Out(\sM)\times\R, \thth1(\R, \sU(\sC)), \sU(\sC))
\endaligned
$$
where $\{\sC, \R, \th\}$ is the flow of weights on $\sM$. If $\a$ is an
 outer action of  a countable discrete group $G$ on {$\sM$}, then its modulus
$\mod(\a)\in \Hom(G, \Aut_\th(\sC))$, $N=\a\inv(\cntr(M))$ 
and the pull back 
$$
\Obm(\a)=\a^*(\Obm(\sM))\in \thasout(G\times\R, N, \sU(\sC)))
$$
to be called  the {\it modular obstruction} of $\a$ are invariants of the outer conjugacy
class of the outer action $\a$.

We prove that if the factor $\sM$ is approximately finite dimensional and $G$ is amenable,
then the invariants uniquely determine the outer conjugacy class of $\a$ and the
every invariant occurs as the invariant of an outer action $\a$ of $G$ on $\sM$. In the
case that $\sM$ is a factor of type  {\threel}, $0<\la\leq 1$, the modular obstruction 
group $\thasout(G\times\R, N, \sU(\sC))$ and the modular obstruction $\Obm(\a)$ take
simpler forms. These together with examples will be discussed in the forthcoming
paper, \cite {KtT2}.}
\endabstract
\endtopmatter

\document

\head{\bf \S 0. Introduction}
\endhead
With the successful completion of the cocycle conjugacy classification of amenable
discrete group actions on {\AFD} factors by many hands over more than two decades, \cite
{C3, J1, JT, O, ST1, ST2, KwST, KtST1}, it is only naturally to consider the outer conjugacy
classification of amenable discrete group outer actions on {\AFD} factors. In fact, the
work on the program has been already started by the pioneering works of Connes,
\cite{Cnn 3, 4, 6}, Jones \cite{J1} and Ocneanu \cite{Ocn}. In this article, we complete the outer conjugacy
classification of discrete amenable group outer actions on AFD factors. The cases of type I, {\twoone} and
{\twoinf} with additional technical assumption were already completed by Jones, \cite{J1}, and Ocneanu,
\cite{Ocn}, so the case of type \threee\ will be mainly considered although the technical
assumption in the case of type {\twoinf} placed in the work of Ocneanu \cite{Ocn} must be
removed. 

As in the case of the cocycle conjugacy classification, we first associate invariants which
are intrinsic to any factor $\sM$, the flow of weights, the modulus, the characteristic
square and the modular obstruction $\Obm(\sM)$. Then the outer conjugacy invariants are
given by the pull back of these intrinsic quantities of the factor by the outer action. To be
more precise, let $\sM$ be a separable factor. Associated with $\sM$ is the characteristic
square:
$$
\cSQN{\T}{\sU(\sC)}{\tB_\th^1}{\sU(\sM)}{\tsU(\sM)}{\tZ_\th^1}{\Int(\sM)}{\cntr(\sM)}
{\tH_\th^1}
$$
which is equivariant under $\Aut(\sM)\times \R$. The middle vertical exact sequence is
the source of the intrinsic invariant:
$$
\Theta(\sM)\in \La_{\mod\times \th}(\Aut(\sM)\times \R, \cntr(\sM), \sU(\sC)).
$$
To avoid heavy notations and to see the essential mechanism governing the above exact
characteristic square, let us consider the situation that a group $H$ equipped with a
distinguished pair of normal subgroups $M\subset L\subset H$ which acts on the ergodic
flow $\{\sC, \R, \th\}$, i.e., the action $\a$ of $H$ on $\sC$ commutes with the flow 
$\th$. Assume that the normal subgroup $L$ does not act on $\sC$, i.e., $L\subset
\Ker(\a)$, so that the action $\a$ factors through the both quotient groups $Q=H/L$ and
$G=H/M$. In the case that $H=\Aut(\sM)$, the groups $L$ and $M$ stand for
$\cntr(\sM)$ and $\Int(\sM)$, therefore $Q=\Outt(\sM)$ and $G=\Out(\sM)$. Let $\wtH,
\wtG$ and  $\wtQ$ denote respectively the product groups $H\times \R, G\times \R$ and $Q\times \R$. We
denote the unitary group $\sU(\sC)$ simply by $A$ for the simplicity. In the case that $H=\Aut(\sM)$, then 
we require appropriate Borelness for mappings. But $Q$ can fail to have a reasonable Borel structure, so
we treat $Q$ as a discrete group. On the product group $\wtQ=Q\times \R$, we consider the product Borel
structure as well as the product topology. 

In this circumstance, we will see that  each characteristic cocycle $(\la, \mu)\in \tZ_\a(\wtH,
L, A)$ gives rise to an $\wtH$-equivariant exact square:
$$\eightpoint\CD
@.1@.1@.1\\
@.@VVV@VVV@VVV\\
1@>>>\T@>>>A@>\part>>\tB@>>>1\\
@.@VVV@ViVV@VVV\\
1@>>>U=E^\th@>>>E@>\part_\th>>\tZ@>>>1\\
@.@VVV@VjVV@VVV\\
1@>>>K@>>>L@>\dot\part>>\tH@>>>1\\
@.@VVV@VVV@VVV\\
@.1@.1@.1\\
\endCD
$$
with $E=A\times _\mu L$. The subgroup $K$ of $L$ is normal in $H$ and depends on 
the characteristic invariant $\chi=[\la, \mu]\in \La_\a(\wtH, L, A)$. We denote it by $K(\chi)$
or $K(\la, \mu)$ to indicate the dependence of $K$ on $\chi$ or $(\la, \mu)$. We then define a
subgroup
$\tZ_\a(\wtH, L, M, A)$ of 
$\tZ_\a(\wtH, L, A)$ to be the subgroup consisting of those $(\la, \mu)\in \tZ_\a(\wtH, L, A)$
such that $M\subset K(\la, \mu)$ and $\La_\a(\wtH, L, M, A)$ to be 
$$
\La_\a(\wtH, L, M, A)=\{\chi\in \La_\a(\wtH, L, A): M \subset K(\chi)\}.
$$

In order to study the outer conjugacy class of an outer action $\a$ of $G$ on a factor $\sM$,
we need to fix a cross-section $\fs: Q\mapsto G$ of the quotient map $\pi: G\mapsto Q$ with
kernel $N=L/M=\Ker(\pi)$ and also to restrict the group of $A$-valued 3-cocycles on $\wtQ$
to the group $\tzasth(\wtQ, A)$ of {\it standard } cocycles and a smaller coboundary group:
$$
\tbasth(\wtQ, A)=\part_\wtQ (\tbath(Q, A))
$$
and to form the quotient group:
$$
\thasth(\wtQ, A)=\tzasth(\wtQ, A)/\tbasth(\wtQ, A).
$$
The cross-section $\fs$ gives rise to a link between the group $\thasth(\wtQ, A)$ and the group
$\Hom_G(N, \thth1)$ of equivariant homomorphisms which in turn allows us to define the fiber
product:
$$
\thasout(G\times\R, N, A)=\thasth(\wtQ, A)*_\fs\Hom_G(N, \thth1).
$$
We then show that this group falls in the modified Huesbshmann - Jones - Ratcliffe exact sequence
which sits next to the Huebshmann - Jones - Ratcliffe exact sequence: 
$$\CD
1@.1\\
@VVV@VVV\\
\tH^1(Q, \T)@>\inf>>\tH^1(G, \T)\\
@V\inf VV@V\inf VV\\
\Hom(H, \T)@=\Hom(H, \T)\\
@V\res VV@V\res VV\\
\endCD
$$
$$\CD
\Hom_H(L, \T)@=\Hom_H(L, \T)\\
@V\part VV@V\part VV\\
\tH^2(Q, \T)@>\inf >>\tH^2(G, \T)\\
@V\inf VV@V\inf VV\\
\tH^2(H, \T)@=\tH^2(H, \T)\\
@V\Res VV@V\res VV\\
\La_\a(\wtH, L, M, A)@>\res >>\La(H, M, \T)\\
@V\d VV@V\dhjr VV\\
\thaout(G, N, A)@>\part>>\ththr(G, \T)\\
@V\Inf VV@V\inf VV\\
\tH^3(H, \T)@=\tH^3(H, \T)\\
\endCD
$$

An action $\a$ of $H$ on a factor $\sM$ with $M=\a\inv(\Int(\sM))$ and $L=\a\inv(\cntr(\sM))$
gives rise naturally to the modular characteristic invariant $\chim(\a)\in \La_\a(\wtH, L, M, A)$
and 
$$
\Obm(\a)=\d(\chim(\a))\in \thasout(G\times\R, N, A)=\thasth(\wtQ, A)*_\fs\Hom_G(N, \thth1).
$$
The cohomology element $\Obm(\a)$ will be called the {\it modular obstruction} of the outer action $\a$ of
$G$.

In the original setting, $H=\Aut(\sM)$, the corresponding $\Obm(\a)$ will be denoted by $\Obm(\sM)$ and
called the {\it intrinsic modular obstruction of} {$\sM$}.

In this article, we will prove the following outer conjugacy classification:
\proclaim{Theorem} {\rm i)} If $\a$ is an outer action of a group $G$ on a factor $\sM$, then the
pair of the modulus $\mod(\a)\in \Hom(G, \Aut_\th(\sC))$ of $\a$ and the modular obstruction{\rm:}
$$
\Obm(\a)\in  \thasout(G\times\R, N, A)
$$
is an outer conjugacy invariant of $\a$ with $N=\a\inv(\cntr(\sM))$.

{\rm ii)} If $G$ is a countable discrete amenable group and $\sM$ is an
approximately finite dimensional factor, then the pair 
$(\Obm(\a), \mod(\a))$ is a complete invariant for the outer conjugacy class of $\a$.

{\rm ii)}  With a countable discrete amenable group $G$ and an {\AFD} factor
$\sM$ fixed, every triplet occurs as the invariant of an outer action of  $G$ on the
$\sM$.
\endproclaim

Contrary to the case of the cocycle conjugacy classification, the outer conjugacy classification of 
outer actions of a countable discrete amenable group on an AFD factor will be carried out by a unified
approach without splitting the case base on the type of the base factor. Indeed, the theory is very much
cohomological and therefore algebraic. Nevertheless,  our classification does not fall in the traditional
classification doctrine of Mackey, we will follow the strategy proposed in an earlier work of Katayama -
Sutherland - Takesaki, \cite {KtST1}. Namely, we first introduce a standard Borel structure to the space of outer
actions of a countable discrete group $G$ on a separable factor $\sM$ and associate functorially
invariants in the Borel fashion. 

Most of the mathematical work of this article was carried out during the authors' stay at
the Department of Mathematics, University of Rome ``La Sapienza", in the spring of 2000,
while the second named author stayed there for the entire academic year of 99/00, we
would like to record here our gratitude to Professor S. Doplicher and his colleagues in
Rome for their hospitality extended to the authors. The first named author also would like
to express his thanks to the National Science Foundation for supporting his visit to Rome
from Osaka in Japan.

\head{\bf \S1. Preliminary and Notations}
\endhead
Let {$\sM$} be a separable factor and $G$ a separable locally compact group. We mean by
an {\it outer action} $\a$ of $G$ on $\sM$ a Borel map from $G$ into the group
$\Aut(\sM)$ of automorphisms of {$\sM$} such that 
$$\aligned
&\a_g\circ \a_h\equiv \a_{gh}\quad \mod\ \Int(\sM),\quad g, h \in G,\\
\endaligned\tag1.1
$$
where $\Int(\sM)$ means the group of inner automorphisms.  
 If in addition the following happens
$$\aligned
\a_g&\not\equiv \id\quad \mod\ \Int(\sM)\quad \text{unless } g=1,
\endaligned
$$
then the outer action $\a$ is called {\it free}.

{\smc Remark.} One should not confused an outer action with a {\it free} action of $G$ on {$\sM$}. A
free action of a discrete group $G$ is by definition that a homomorphism $\a: g\in
G\mapsto \a_g\in \Aut(\sM)$ such that $\a_g \notin \Int(\sM), g \neq 1.$ There is no good
definition for the freeness of an action $\a$ of a continuous group $G$. Although one
might take the triviality $\sM'\cap \sM\rtimes_\a G=\C$ of the relative commutant of the
original factor {$\sM$} in the crossed-product as the definition of the freeness of $\a$, which
is an easy consequence of the freeness of $\a$ in the discrete case.

Let $\{\tM, \R, \th, \tau\}$
be the non-commutative flow of weights on {$\sM$} in the sense of Falcone - Takesaki, \cite{FT1},
and $\{\sC, \R, \th\}$ be the Connes - Takesaki flow of weights, \cite{CT}, so that $\sC$
is the center  of $\tM$ and the flow $\{\sC, \R, \th\}$ is the restriction of the non-commutative flow
of weights. The {\vna} $\tM$ is generated by $\sM$ together with one parameter unitary
groups $\{\f^{it}: \f\in \fW_0(\sM), t \in \R\}$, where $\fW_0(\sM)$ means the set of all
{\fwts} on {$\sM$} and $\{\f^{it}\}$'s are related by the Connes cocycle derivatives:
$$
\f^{it}\p^{-it}=(D\f: D\p)_t, \quad \f, \p \in \fW_0(\sM),\ t \in \R.\tag1.2
$$ 
The non-commutative flow $\th$ then acts on $\tM$ by
$$\aligned
&\th_s(x)=x,\quad x\in \sM;\\
&\th_s(\f^{it})=(e^{-s}\f)^{it},\quad \f\in \fW_0(\sM), 
\endaligned \qquad s, t \in \R.\tag1.3
$$
Associated with the non-commutative flow of weights is the extended unitary group
$\tsU(\sM)=\{u\in \sU(\tM): u\sM u^*=\sM\}$. Each $u\in \tsU(\sM)$ gives rise to an
automorphism $\tAd(u)=\Ad(u)|_\sM$ of $\sM$. The set of such automorphisms will be
denoted by $\cntr(\sM)$ and it is a normal subgroup of $\Aut(\sM)$.

An important property of the non-commutative flow
of weights is the relative commutant of $\sM$ in $\tM$:
$$
\sM'\cap \tM=\sC.
$$
A continuous one parameter family $\{u_s\in \sU(\tM): s \in \R\}$ is called {\it $\th$-one
cocycle} if 
$$
u_{s+t}=u_s\th_s(u_t),\quad s, t \in \R.
$$
The set of all $\th$-one cocycles in $\sC$ form a group relative to the pointwise product in 
$\sC$ and is denoted by $\tZ_\th^1$.  

The action $\th$ on $\tM$ is known to be stable in the sense that every $\th$-one cocycle
$\{u_s\}$ is coboundary, i.e.,  there exists $v\in \sU(\tM)$ such that 
$$
u_s=\th_s(v)v^*=(\part v)_s, \quad s \in \R.
$$
The set $\{\part v: v \in \sU(\sC)\}$ of coboundaries is a subgroup of $\tZ_\th^1$ and
denoted by $\tB_\th^1$. The quotient group $\tH_\th^1= \tZ_\th^1/\tB_\th^1$ is an
abelian group which is the first cohomology group of the flow of weights.
The elements of extended unitary
group
$\tsU(\sM)$ are then characterized by the fact that for $u\in \sU(\tM)$:
$$
u\in \tsU(\sM)\quad \Leftrightarrow\quad (\part u)_t \in \sC, t \in \R.
$$
Therefore the map $\part: v\in \tsU(\sM) \mapsto \part v \in \tZ_\th^1$ is surjective. An
important fact about this map is that the exact sequence 
$$
1\ \longrightarrow\ \sU(\sM)\ \longrightarrow\ \tsU(\sM)\ \overset \part \to \longrightarrow\
\tZ_\th^1\ \longrightarrow\ 1
$$
splits equivariantly as soon as a {\fwt} $\f$ is fixed, i.e., to each faithful semi-finite normal weight $\f$ there
corresponds a homomorphism
$b_\f: c\in \tZ_\th^1\mapsto b_\f(c)\in \tsU(\sM)$ such that
$$\aligned
&\tAd(b_\f(c))=\sig_c^\f \quad \text{if}\ \f\ \text{is dominant}, c \in \tZ_\th^1;\\
&b_{\a(\f)}(c)=\a\scirc b_\f\scirc \a\inv(c), \quad c \in \tZ_\th^1,\ \a \in \Aut(\sM);\\
&(D\f: D\p)_c=b_\p(c)b_\f(c^*),\quad \f, \p \in \fW_0(\sM).
\endaligned\tag1.4
$$
This was proven by Falcone - Takesaki \cite{FT2} among other things. \footnote{The coboundary operation
in \cite{FT2} was defined differently as $\part u_s=u\th_s(u^*)$, so in our case the map $b_\f: c \in \tzth1
\mapsto b_\f(c)\in \tsU(\sM)$ behaves as described here.}

The group $\cntr(\sM)$ of ``extended modular" automorphisms is a normal subgroup of
$\Aut(\sM)$ but not closed in the case of type \three0. Nevertheless it is a Borel subgroup
so that its inverse image $N(\a)=\a\inv(\cntr(\sM))$, denoted simply by $N$ in the case that
$\a$ is fixed, is a normal Borel subgroup of the original group $G$. Thus the quotient
group $Q=G/N$ cannot be expected to be a good topological group in general unless $G$
is discrete. Thus we consider mainly discrete groups. Other than the definition of the
invariants of $\a$ we do not have any substantial result on continuous groups any way at
the moment. Interested readers are challenged to go further in the direction of the cocycle
conjugacy problem of one parameter automorphism groups: clearly the very first step
toward the continuous group actions on a factor.

\head{\bf \S2. Modified  Huebschmann Jones Ratcliffe Exact Sequence}
\endhead

We recall the Huebschmann - Jones - Ratcliffe exact sequence, \cite{Hb, J1, Rc}: 
$$\eightpoint\aligned
1\longrightarrow\ &\text{\rm H}^1(Q, A)\overset {\pi^*}\to\longrightarrow
\text{\rm H}^1(G, A)\longrightarrow\text{\rm H}^1(N, A)^G\longrightarrow\\
&\longrightarrow\text{\rm H}^2(Q, A)\longrightarrow \tH^2(G, A)\longrightarrow
\La(G, N, A)\overset{\d }\to\longrightarrow\text{\rm H}^3(Q,
A)\overset {\pi^*}\to\longrightarrow\text{\rm H}^3(G, A),
\endaligned\tag2.1
$$
where either i)  $G$ is a separable locally compact group acting on a separable abelian
{\vna} $\sC$ with $A=\sU(\sC)$ and $N$ a Borel normal subgroup, or ii) $G$ is a discrete
group and $N$ a normal subgroup. We need the second case because the automorphism
group $\Aut(\sM)$ of a separable factor $\sM$ and the normal subgroup $\cntr(\sM)$ will be 
taken as the groups $G$ and $N$. If $\cntr(\sM)$ is not closed as in the case of an AFD
factor $\sM$, then the quotient group $\Aut(\sM)/\cntr(\sM)$ does not have a good
topological property beyond the discrete group structure. 

We are interested in the exactness at $\tH_\a^3(Q, A)$. In particular, we need an explicit
construction of $[\la, \mu]\in \La(G, N, A)$ such that $\d[\la, \mu]=[c]$ for those
$c\in \tZ_\a^3(Q, A)$ with $\pi^*(c)\in \tB_\a^3(G, A)$ in terms of a cochain $\mu\in
\tC^2(G, A)$ with $\pi^*(c)=\part_G^\a(\mu)$. In the situation where Polish topologies
are available on $G, N$ and $Q$, we assume or demand that all cocycles and cochains are
Borel and when it is appropriate equalities are considered modulo null sets relative to the
relevant measures. This kind of restrictions requires  us nailing down several objects
explicitly rather than relying on mere existence of the required objects through abstract
mechanism.

Given a cocycle $(\la, \mu)\in \tZ_\a(G, N, A)$, we have a $G$-equivariant  exact sequence:
$$\CD
E: 1@>>>A@>i_A>>E=A\times_\mu N@>j_E>{\underset{\se}\to \longleftarrow}>N@>>>1
\endCD
$$
along with a cross-section $\se$ such that
$$\aligned
&\se(m)\se(n)=\mu(m, n)\se(mn),\quad  m, n \in N;\\
&\a_g(\se(g\inv mg))=\la(m, g)\se(m), \text{or equivalently}\\
&\a_g(\se(m))=\la(gmg\inv, g) \se (gmg\inv), \quad g\in G.
\endaligned
$$
Choose a cross-section $\sp$ of $\pi$:
$$\CD
1@>>>N@>>>G@>\pi>\underset{\sp}\to\longleftarrow>Q@>>>1,
\endCD
$$
which generates the cocycle $\fnn\in \tZ(Q, N)$:
$$
\sp(p)\sp(q)=\fnn(p, q)\sp(pq), \quad p, q \in Q.
$$
Then we get the associated three cocycle $c_E\in \tZ_\a^3(Q, A)$ given by the
following:
$$\aligned
c_E&(p, q, r)=(\part_Q (\se\scirc \fnn))(p, q, r)\\
&=\a_{\sp(p)}(\se(\fnn(q, r)))\se(\fnn(p,
qr))(\se(\fnn(p, q))\se(pq, r)))\inv,\\
\endaligned
$$
which is expressed in terms of $(\la, \mu)$ and $\fnn$ directly:
$$\aligned
c^{\la, \mu}&(p, q, r)=\la(\sp(p)\fnn(q, r)\sp(p)\inv, \sp(p))\\
&\hskip 1in\times\mu(\sp(p)\fnn(q, r)\sp(p)\inv, \fnn(p, qr))\\
&\hskip 1in\times \mu(\fnn(p, q), \fnn(p, qr))\inv.
\endaligned\tag2.2
$$
This can be shown by a direct computation from the definition, which we leave to the reader.
We  denote the cohomology class $[c_E]\in \tH_\a^3(Q, A)$ of $c_E$ by $\d([\la,
\mu])$, which does not depends on the choices of the cross-sections $\sp$ and $\se$ but
only on the cohomology class of $(\la, \mu)$.

\proclaim{Lemma 2.1} The image $\d(\La(G, N, A))$ in
$\ththr(Q, A)$ consists of  precisely those $[\xi]\in \ththra(Q, A)$ such
that 
$\pi^*([\xi])=1$ in $\ththra(G, A)$. More precisely if a cochain  $\mu\in \tC_\a^2(G,
A)$ gives $\partg \mu=\pi^*(\xi)$, then 
$$\aligned
&\la(m, g)=\mu(g, g\inv m g)\mu(m, g)\inv, \quad m\in N, g\in G.
\endaligned\tag2.3
$$
together with the restriction $(i_N)_*(\mu)$ gives an element of $\tZ(G, N, A)$
such that $[\xi]=\d[\la, \mu]$ where $i_N$ is the embedding of $N$ into $G$,
i.e.,
$$\CD
1@>>>N@>i_N>>G@>\pi>>Q@>>>1.
\endCD
$$
The cochain $f\in \tcatw(Q, A)$ given by
$$
f(p, q)=\mu(\sp(p), \sp(q)) \mu(\fnn(p, q),
\sp(pq))\inv\in A,\tag2.4
$$
relates the original  cocycle $\xi\in  \tzthra(Q, A)$ and the new cocycle $c^{\la, \mu}$
in the following way{\rm:}
$$
\xi=(\part_Q f)c^{\la, \mu}.
$$
\endproclaim
\demo{Proof} First we construct a $G$-equivariant
exact sequence:
$$\CD
E: \quad 1@>>>A@>>>E@>>>N@>>>1
\endCD
$$
from the data $\pi^*(\xi)=\part_G(\mu)\in \tB_a^3(G, A)$ with $\mu\in \tC_\a^2(G, A)$.
 Let $B=
A^Q$ be the abelian group of all $A$-valued, (Borel if applicable), functions on $Q$
on which $Q$ acts by:
$$
(\a_p(b))(q)=\a_p(b(qp)),\quad p, q \in Q,\ b \in B=A^Q.
$$
Viewing $A$ as the subgroup of $B$ consisting of all constant
functions, we get an exact sequence:
$$\CD
1@>>>A@>i>>B@>j>>C@>>>1.
\endCD
$$
 By the cocycle identity,  we have
$$\aligned
&\xi(q, r, s)=\ainv{p}\bigg(\xi(pq, r, s)\xi(p, qr, s)\inv\xi(p, q,
rs)\xi(p, q, r)\inv\bigg)
\endaligned
$$
gives the following with $\eta(p, q, r)=\ainv{p}(\xi(p, q, r))$
viewed as an element of $\tC^2(Q, B)$ as a function of $p$:
$$\aligned
\xi(q, r, s)&=\a_q(\eta(pq, r, s))\eta(p, qr, s)\inv\eta(p, q, rs)
\eta(p, q, r)\inv\\
&=\text{Constant in}\ p.
\endaligned
$$
Hence $j_*(\part_Q \eta)=\part_Q(j_*(\eta))=1$, so that 
$j_*(\eta)\in \tZ^2(Q, C)$ and therefore we get an exact
sequence based on $j_*(\eta\inv)$:
$$\CD
1@>>>C@>>>D@>\sig>\underset{\fs_\sig}\to
\longleftarrow>Q@>>>1,
\endCD
$$
with $D=C\rtimes_{\a, j_*(\eta\inv)}Q$ and a cross-section $\fs_\sig: p\in
Q \mapsto (1, p)\in D$ such that
$$
\fs_\sig(p)\fs_\sig(q)=j(\eta(p, q)\inv)\fs_\sig(pq), \quad p, q \in Q.
$$

With $\mu\in \tC^2(G, A)$ such that $\pi^*(\xi)=\partg \mu$, we get
$$\aligned
\xi(\pi(g),\ &\pi(h), \pi(k))=\a_g(\mu(h, k))\mu(gh, k)\inv \mu(g,
hk)\mu(g, h)\inv\\
&=\a_g(\eta(p\pi(g), \pi(h), \pi(k)))\eta(p, \pi(gh), \pi(k))\inv \\
&\hskip .5in\times\eta(p, \pi(g), \pi(hk))\eta(p, \pi(g), \pi(h))\inv.
\endaligned
$$
Thus, $\z=i_*(\mu)\pi^*(\eta\inv)\in \tZ^2(G, B)$, which
allows  us to create an exact sequence:
$$\CD
1@>>>B@>>>F@>\tilde
\sig>\underset{\fs_{\tilde\sig}}\to\longleftarrow>G@>>>1
\endCD
$$
with $F=B\rtimes_{\a, \z}G$ and the cross-section $\fs_{\tilde\sig}$, given by $\fs_\tsig(g)=(1, g),
g\in G,$ such that
$$
\fs_{\tsig}(p)\fs_\tsig(q)=j(\zeta(p, q))\fs_\tsig(pq).
$$
Since 
$j_*(\z)=j_*(\eta\inv)$,  the next diagram is commutative:
$$\CD
1@>>>B@>>>F@>\tsig>>G@>>>1\\
@. @Vj VV@V(j\times \pi)VV@V\pi VV\\
1@>>>C@>>>D@>\sig>>Q@>>>1
\endCD
$$
With $E=\Ker(j\times \pi)$, we get the expanded commutative
diagram:
$$\CD
@.1@.1@.1\\
@. @VVV@VVV@VVV\\
1@>>>A@>>>E@>>>N@>>>1\\
@. @Vi VV@VVV@Vi_N VV\\
\endCD
$$
$$\CD
1@>>>B@>>>F@>\tilde\sig>
\underset{\fs_{\tilde\sig}}\to\longleftarrow>G@>>>1\\ 
@. @Vj VV@V(j\times
\pi)VV@V\pi V\uparrow\sp V\\
1@>>>C@>>>D@>\sig>\underset{\fs_\sig}\to\longleftarrow>Q@>>>1\\
@. @VVV@VVV@VVV\\ @.1@.1@.1\\
\endCD
$$
The construction of the above diagram came equipped with
cross-sections, $\fs_{\tilde\sig}$ and $\fs_\sig$ such that
$$
\fs_\sig\scirc \pi=(j\times \pi)\scirc \fs_{\tilde\sig}.
$$

We now look at the extension $E$, which is the kernel
$\Ker(j\times \pi)$. An element $(b, g)\in F$ belongs to $E$ if and
only if $j(b)=1$ and $\pi(g)=1$; if and only if $(b, g)\in
A\times N$. For $m, n\in N$ we have
$$
\fs_{\tilde \sig}(m)\fs_{\tilde\sig}(n)=\mu(m, n)
\eta(\pi(m), \pi(n))\inv\fs_{\tilde\sig}(mn)=\mu(m, n)
\fs_{\tilde \sig}(mn),
$$
so that we get
$$
E=A\times_{\mu}N.
$$
As  $E$ is a normal subgroup of $F$, each $\fs_{\tilde\sig}(g),
g\in G,$ normalizes $E$. Since the value of the two-cocycle
$\z$ belongs to $B$ each element of which commutes with
$A$ and $\fs_{\tilde\sig}(N)$, the restriction of
$\a_g=\Ad(\fs_{\tilde\sig}(g))$ to $E$ gives rise to an honest
action of $G$ which is consistent with the original action of
$G$ on $A$. Thus we obtain a $G$-equivariant exact sequence:
$$\CD
E:\quad
1@>>>A@>>>E@>{\tsig|_E}>\underset{\fs_\tsig|_N}\to\longleftarrow>N@>>>1.
\endCD
$$
Now  we compare the original $[\xi]$ and $[c_E]$ in
$\tH^3(Q, A)$. The cross-section $\fs_\tsig$ takes $N$ into $E$ so that its
restriction $\sts|_N$ is a cross-section for $\tsig|_E$.
 The associated three cocycle $c_E\in \tZ^3(Q, A)$ is obtained by:
$$\aligned
c_E&(p, q, r)=\part_Q (\fs_\tsig\scirc \fnn)(p, q, r).
\endaligned
$$
Consider the map $\fs=\fs_\tsig\scirc\sp: Q\mapsto E$ and compute
$$\aligned
( \part_Q &\fs)(p, q)=
\fs_{\tilde\sig}(\sp(p))\sts(\sp(q))\sts(\sp(pq))\inv\\
&=\mu(\sp(p), \sp(q)) \eta(p, q)\inv\sts(\sp(p)\sp(q))\sts(\sp(pq))\inv\\
&=\mu(\sp(p), \sp(q)) \eta(p, q)\inv\sts(\fnn(p, q)\sp(pq))
\sts(\sp(pq))\inv\\
&=\mu(\sp(p), \sp(q))\eta(p, q)\inv\mu(\fnn(p, q),
\sp(pq))\inv\eta(\pi(\fnn(p, q)), pq)\\
&\hskip .5in \times \sts(\fnn(p,
q))\sts(\sp(pq))\sts(\sp(pq))\inv\\
&=\mu(\sp(p), \sp(q)) \mu(\fnn(p, q),
\sp(pq))\inv\eta(p, q)\inv\sts(\fnn(p, q)).
\endaligned
$$
Thus with 
$$
f(p, q)=\mu(\sp(p), \sp(q)) \mu(\fnn(p, q),
\sp(pq))\inv\in A,\tag2.4
$$
we get
$$\aligned
1&=(\part_Q\partq \fs)(p, q, r)=\part_Q f(p, q, r)\part_Q\eta(p, q,
r)\inv\part_Q(\sts\scirc
\fnn)(p, q, r)\\
&=\part_Qf(p, q, r)\xi(p, q, r)\inv c_E(p, q, r),
\endaligned
$$
so that
$$
[\xi]=[c_E]\in \d(\La(G, N, A)).
$$
Now we compute the associated characteristic cocycle $(\la, \mu)\in \tZ(G, N, A)$:
$$\aligned
\la(m, g)&\sts(m)=\a_g(\sts(g\inv m g))=\sts(g)\sts(g\inv mg)\sts(g)\inv\\
&=\mu(g, g\inv mg)\eta(\pi(g), \pi(g\inv m  g))\inv\sts(mg)\sts(g)\inv\\
&=\mu(g, g\inv mg)\mu(m, g)\inv \eta(\pi(m), \pi(g))\sts(m)\sts(g)\sts(g)\inv\\
&=\mu(g,  g\inv mg)\mu(m, g)\inv\sts(m),
\endaligned
$$
which proves (2.3).

As we will need only the construction of a $G$-equivariant short exact sequence
from the cochain $\mu\in \tC_\a^2(G, A)$ with $\pi^*(\xi)=\partg \mu$, we
leave the proof for the converse to the reader. It is a direct computation. 
\QED
\enddemo

In the sequel, the group $G$ appears as the quotient group of another group $H$ by a
normal subgroup $M$, i.e., $G=H/M$. Let $\pig$ be the quotient map
$\pig: H\mapsto G$. Set $L=\pi\inv(N)$ and 
$$
\wtH=H\times \R, \quad \wtG=G\times \R, \quad \text{and }\quad
\wtQ=Q\times\R.\tag2.5
$$
Whenever an action $\a$ of $\wtH$ on a group $E$ is given, we denote the restriction of $\a$ to $\R$ 
by $\th$.  When an action $\a$ of the group $H$ is given and the cross-sections
$\sh\!: g\in G\mapsto \sh(g)\in H$ for $\pig,\ \fs\!: p\in Q\mapsto \fs(p)\in G$
for the quotient map $\pi\!: g\in G\mapsto \pi(g)=gN\in Q$  and $\dfs\!: p\in
Q\mapsto
\dfs(p)=\sh(\fs(p))\in H$ for the map $\dpi=\pi\scirc \pig$ are specified, we use the 
abbreviated notations:
$$
\a_g=\a_{\sh(g)},\quad g\in G;\quad \a_p=\a_{\dfs(p)}, \quad p\in Q,
$$
which satisfy:
$$\aligned
\a_g\scirc \a_h&=\a_{\fnm(g, h)}\scirc \a_{gh}, \quad g, h \in G;\\
\a_p\scirc \a_q&=\a_{\fnl(p, q)}\scirc \a_{pq}, \quad p, q \in Q,
\endaligned\tag2.6
$$
where
$$\aligned
\fnm(g, h)&=\sh(g)\sh(h)\sh(gh)\inv\in M, \quad g, h \in G;\\
\fnl(p, q)&=\dfs(p)\dfs(q)\dfs(pq)\inv\in L, \quad p, q \in Q.
\endaligned\tag2.7
$$

We examine the last half of the {\hjr}:
$$\CD
\thatw(\wtH, A)@>\res>>\La_\a(\wtH, L,
A)@>\d_{\scriptscriptstyle\text{HJR}}>>\thath(\wtQ, A)@>\inf>>\thath(\wtH, A).
\endCD
$$
First we show:
\proclaim{Lemma 2.2}  For each $\mu'\in\tzatw(\wtH, A)$, there is
an element $\mu\in \tzatw(\wtH, A)$ such that $\mu'$ and $\mu$ are cohomologous
and $\mu$ satisfies the condition{\rm:}
$$
\mu(\tilde h, \tilde k)=\mu_H(h, k)\a_h(d_\mu(s; k)),
\quad \tilde h=(h, s), \tilde k=(k, t)\in \wtH=H\times \R,\tag2.8
$$
where
$$\aligned
\mu_H\in \tzatw&(H, A),\quad d(\ \cdot\ ; h)\in \tzth1(\R, A);\\
\th_s(\mu_H(h, k))&\mu_H(h, k)^*=d(s; h) \a_h(d(s; k))
d(s; hk)^*;\\
\endaligned\tag2.9
$$
equivalently
$$
\part_\th \muh =\part_H d.\tag2.9$'$
$$

\endproclaim
\demo{Proof} We recall $\tH_\th^2(\R, A)=\{1\}$. So
we may and do assume that $\mu'(s, t)=1, s, t\in \R$.
Consider the group extension:
$$\CD
1@>>>A@>i>>F=A\times_{\mu'} \wtH@>j>>\wtH@>>>1.
\endCD
$$
The assumption on the restriction $\mu'|_{\R\times\R}$ allows us to find a one parameter
subgroup $\{u(s)\!: s\in \R\}$ of $  F $ with $j(u(s))=s, s \in \R$.
Choose a cross-section
$\fs_j'\!: h\in H\mapsto \fs_j'(h)\in F$ of the map $j$ such that
$$
\fs_j'(h)\fs_j'(k)=\mu'(h, k)\fs_j'(hk), \quad h, k\in H.
$$
Now set 
$$
\sj(h, s)=\fs_j'(h)u(s), \quad (h, s)\in \wtH.
$$
Now we compute the associated 2-cocycle $\mu\!:$
$$\aligned
\mu(h, s; k, t)&=\sj(h, s)\sj(k, t)\sj(hk, s+t)\inv\\
&=\fs_j'(h)u(s)\fs_j'(k)u(t)\{\fs_j'(hk)u(s+t)\}\inv\\
&=\fs_j'(h)\mu'(s; k)\fs_j'(k, s)u(t)\{\fs_j'(hk)u(s+t)\}\inv\\
&=\fs_j'(h)\mu'(s; k)\mu'(k; s)\inv
\fs_j'(k)u(s)u(t)\{\fs_j'(hk)u(s+t)\}\inv\\ &=\a_h\Big(\mu'(s;
k)\mu' (k; s) \inv\Big) \fs_j'(h)\fs_j'(k)\fs_j'(hk)\inv\\
&=\a_h\Big(\mu'(s; k)\mu' (k; a) \inv\Big)\mu'(h; k)
\endaligned
$$
for each $(h, s), (k, t)\in \wtH$. Setting
$$
\mu_H=\mu'|_H \quad \text{and}\quad d(s; h)=\mu'(s; h)
\mu'(h; s)^*
$$
we obtain the first formula and also
$$
d(s+t; h)=d(s; h)\th_s(d(t; h)), \quad s, t\in \R, h\in H.
$$

We next check the second identity which follows from the
cocycle identity for $\mu$ as seen below:
$$\aligned
1&=\a_{\tilde g}(\mu(\tilde h, \tilde k))\mu(\tilde g\tilde h,
\tilde k)^*\mu(\tilde g, \tilde h\tilde k)\mu(\tilde g, \tilde h)^*,\\
&\hskip 1.5in\tilde g=(g, s), \tilde h=(h, t), \tilde k=(k, u)\in \wtH,\\
&=\a_{\tilde g}\Big(\a_h(d_\mu(s; k))\muh(h,
k)\Big)\a_{gh}(d_\mu(s+t; k)^*)\\
&\hskip.5in\times\muh(gh, k)^*\a_g(d_\mu(s; hk))\muh(g,
hk)\a_g(d_\mu(s; h)^*)\muh(g, h)^*\\
&=\a_g\Big(\th_s(\a_h(d_\mu(  s  ; k))\muh(h, k))\a_h(d_\mu(s+t;
k)^*)\Big)
\muh(gh, k)^*\\
&\hskip .5in\times \a_g(d_\mu(s; hk))\muh(g, hk)\a_g(d_\mu(s;
h)^*)\muh(g, h)^*\\
&=\a_g\Big(\th_s(\a_h(d_\mu(  s ; k))\muh(h, k))\a_h(d_\mu(s; k)^*\th_s(d_\mu(t;
k)^*))\Big)
\muh(gh, k)^*\\
&\hskip .5in\times \a_g(d_\mu(s; hk))\muh(g, hk)\a_g(d_\mu(s;
h)^*)\muh(g, h)^*\\
&=\a_g\Big(\th_s(\muh(h, k))\a_h(d_\mu(s; k))^*\Big)
\muh(gh, k)^*\\
&\hskip .5in\times \a_g(d_\mu(s; hk))\muh(g, hk)\a_g(d_\mu(s;
h)^*)\muh(g, h)^*\\
&=\a_g\Big(\th_s(\muh(h, k))\muh(h, k)^*\a_h(d_\mu(s; k)^*)\Big)
\a_g(\muh(h, k))\muh(gh, k)^*\\
&\hskip .5in\times \a_g(d_\mu(s; hk))\muh(g, hk)\a_g(d_\mu(s;
h)^*)\muh(g, h)^*\\
&=\a_g\Big(\th_s(\muh(h, k))\muh(h, k)^*\a_h(d_\mu(s;
k)^*)d_\mu(s; hk)d_\mu(s; h)^*\Big);\\
&\hskip.2in\th_s(\muh(h, k))\muh(h, k)^*=d_\mu(s; h)\a_h(d(s; k))d(s;
hk)^*.
\endaligned
$$
This proves the lemma.
\QED
\enddemo
{\smc Definition 2.3.} A cocycle $\mu\in\tzatw(\wtH, A)$ of the form (2.8) will be called
{\it standard} and $d_\mu$ and $\mu_H$ in (1) will be called naturally the
$\R$-part and the $H$-part of the cocycle $\mu$.

\proclaim{Lemma 2.4} {\rm i)} If $\mu\in \tzatw(\wtH, A)$ is standard, then the $(\la_\mu,
\mu)=\res(\mu)\in \tZ_\a(\wtH, L, A)$ is given by the following{\rm:}
$$\aligned
\la_\mu(m; \tilde g)&=\a_g(d_\mu(s; g\inv mg))\muh(g, g\inv mg)\mu_H(m; g)^*, \\
&\hskip 1.5in
\tilde g=(g, s)\in \wtH, \ m\in L.
\endaligned\tag2.10
$$

{\rm ii)} If $(\la, \mu)\in \tZ_\a(\wtH, L, A)$, then $c=c_{\la,
\mu}=\d_{\scriptscriptstyle\text{\rm HJR}}(\la,
\mu)$ is given by{\rm:}
$$\aligned
c(\tp, \tq, \tr)&=\a_p\Big(\la(\fnl(q, r); s)\Big) \la(\dfs(p)\fnl(q, r))\dfs(p)\inv, \dfs(p))\\
&\hskip .5in\times\mu(\dfs(p)\fnl(q, r)\dfs(p)\inv, \fnl(p, qr)) \\
&\hskip .5in\times
\Big\{\mu(\fnl(p, q), \fnl(pq, r))\Big\}^*
\endaligned\tag2.2$'$
$$
for each triplet $\tp=(p, s), \tq=(q, t), \tr=(r, u)\in \wtQ$.

{\rm iii)} If $(\la, \mu)=(\la_\mu, \mu)=\res(\mu)$ with $\mu \in \tzatw(\wtQ, A)$ standard,
then the $3$-cocycle
$$
c=c_\mu=\dhjr(\la_\mu, \mu)
$$
is cobounded by $f\in\tcatw(\wtQ, A)$ given by{\rm:}
$$\aligned 
f(\tp, \tq)&=\mu(\dfs(\tp); \dfs(\tq))^*\mu(\fnl(p, q); \dfs(\tp \tq))\\
&=\mu(\dfs(p), s; \dfs(q), t)^*\mu(\fnl(p, q); \dfs(pq), s+t)\\
&=\a_p(d_\mu(s; \dfs(q))^*)\muh(\dfs(p), \dfs(q))^*\muh(\fnl(p, q); \dfs(pq))\in A,
\endaligned
$$
where $\dfs$ is a cross-section of the quotient homomorphism $\dot \pi\colon\ H\mapsto
Q=H/L$.
\endproclaim
\demo{Proof} i) The 2-cocyle $\mu\in \tzatw(\wtH, A)$ gives rise to the following commutative
diagram of exact sequences equipped with cross-sections $\sh$ and $\sj$:
$$\CD
1@>>>A@>>>
F=A\times_\mu \wtH@>\tj>\underset{\tilde
\fs_H}\to\longleftarrow >\wtH@>>>1\\
@.@|@A\bigcup AA@A\bigcup AA\\
1@>>>A@>>>E=A\times_\mu
L@>j>\underset{\sj}\to\longleftarrow >L@>>>1
\endCD
$$
The action $\a$ of $\wtH$ on $E$ is given by
$\a_{\tilde g}=\Ad(\sh(\tilde g))|_E, \tilde g=(g, s)\in H,$  viewing $E$ as a submodule of
$ F$, where the cross-section $\sh$ is given by
$$
\sh(g)=(1, g)\in  F=A\times _\mu \wtH,\quad g\in H.
$$ 
The action $\th$ of $\R$  on $E$ is given by:
$$
\th_s(a, m)=(\th_s(a)\la_\mu(m; s), m), \quad s\in \R, (a, m)\in
E=A\times_\mu L.
$$
Hence $\th_s(\sj(m))=\la_\mu(m; s)\sj(m), m\in L, s \in \R$. Now the cocycle
$$
\res(\mu)=(\la_\mu, \mu)\in \tZ(\wtH, L, A)
$$
is given by
$$\aligned
\la_\mu(m, \tilde g)&\sj(m)=\a_g(\sj(\tilde g\inv m \tilde g))=\sh(\tilde g)\sj(\tilde g\inv
m\tilde g)\sh(\tilde g)\inv\\
&=\sh(\tilde g)\sh(\tilde g\inv m\tilde g)\sh(\tilde g)\inv=\mu(\tilde g, \tilde g\inv
m\tilde g)\sh(m\tilde g)\sh(\tilde g)\inv\\
&=\mu(\tilde g, \tilde g\inv
m\tilde g)\mu(m, \tilde g)\inv \sh(m)\sh(\tilde g)\sh(\tilde g)\inv\\
&=\mu(\tilde g, \tilde g\inv
m\tilde g)\mu(m, \tilde g)\inv \sh(m);\\
&\la_\mu(m, \tilde g)=\mu(\tilde g, \tilde g\inv
m\tilde g)\mu(m, \tilde g)\inv
\endaligned 
$$
for $m \in L, \tilde g=(g, s)\in\wtH$. As $\mu$ is standard, we get further simplification:
$$\aligned
\la_\mu(m; g, s)&=\a_g(d(s; g\inv m g))\muh(g, g\inv mg)\muh(m, g)^*,\ (g, s)\in \wtH, m\in L.
\endaligned
$$

ii) Now suppose $(\la, \mu)\in \tZ_\a(\wtH, L, A)$. 
The associated $A$-valued 3-cocycle $c=c_{\la, \mu}$ is given by (2.2) and the formula
(2.2$'$) follows from (2.2) and the cocycle identity, see \cite{ST2: (1.7), page 411} for
$\la$:
$$\aligned
\la(\dfs(\tp)&\fnl(\tq, \tr)\dfs(\tp)\inv; \dfs(\tp))=\la(\dfs(p)\fnl(q, r)\dfs(p)\inv; \dfs(p), s)\\
&=\a_{\dfs(p)}(\la(\fnl(q, r)); s))\la(\dfs(p)\fnl(q, r)\dfs(p)\inv; \dfs(p))
\endaligned
$$
for $\tp=(p, s), \tq=(q, t), \tr=(r, u)\in \wtQ$ because $\dfs(\tp)=(\dfs(p), s)$.

iii) Now assume $(\la, \mu)=(\la_\mu, \mu)$ with $\mu\in \tzatw(\wtH, A)$ standard.
First we compute the associated cocycle $c=c_\mu$:
$$\allowdisplaybreaks\aligned
c_{\mu}&(\tp, \tq, \tr)=\a_{p}(\la_\mu(\fnl(q, r); s))\la_\mu(\dfs(p)\fnl(q,
r))\dfs(p)\inv, \dfs(p))\\
&\hskip .5in\times\mu(\dfs(p)\fnl(q, r)\dfs(p)\inv,
\fnl(p, qr))\\
&\hskip.5in\times 
\Big\{\mu(\fnl(p, q), \fnl(pq, r))\Big\}^*\\
&=\a_p\Big(\mu(s;  \fnl(q, r))\mu(\fnl(q, r); s)^*\Big)\\
&\hskip.5in\times
\mu(\dfs(p); \dfs(p)\inv \dfs(p)
\fnl(q, r)\dfs(p)\inv\dfs(p))\\
&\hskip .5in\times \mu(\dfs(p)\fnl(q, r)\dfs(p)\inv; \dfs(p))^*\\
&\hskip .5in\times\mu(\dfs(p)\fnl(q, r)\dfs(p)\inv;
\fnl(p, qr))\\
&\hskip.5in\times \Big\{\mu(\fnl(p, q); \fnl(pq, r))\Big\}^*\\
&=\a_p\Big(d_\mu(s;  \fnl(q, r))\Big)\muh(\dfs(p);  \fnl(q, r))\\
&\hskip.5in\times
\muh(\dfs(p)\fnl(q, r)\dfs(p)\inv ; \dfs(p))^*\\ 
&\hskip .5in\times\muh(\dfs(p)\fnl(q, r)\dfs(p)\inv;
\fnl(p, qr))\\
&\hskip.5in\times \Big\{\muh(\fnl(p, q); \fnl(pq, r))\Big\}^*.
\endaligned
$$
We now compute the coboundary of $f$:
$$\aligned
\part_\wtQ f(\tp&, \tq, \tr)=\a_\tp(f(\tq, \tr))f(\tp, \tq\tr)\{f(\tp, \tq)f(\tp\tq, \tr)\}^*\\
&=\a_{\tp}\Big(\mu(\dfs(\tq); \dfs(\tr))^*\mu(\fnl(q, r); \dfs(\tq\tr))\Big)\\
&\hskip.2in\times 
\mu(\dfs(\tp);
\dfs(\tq\tr))^*\mu(\fnl(p, qr); \dfs(pqr))\\
&\hskip .2in\times \Big\{\mu(\dfs(\tp); \dfs(\tq))^*\mu(\fnl(p, q); \dfs(\tp\tq))\\
&\hskip.2in\times
\mu(\dfs(\tp\tq); \dfs(\tr))^*\mu(\fnl(pq, r); \dfs(\tp\tq\tr))\Big\}^*\\
&=\mu(\dfs(\tp);\dfs(\tq)\dfs(\tr)) \mu(\dfs(\tp), \dfs(\tq))^*\mu(\dfs(\tp)\dfs(\tq); \dfs(\tr)  )^* \\
&\hskip .2in\times 
\mu(\dfs(\tp); \fnl(q, r)\dfs(\tq\tr))^*\\
\endaligned
$$
$$\aligned
&\hskip.2in\times
\Big\{\mu(\dfs(\tp); \fnl(q,
r))\mu(\dfs(\tp)\fnl(q, r); \dfs(\tq\tr))\Big\}\\
&\hskip .2in\times\mu(\dfs(\tp); \dfs(\tq\tr))^*\mu(\fnl(p, qr);
\dfs(\tp\tq\tr))\\
&\hskip .2in\times \Big\{\mu(\dfs(\tp); \dfs(\tq))^*\mu(\fnl(p, q); \dfs(\tp\tq))\\
&\hskip.2in\times
\mu(\dfs(pq); \dfs(r))^*\mu(\fnl(pq, r); \dfs(pqr))\Big\}^*\\
&=\mu(\dfs(\tp);\fnl(q, r)\dfs(\tq\tr)) \mu(\dfs(\tp); \dfs(\tq))^*\mu(\dfs(\tp)\dfs(\tq);
\dfs(\tr) )^* \\ &\hskip .2in\times \mu(\dfs(\tp); \fnl(q, r)\dfs(\tq\tr))^*\\
&\hskip.2in\times
\Big\{\mu(\dfs(\tp);
\fnl(q, r))\mu(\dfs(\tp)\fnl(q, r); \dfs(\tq\tr))\Big\}\\
&\hskip .2in\times\mu(\dfs(\tp); \dfs(\tq\tr))^*\mu(\fnl(p, qr); \dfs(\tp\tq\tr))\\
&\hskip .2in\times \Big\{\mu(\dfs(\tp); \dfs(\tq))\mu(\fnl(p, q); \dfs(\tp\tq))^*\\
&\hskip.2in\times
\mu(\dfs(\tp\tq); \dfs(\tr))\mu(\fnl(pq, r); \dfs(\tp\tq\tr))^*\Big\}\\
&= \mu(\dfs(\tp)\dfs(\tq); \dfs(\tr))^*\Big\{\mu(\dfs(\tp); \fnl(q,
r))\mu(\dfs(\tp)\fnl(q, r); \dfs(\tq\tr))\Big\}\\
&\hskip .2in\times\mu(\dfs(\tp); \dfs(\tq\tr))^*\mu(\fnl(p, qr);
\dfs(\tp\tq\tr))\\
&\hskip .2in\times
\Big\{\mu(\fnl(p, q); \dfs(\tp\tq))^*
\mu(\dfs(\tp\tq); \dfs(\tr))\mu(\fnl(pq, r); \dfs(\tp\tq\tr))^*\Big\}.
\endaligned
$$
We compute some terms below:
$$\aligned
\mu(\dfs(\tp)&\fnl(q, r); \dfs(\tq\tr))=\mu(\dfs(\tp)\fnl(q, r)\dfs(\tp)\inv\dfs(\tp); \dfs(\tq\tr))\\
&=\mu(\dfs(\tp)\fnl(q, r)\dfs(\tp)\inv; \dfs(\tp))^*\\
&\hskip.5in\times
\mu(\dfs(\tp)\fnl(q, r)\dfs(\tp)\inv; \dfs(\tp)\dfs(\tq\tr))\mu(\dfs(\tp); \dfs(\tq\tr))\\
&=\mu(\dfs(\tp)\fnl(q, r)\dfs(\tp)\inv; \dfs(\tp))^*
\mu(\dfs(\tp); \dfs(\tq\tr))\\
&\hskip.5in\times
\mu(\dfs(\tp)\fnl(q, r)\dfs(\tp)\inv; \fnl(p, qr)\dfs(\tp\tq\tr))\\
&=\mu(\dfs(\tp)\fnl(q, r)\dfs(\tp)\inv; \dfs(\tp))^*
\mu(\dfs(\tp); \dfs(\tq\tr))\mu(\fnl(p, qr); \dfs(\tp\tq\tr))^*\\
&\hskip.5in\times \mu(\dfs(\tp)\fnl(q, r)\dfs(\tp)\inv; \fnl(p, qr))\\
&\hskip.5in\times\mu(\dfs(\tp)\fnl(q, r)\dfs(\tp)\inv 
\fnl(q, r); \dfs(\tp\tq\tr));\\
\mu(\dfs(\tp)&\dfs(\tq); \dfs(\tr))=\mu(\fnl(p, q)\dfs(\tp\tq); \dfs(\tr))\\
\endaligned
$$
$$\aligned&=\mu(\fnl(p, q);\dfs(\tp\tq))^*\mu(\fnl(p, q); \dfs(\tp\tq)\dfs(\tr))\mu(\dfs(\tp\tq); \dfs(\tr))\\
&=\mu(\fnl(p, q);\dfs(\tp\tq))^*\mu(\dfs(\tp\tq); \dfs(\tr))\mu(\fnl(p, q); \fnl(pq,
r)\dfs(\tp\tq\tr))\\ &=\mu(\fnl(p, q);\dfs(\tp\tq))^*\mu(\dfs(\tp\tq); \dfs(\tr))\mu(\fnl(pq, r);
\dfs(\tp\tq\tr))^*\\ &\hskip .5in\times\mu(\fnl(p, q); \fnl(pq, r))\mu(\fnl(p, q)\fnl(pq, r);
\dfs(\tp\tq\tr)).
\endaligned
$$
We then substitute the above expression in the original calculation:
$$\allowdisplaybreaks\aligned
\part_\wtQ f&(\tp, \tq, \tr)=  \mu(\dfs(\tp)\dfs(\tq), \dfs(\tr))^*\mu(\dfs(\tp), \fnl(q, r))\\
&\hskip.2in\times
\mu(\dfs(\tp)\fnl(q, r), \dfs(\tq\tr))\mu(\dfs(\tp), \dfs(\tq\tr))^*\mu(\fnl(p, qr),
\dfs(\tp\tq\tr))\\
&\hskip .2in\times
\mu(\fnl(p, q), \dfs(\tp\tq))^*
\mu(\dfs(\tp\tq), \dfs(\tr))
\mu(\fnl(pq, r), \dfs(\tp\tq\tr))^*\\
&=\mu(\fnl(p, q),\dfs(\tp\tq)) \mu(\dfs(\tp\tq), \dfs(\tr))^*\\
&\hskip.2in\times
\mu(\fnl(pq, r), \dfs(\tp\tq\tr)) \mu(\fnl(p, q), \fnl(pq, r))^*\\
&\hskip.2in\times
\mu(\fnl(p, q)\fnl(pq, r),
\dfs(\tp\tq\tr))^*\\
&\hskip.2in\times\mu(\dfs(\tp), \fnl(q, r))
\mu(\dfs(\tp)\fnl(q, r)\dfs(\tp)^{  -1 }, \dfs(\tp))\\
&\hskip.2in\times
\mu(\dfs(\tp), \dfs(\tq\tr))\mu(\fnl(p, qr), \dfs(\tp\tq\tr))^*\\
&\hskip.2in\times
\mu(\dfs(\tp)\fnl(q, r)\dfs(\tp)\inv, \fnl(p, qr))\\
&\hskip.2in\times \mu(\dfs(\tp)\fnl(q, r)\dfs(\tp)\inv 
\fnl(q, r), \dfs(\tp\tq\tr))\\
&\hskip .2in\times\mu(\dfs(\tp), \dfs(\tq\tr))^*\mu(\fnl(p, qr),
\dfs(\tp\tq\tr))\\
&\hskip .2in\times \mu(\fnl(p, q), \dfs(\tp\tq))^*
\mu(\dfs(\tp\tq), \dfs(\tr))
\mu(\fnl(pq, r), \dfs(\tp\tq\tr))^*\\
&=  \mu(\fnl(p, q), \fnl(pq, r))^*\mu(\dfs(\tp), \fnl(q, r))\\
&\hskip.2in\times
\mu(\dfs(\tp)\fnl(q, r)\dfs(\tp)\inv, \dfs(\tp))\\
&\hskip.2in\times \mu(\dfs(\tp)\fnl(q, r)\dfs(\tp)\inv, \fnl(p, qr))\\
&=\a_p\Big(d_\mu(s;  \fnl(q, r))\Big)\muh(\dfs(p);  \fnl(q, r))\\
&\hskip.2in\times
\muh(\dfs(p)\fnl(q, r)\dfs(p)\inv ; \dfs(p))^*\\ 
&\hskip .2in\times\muh(\dfs(p)\fnl(q, r)\dfs(p)\inv;
\fnl(p, qr))\\
&\hskip.2in\times 
\Big\{\muh(\fnl(p, q); \fnl(pq, r))\Big\}^*\\
&=c_\mu(\tp, \tq, \tr) \\
\endaligned
$$
for each triplet $\tp=(p, s), \tq, \tr \in \wtQ$. This completes the proof.
\QED
\enddemo

\proclaim{Lemma 2.5} {\rm i)} Every cohomology class $[c]\in \ththra(\wtQ, A)$ can
be represented by a cocycle $c$ of the form{\rm:}
$$\aligned
c(\tp&, \tq, \tr)=\a_p(d_c(s; q, r))c_Q(p, q, r),\\
&\hskip .5in\tp=(p, s), \tq=(q, t), \tr=(r, u)\in \wtQ,
\endaligned\tag2.11
$$
where $c_Q\in \tzthra(Q, A)$ and $d_c(\cdot, q, r)\in \tzth1$.

{\rm ii)} Given a function $d: \R\times Q^2\mapsto A$ and $c_Q\in \tzath(Q, A)$, the
function $c$ given by{\rm:}
$$
c(\tp, \tq, \tr)=\a_p(d(s; q, r))c_Q(p, q, r)
$$
is an element of $\tzath(\wtQ, A)$ if and only if
\roster 
\item"a)" for each fixed $q, r\in Q$, $d(\cdot, q, r)$ is an $\R$-cocycle, i.e., 
$$
d(s+t, q, r)=d(s, q, r)\th_s(d(t, q, r)),\quad s, t \in \R, q, r \in Q;
$$
\item"b)"
$c_Q$ and $d$ are linked by the
following formula{\rm:}
$$\aligned
\th_s(c_Q(p, q, r))&c_Q(p, q, r)^*\\
&=\a_p(d(s; q, r))d(s; p, qr)\{d(s; p, q)d(s; pq, r)\}^*
\endaligned\tag2.12
$$
for each $s\in \R, p, q, r\in Q$, i.e.,
$$
\partq d=\part_\th c_Q.\tag2.12$'$
$$
\endroster

{\rm iii)} For a cocycle $c\in \tZ_\a^3(\wtQ, A)$ of the form {\rm (2.11)} the following
are  equivalent{\rm:}
\roster
\item"a)" There exists $a\in \tC_\a^2(Q, A)$ such that 
$$\aligned
c=\part_\wtQ a;
\endaligned
$$
\item"b)" There exists $a\in \tC_\a^2(Q, A)$ such that
$$
d_c(s; q, r)=\th_s(a(q, r))a(q, r)^*, q, r \in Q, s \in \R;\quad c_Q=\part_Q a.
$$
\endroster
\endproclaim
\demo{Proof} The assertion (i) follows from the fact that the additive real line $\R$ has
trivial second and third cohomologies. Every 3-cocyle we encounter in this paper will be
of this form without perturbation anyway. So we omit the proof. 

ii) This follows directly from the cocycle identity for $c$. We omit the detail. 

iii) This equivalence again follows from a direct easy computation.
\QED
\enddemo

{\smc Definition 2.6.} A cocycle $c\in \tzath(\wtQ, \R)$ of the form (2.11) will be called
{\it standard}. We will concentrate on the subgroup $\tzasth(\wtQ, A)$ of all standard
cocyles in $\tzath(\wtQ, A)$. The index ``s" stands for ``standard".  We then set
$$
\thasth(\wtQ, A)=\tzasth(\wtQ, A)/\part_{\wtQ}(\tcatw(Q, A)).
$$
The coboundary group $\tbasth(\wtQ, A)=\part_{\wtQ}(\tcatw(Q, A))$ is a subgroup of 
the usual third coboundary group $\tbath(\wtQ, A)=\part_\wtQ(\tcatw(\wtQ, A))$, so that
we have a natural surjective homomorphism: 
$$\CD
\thasth(\wtQ, A)@>>> \thath(\wtQ, A)
\endCD
$$

The fixed cross-section
$\fs\!: Q\mapsto G$ allows us to consider the fiber product
$$
\thasth(\wtQ, A)*_\fs\Hom_G(N, \thth1)
$$
consisting of those pairs $([c], \nu)\in \thasth(\wtQ, A)\times \Hom_G(N, \thth1)$
such that
$$\aligned
[d_c(\cdot, q, r)]= \nu(\fnn(q, r))\quad \text{in }\ \thth1, \quad q, r \in Q.
\endaligned
$$
The group $\thasth(\wtQ, A)*_\fs\Hom_G(N, \thth1)$ will be denoted by $\thasout(G\times\R,
N, A)$ for short. The suffix ``$\fs$" is placed to indicate that this fiber product depends
heavily on the cocyle $\fnn$ hence on the cross-section $\fs$. As mentioned earlier, the
invariant for outer actions of $G$ must respect the cross-section $\fs$ because
a change in the cross-section results an alteration on the outer conjugacy class
from the original outer conjugacy class. Before stating the main theorem of
the section, we still need some preparation.

\proclaim{Theorem 2.7} Suppose that $\{\sC, \R, \th\}$ is an ergodic flow and a
homomorphism $\a${\rm:} $g\in H\mapsto \a_g\in \Aut_\th(\sC)$, the group of
automorphisms of $\sC$ commuting with $\th$. Assume the following{\rm:}
\roster
\item"i)" a pair of normal subgroup $M\subset L\subset H$ is given{\rm;}
\item"ii)" the subgroup $L$, hence $M$ as well, acts trivially on $\sC$, i.e., $L\subset
\Ker(\a)${\rm;}
\item"iii)" with $G=H/M, N=L/M$ and $Q=H/L$, let $\pig\,\coron\ H\mapsto G$,
$\pi\,\coron\ G\mapsto Q$ and $\dpi=\pi\scirc \pig\,\coron\ H\mapsto Q$ be the quotient maps
such that
$$
\Ker(\pig)=M, \quad \Ker(\pi)=N\quad \text{and}\quad \Ker(\dpi)=L;
$$
\item"iv)" Fix a cross-section $\fs\!:\ Q\mapsto G$ of the map $\pi$ and
choose  cross-sections $\sh\coron\ G\mapsto H$ and $\dfs\!:\ Q\mapsto
H$ in such a way that
$$
\dfs=\sh\scirc \fs.
$$
\endroster
Set $\wtH=H\times \R, \wtG=G\times \R$ and $\wtQ=Q\times \R$.
Let A denote the unitary group $\sU(\sC)$ of $\sC.$ $  $ Under the above setting, there is a
natural exact sequence which sits next to the {\rm Huebshmann-Jones-Ratcliffe} exact
sequence{\rm:}
$$\CD
\thtw(H, \T)@=\thtw(H, \T)\\
@V\Res VV@V\res VV\\
\La(\wtH, L, M, A)@>>>\La_\a(H, M, \T)\\
@V\d VV@V\dhjr VV\\
\thasout(G\times\R, N, A)@>\part>>\ththr(G, \T)\\
@V\Inf VV@V\inf VV\\
\ththr(H, \T)@=\ththr(H, \T)
\endCD\tag2.13
$$
\endproclaim
We need some preparation.

\proclaim{Lemma 2.8} To each characteristic cocycle $(\la, \mu)\in \tZ_\a(\wtH, L, A)$, there
corresponds uniquely an $\wtH$-equivariant  exact square{\rm:}
$$\eightpoint\CD@.1@.1@.1\\
@.@VVV@VVV@VVV\\
1@>>>\T@>>>A@>\part_\th>>\tB@>>>1\\
@.@VVV@V i VV@VVV\\
1@>>>U=E^\th@>>>E@>\tpart_\th>>\tZ@>>>1\\
@.@VVV@V j VV@V\piz VV\\
1@>>>K@>>>L@>\dot\part>>\tH@>>>1\\
@.@VVV@VVV@VVV\\
@.1@.1@.1
\endCD\tag2.14
$$
with $E=A\times_\mu L$.
\endproclaim
\demo{Proof} The cocycle $(\la, \mu)$ gives an $\wtH$-equivariant exact sequence:
$$\CD
E: 1@>>>A@>i>>E=A\times_\mu L@>j>>L@>>>1.
\endCD
$$
With $U=E^\th$, the fixed point subgroup of $E$ under the action $\th$ of
$\R$, we set
$$\aligned
&\quad\tB=A/\T\cong \tbth1(\R, A);\quad \tZ=E/U;\\
& K= K(E)=j(U)\cong U/\T;\quad\tH=\tZ/\tB.
\endaligned
$$
As the real line $\R$ does not act on the group $H$, we have 
$$\aligned
j(\th_s(x)x\inv)&=j(x)j(x)\inv=1, \quad x\in E,\ s\in \R;\\
&\th_s(x)x\inv =(\part_\th x)_s\in A,
\endaligned
$$
and $a:  s\in \R\mapsto a_s=(\part_\th x)_s\in A$ is a cocycle, a member of $\tzth1(\R, A)$. Thus
$\tZ\subset \tzth1(\R, A)$ and naturally $\tH\subset \thth1(\R, A)$. The map $\part_\th$ can be viewed 
either as the quotient map: $E\mapsto \tZ$ or the coboundary map described above.
Now it is clear that these groups
$\T, A, \cdots, \tH$ form a commutative exact sequare of (2.14) on which $\wtH$ acts.
\QED
\enddemo
We will denote the subgroup $K$ of $L$ in (2.14) by $K(\la, \mu)$ or $K(\chi)$ to indicate the dependence of
$K$ on the cocycle $(\la, \mu)\in \tZ_\a(\wtH, L, A)$ or the charactieristic invariant
$\chi=[\la, \mu]\in \La_\a(\wtH, L, A)$. We then define the subgroups: 
$$\aligned
\tZ_\a(\wtH, L, M, A)&=\{(\la, \mu)\in \tZ_\a(\wtH, L, A): K(\la, \mu)\supset M\};\\
\La_\a(\wtH, L, M, A)&=\{\chi\in \La_\a(\wtH, L, A):\ K(\chi)\supset M\}.
\endaligned\tag2.15
$$
 A cocycle $(\la, \mu)\in \tZ(\wtH, L, A)$ belongs to the subgroup
 $\tZ_\a(\wtH, L, M, A)$ if and only if the cocycle $(\la, \mu)$ satisfies the
conditions:
$$
(\la|_{M\times \wtH}, \mu|_{M\times M})\in \tZ(\wtH, M, \T); \quad \la(m; s)=1,\quad s\in \R, m \in
M.
$$

 Let $\pr_H: \wtH\mapsto H$ be the projection map from $\wtH=H\times
\R$ to the $H$-component and $i_{A, \T}: \T\mapsto A$ be the canonical embedding of
$\T$ to $A$. Finally, let $i_{L, M}: M\mapsto L$ be the embedding of $M$ into $L$. Then we have
naturally:
$$\CD
@.\La_\a(\wtH, L, A)@>i^*_{L, M}>>\La(\wtH, M, A)\\
@.@.@|\\
\La(H, M, \T)@>\pr_H^*>>\La(\wtH, M, \T)@>(i_{A, \T})_*>>\La(\wtH, M, A)
\endCD
$$
In terms of these maps, we can restate the subgroup $\La_\a(\wtH, L, M, A)$ in the
following way:
$$
\La_\a(\wtH, L, M, A)=(i_{L, M}^*)\inv((i_{A, \T})_*\scirc \pr_H^*(\La(H, M, \T)).
$$

The above maps also generates the following chain:
$$\CD
\thtw(H, \T)@>\pr_H^*>>\thtw(\wtH, \T)@>(i_{A, \T})_*>>\thatw(\wtH, A)
@>\res>>\La(\wtH, L, A)
\endCD
$$
and the range of the composed map 
$$
\Res=\res\scirc(i_{A, \T})_*\scirc\pr_H^*:\thtw(H, \T)\mapsto \La_\a(\wtH, L, A)
$$ 
is contained in the group $\La_\a(\wtH, L, M, A)$ defined above, which
generates the maps:
$$\CD
\thtw(H, \T)@>\Res>>\La_\a(\wtH, L, M, A).
\endCD
$$
Coming back to the orignial situation that $H=\Aut(\sM), M=\Int(\sM)$ and $L=\cntr(\sM)$, we know
that each element of $\Res(\thtw(H, \T))$ gives a perturbation of the action of
$\Aut(\sM)$ on $\sM$ differ by $\Int(\sM)$. Hence we must be concerned with the
quotient group 
$$
\La_\a(\wtH, L, M, A)/\Res(\thtw(H, \T)).
$$

The map $\dhjr=\d$ in the the Huebschmann - Jones - Ratcliffe exact sequence: 
$$\eightpoint\aligned
1\longrightarrow\ &\text{\rm H}^1(\wtQ, A)\overset {\pi^*}\to\longrightarrow
\text{\rm H}^1(\wtH, A)\longrightarrow\text{\rm H}^1(L, A)^\wtH\longrightarrow\\
&\longrightarrow\text{\rm H}^2(\wtQ, A)\longrightarrow \tH^2(\wtH, A)\longrightarrow
\La(\wtH, L, A)\overset{\d }\to\longrightarrow\text{\rm H}^3(\wtQ,
A)\overset {\pi^*}\to\longrightarrow\text{\rm H}^3(\wtH, A),
\endaligned
$$
 to be abbreviated  the {\bf HJR-exact sequence},  \cite{Hb, J1, Rc},
gives a natural map $\dhjr\coron\ \La_\a(\wtH, L, M, A)\mapsto \thath(\wtQ, A)$. 
The map $\dhjr$ will  be called the {\bf HJR} map and the modified
HJR map $\d$ relevant to our discussion  will be constructed along with the other
two maps:
$$\aligned
\part:\ &\thasout(G\times\R, N, A)\longrightarrow \ththr(G, \T);\\
\Inf:\ &\thasout(G\times\R, N, A)\longrightarrow \ththr(H, \T).
\endaligned
$$

\subhead\nofrills{Construction of the modified HJR-map $\d$:}
\endsubhead\
First we fix a cocycle $(\la,  \mu)\in \La_\a(\wtH, L, M, A)$ and consider the
corresponding  crossed extension $E$:
$$\CD
1@>>>A@>i>>E@>j>>L@>>>1.
\endCD
$$
As $M\subset K(\la,  \mu)$, with $V=M\times_\mu \T$ and $F=E/V$
we have an $\wtH$-equivariant exact square: 
$$\CD@.1@.1@.1\\
@.@VVV@VVV@VVV\\
1@>>>\T@>>>A@>\part_\th>>\tB@>>>1\\
@.@VVV@V i VV@VVV\\
1@>>>V@>>>E@>>>F@>>>1\\
@.@VVV@V j V\big\uparrow\sj  V@V\pi_N VV\\
1@>>>M@>>>L@>\pig|_L>\underset{\fsg}\to\longleftarrow 
>N@>>>1\\
@.@VVV@VVV@VVV\\
@.1@.1@.1
\endCD
$$
As $M\subset K(\la, \mu)$, we get a $G$-equivariant
homomorphism $\nu_\chi: N\mapsto \tH\subset \thth1(\R,
A)$, where $G=H/M$, i.e., $\nu_\chi\in \Hom_G(N, \thth1(\R,
A)).$ 

\proclaim{Lemma 2.9} Fix $(\la, \mu)\in \tZ_\a(\wtH, L, M, A)$.

 {\rm i)} For  the cross-section $\sj\!: m\in L\mapsto (1, m)\in E=A\times_\mu L$  of
the map $j\!: E\mapsto L$ associated with the cocycle
$(\la, \mu)$, the cocycle $c=c^{\la, \mu}$ given by
{\rm(2.2$'$)} is standard with
$$\aligned
d_c(s; q, r)&=\la(\fnl(q, r); s), \quad q, r \in Q, s\in \R;\\
c_Q(p, q, r)&=\la(\dfs(p)\fnl(q, r)\dfs(p)\inv; \dfs(p))\mu(\dfs(p)\fnl(q, r)\dfs(p)\inv,
\fnl(p, qr))\\
&\hskip 1in\times
\{\mu(\fnl(p, q), \fnl(pq, r)\}^*, \quad p, q, r\in Q.
\endaligned
$$

{\rm ii)} 
$$\aligned
([c_\chi],  \nu_\chi)\in \thasth(\wtQ, A)*_{\text
s}\Hom_G(N, \thth1(\R, A));
\endaligned
$$
\endproclaim
\demo{Proof} i) This is obvious from the formula (2.2$'$).

ii) 
The cross-section $\sh$ gives an $M$-valued 2-cocycle:
$$
\fnm(g, h)=\sh(g)\sh(h)\sh(gh)\inv\in M, \quad g, h \in G,
$$
which allows us  to relate $\sh(\fnn(q, r))$ and $\fnl(q, r)$:
$$\aligned
\pig(\fnl(q, r))&=\pig\Big(\dot
\fs(q)\dot \fs(r)\dot \fs(qr)\inv\Big)\\
&=\fs(q)\fs(r)\fs(qr)\inv=\fnn(q, r).
\endaligned
$$
Hence $\fnl(q, r)\equiv \sh(\fnn(q, r))\ \mod\ M$. As for each $m\in M, \ell\in L$ and 
$s\in \R$ we have
$$\aligned
\la(m\ell; s)&=\th_s(\mu(m; \ell)^*)\mu(m; \ell)\la(m; s)\la(\ell; s)\\
&=\th_s(\mu(m; \ell)^*)\mu(m; \ell)\la(\ell; s),
\endaligned
$$
we get
$$
[\la(m\ell;\ \cdot)]=[\la(\ell;\ \cdot)]\quad \text{in}\ \thth1(\R, A)\quad
\text{for every}\ m\in M, \ell \in L.
$$
Thus $[\la(\fnl(q, r); \ \cdot)]=\nu_\chi(\fnn(q, r))\in \thth1(\R, A)$, which precisely
means that $([c^{\la, \mu}], \nu)\in \thasout(G\times\R, N, A)$.
\QED
\enddemo

Thus we obtain an element: 
$$
\d(\chi)=([c^{\la, \mu}], \nu_\chi)\in \thasout(G\times\R, N, A),
$$
and therefore the map 
$$
\d: \La_\a(\wtH, L, M, A)\mapsto \thasout(G\times\R, N, A).
$$
We will call $\d$ the {\bf modified HJR-map}. To distinguish this modified HJR map
from the original HJR map, we denote the original one by $\dhjr$ and the modified one
simply by $\d$. 
Note that the map $\d$ does not depend on the choice of the
section $\sh\,\coron G\longmapsto H,$ but depends on the choice of 
the section $\fs\,\coron Q\longmapsto G. $
We  now begin the proof of Theorem 2.7.
\subhead\nofrills{$\pmb{\Ker(\d)=\text{Im}(\Res)}$:} 
\endsubhead\
First assume that
$$
\mu\in \tztw(H, \T)\quad\text{and}\quad\chi=\Res([\mu]) \in  \La(\wtH, L, M, A).
$$
Then we have a
commutative diagram of exact sequences:
$$\CD
1@>>>\T@>>>F=\T\times_\mu
H@>j_H>\underset{\fs\!_F}\to\longleftarrow >H@>>>1\\ 
@.@VVV@V\bigcap VV@|\\
1@>>>A@>>>\widetilde
F=A\times_\mu H@>\tj>\underset{\tilde
\fs\!_{ F}}\to\longleftarrow >H@>>>1\\
@.@|@A\bigcup AA@A\bigcup AA\\
1@>>>A@>>>E=A\times_\mu
L@>j>\underset{\sj}\to\longleftarrow >L@>>>1
\endCD
$$
The action $\a$ of $H$ on $E$ is given by
$\a_h=\Ad(\sh(h))|_E, h\in H,$  viewing $E$ as a submodule of
$\tilde F$, where the cross-section $\sh$ is given by
$$
\sh(h)=(1, h)\in \tilde F=A\times _\mu H,\quad h\in H.
$$ 
The action $\th$ of $\R$  on $E$ is given by:
$$
\th_s(a, m)=(\th_s(a), m), \quad s\in \R, (a, m)\in
E=A\times_\mu L.
$$
Hence $\th_s(\sj(m))=\sj(m), m\in L, s \in \R$. Now 
$\res(\mu)=(\la_\mu, \mu)\in \tZ(\wtH, L, A)$ is given by (2.10). As $\mu$ takes
values in $\T$, we have $\mu=\muh$, i.e., $d_\mu=1$. Consequently, $\la_\mu(m; s)=1,
m\in L, s\in
\R$ which entails
$$
\la_\mu(\fnl(p, q); s)=1, p, q \in Q=H/L, s \in \R.
$$
By Lemma 2.4.(iii), the associated 3-cocyle $c_\mu=c_{\la_\mu,\mu}\in
 \tzath(\wtQ, A)$ is co-bounded by $f$ of the form:
$$
f(p, q)=\mu(\dfs(p), \dfs(q))^*\mu(\fnl(p, q), \dfs(pq))\in \T.
$$
This shows that ${\text {\rm Im (Res)}}\subset \Ker(\d)$.

We are now moving to show the reversed inclusion: $\text {\rm Im (Res)}\supset \Ker(\d)$. We
first compare the original {\hjr} and our modified HJR sequence. To this end, we recall
 that the cohomology group $\thasth(\wtQ, A)$ is obtained as the quotient group
of a subgroup $\tzasth(\wtQ, A)$ of $\tzath(\wtQ, A)$ by a subgroup $\tbasth(\wtQ, A)$ of
$\tbath(\wtQ, A)$. Thus we have a natural map: $\thasth(\wtQ, A)\mapsto \thath(\wtQ, A)$.
Consequently, the above {\hjr} applied to our context yields the following commutative diagram:
$$\eightpoint\CD
\tH_\a^2(\wtH, A)@>\res>>\La_\a(\wtH, L, A)@>\dhjr>>\tH_\a^3(\wtQ, A)@>\inf>>
\tH_\a^3(\wtH, A)\\
@AA (i_{A, \T})_* A@AAA@AAA\\
\thtw(H, \T)@>\Res>>\La_\a(\wtH, L, M, A)@>\d>>\thasth(\wtQ, A)
\endCD
$$

Suppose $\chi=[\la, \mu]\in \Ker(\d)\subset \La_\a(\wtH, L, M, A)$. Then 
$$
1=\d(\chi)=([c_\chi], \nu_\chi)\quad\text{ in }\quad
\thasth(\wtQ, A)*_\fs \Hom_G(N, \thth1).
$$
The above assumption also means  $\dhjr(\chi)=1$. Hence the {\hjr} guarantees that the
2-cocyle $\mu$ on $L$ can be extended to $\wtH$ as an $A$-valued 2-cocyle over
$\wtH$ which we denote by $\mu$ again so that $\la=\la_\mu$. To proceed further, we
need the following:

\proclaim{Lemma 2.10} If a {\rm 2}-cocycle $\mu\in \tzatw(\wtH, A)$ is standard,
and
if
$$
\chi=\res([\mu])= [\la_\mu, \mu]\in \La_\a(\wtH, L, A)
$$
generates trivial $\nu_\chi=1$
of
$\Hom_G(N, \thth1(\R, A))$, then there exists a standard $\tmu\in \tzatw(\wtH, A)$ such that
\roster
\item"i)" $\mu\equiv \tmu\quad \mod\ \tbatw(\wtH, A);$
\item"ii)" $\la_\tmu(m; s)=1,\quad m\in L, s\in \R$, i.e., $d_\tmu(s; m)=1$.
\endroster
\endproclaim
\demo{Proof} Let $\piz$ be the quotient map: $c\in\tZ\mapsto [c]\in \tH=\tZ/\tB$.
The condition $\nu_\chi(n)=1, n \in N,$ implies that
$$\aligned
\pi_\tZ(\partth(\sj(m)))&=\dpartth(m)=\nu_\chi(\pig(m))=1 \in \thth1,\quad m\in L;\\
&\partth(\sj(m))\in \tbth1,
\endaligned
$$
so that for each $m\in L$ there exists $a(m)\in A$ such that
$$\aligned
\la(m; s)=\th_s(\sj(m))\sj(m)\inv=\partth(\sj(m))_s=\th_s(a(m))a(m)^*.
\endaligned
$$
Extending the function $a: L\mapsto A$ to the entire $\wtH$ in
such a way that
$$
a(g, s)=a(g), \quad (g, s)\in \wtH=H\times \R,
$$
we define a new 2-cocycle:
$$
\tmu(\tilde g, \tilde h)=a(g)^*\a_{\tilde g}(a(h)^*)\mu(\tilde g, \tilde h)a(gh), \quad
\tilde g, \tilde h\in \wtH,
$$
where $g$ and $h$ are the $H$-component of $\tilde g$ and $\tilde h$ respectively. We then
examine if $\tmu$ remains standard:
$$\aligned
\tmu(g, s&; h, t)=a(g)^*\th_s(\a_g(a(h)^*)\a_g(d_\mu(s; h))\muh(g, h)a(gh)\\
&=\a_g(\th_s(a(h))^*a(h))\a_g(d_\mu(s; h))a(g)^*\a_g(a(h)^*)\muh(g, h)a(gh)\\
&=\a_g\Big(\th_s(a(h))^*a(h))d_\mu(s; h)\Big)a(g)^*\a_g(a(h)^*)\muh(g, h)a(gh).
\endaligned
$$
Therefore with 
$$\aligned
d_\tmu(s; h)&=\th_s(a(h)^*)a(h)d_\mu(s; h), \quad s\in \R, h\in H;\\
\tmu(g, h)&=a(g)^*\a_g(a(h)^*)\muh(g, h) a(gh), \quad g, h \in H,
\endaligned
$$
we confirm that $\tmu$ is standard. Now the corresponding characteristic cocycle have the form:
$$\aligned
\la_\tmu(m&; g, s)=\tmu(g, s; g\inv m g))\tmu(m; g, s)^*\\
&=\a_g(d_\tmu(s; g\inv m g))\tmu_H(g; g\inv mg)\tmu_H(m; g)^*\\
&=\a_g\Big(\th_s(a(g\inv m g)^*)a(g\inv m g)d_\mu(s; g\inv mg)\Big)\\
&\hskip.5in\times a(g)^*\a_g(a(g\inv mg)^*)a(mg)a(m)a(g)a(mg)^*\\
&\hskip.5in\times
\muh(g; g\inv mg)\muh(m; g)^*\\
&=\a_g\Big(\th_s(a(g\inv m g)^*)d_\mu(s; g\inv mg)\Big)\\
&\hskip.5in\times
a(m) \muh(g; g\inv mg)\muh(m; g)^*\\
&=\a_{g, s}(a(g\inv m g)^*)a(m)\la_\mu(m; g, s).
\endaligned
$$
With $g=1$, we get
$$\aligned
& \la_\tmu(m; s)=\th_s(a(m)^*)a(m)\la(m; s)=1, \quad m\in L, g\in H, s\in
\R.
\endaligned
$$
This completes the proof.
\QED
\enddemo

So we replace the original characteristic cocycle $(\la, \mu)$ by the modified one $(\la_\tmu,
\tmu)$ by Lemma 2.10 so that 
$$
\la=\la_\mu \quad \text{and}\quad d_\mu(s; m)=1, \quad m\in L, s \in \R,
$$
and $\mu\in \tzatw(\wtH, A)$ is standard.

Now we use the fact that the HJR map $\dhjr$ pushes $(\la_\mu, \mu)$ to $c_\mu\in
\tzasth(\wtQ, A)\subset \tzath(\wtQ, A)$ which is cobounded by $f$ of  Lemma 2.4 (iii):
$$
f(\tp, \tq)=\a_p(d_\mu(s; \dfs(q))^*)\muh(\dfs(p), \dfs(q))^*\muh(\fnl(p, q); \dfs(pq))\in A.
$$
We examine $\part_\th (f|_Q)$ by making use of the relation between $d_\mu$ and $\muh$ in
the formula (2.9):
$$\aligned
\th_s(f(q, r))&f(q, r)^*\\
&=\th_s\Big(\muh(\dfs(q), \dfs(r))^*\muh(\fnl(q, r); \dfs(qr))\Big)\\
&\hskip .5in\times\{\muh(\dfs(q), \dfs(r))^*\muh(\fnl(q, r); \dfs(qr))\}^*\\
&=\th_s\Big(\muh(\dfs(q), \dfs(r))^*\Big)\muh(\dfs(q), \dfs(r))\\
&\hskip .5in\times\th_s\Big(\muh(\fnl(q, r); \dfs(qr)\Big)\muh(\fnl(q, r); \dfs(qr))^*\\
&=d_\mu(s; \dfs(q))\a_q(d_\mu(s; \dfs(r)))d_\mu(s; \dfs(q)\dfs(r))^*\\
&\hskip .5in\times \{d_\mu(s; \fnl(q, r))d_\mu(s; \dfs(qr))d_\mu(s; \fnl(q, r)\dfs(qr))^*\}^*\\
&=d_\mu(s; \dfs(q))\a_q(d_\mu(s; \dfs(r))) \{d_\mu(s; \fnl(q, r))(d_\mu(s; \dfs(qr)))\}^*.
\endaligned
$$
Next we compare this with $\partq f$ computed in the proof of Lemma 2.4. (iii). Substituting
$p, q, r$ in place of $\tp, \tq, \tr$ in the last expression of $\part_\wtQ f$, we obtain
$$\aligned
(\partq f)(p&, q, r)\\
&=\muh(\dfs(p);  \fnl(q, r))
\muh(\dfs(p)\fnl(q, r)\dfs(p)\inv ; \dfs(p))^*\\ 
&\hskip .3in\times\muh(\dfs(p)\fnl(q, r)\dfs(p)\inv;
\fnl(p, qr)) \Big\{\muh(\fnl(p, q); \fnl(pq, r))\Big\}^*.
\endaligned
$$
Combining the above two coboundary calculations, we obtain:
$$\aligned
(\part_\wtQ (f|_Q))(\tp&, \tq, \tr)\\
&=\a_p(\th_s(f(q, r))f(q, r)^*)(\partq (f|_Q))(p, q, r)\\
&=\a_p\Big(d_\mu(s; \dfs(q))\a_q(d_\mu(s; \dfs(r)))\\
&\hskip.5in\times
 \{d_\mu(s; \fnl(q, r))
  d_\mu(s; \dfs(qr)) \}^*\Big)\\
&\hskip.5in\times 
\muh(\dfs(p);  \fnl(q, r))
\muh(\dfs(p)\fnl(q, r)\dfs(p)\inv ; \dfs(p))^*\\ 
&\hskip .5in\times\muh(\dfs(p)\fnl(q, r)\dfs(p)\inv;
\fnl(p, qr))\\
&\hskip.5in\times 
\Big\{\muh(\fnl(p, q); \fnl(pq, r))\Big\}^*.
\endaligned
$$
Comparing this with $c_\mu$, we conclude
$$\aligned
(\part_\wtQ (f|_Q))(\tp&, \tq, \tr)=\a_p\Big(d_\mu(s; \dfs(q))\a_q(d_\mu(s; \dfs(r)))
d_\mu(s; \dfs(qr))^*\Big)c_\mu(\tp, \tq, \tr).
\endaligned
$$
Now we use the assumption that $\d(\la, \mu)=c_\mu\in \tbasth(\wtQ, A)$, which means the
existence of a new cochain $\xi\in \tC_\a^2(Q, A)$ such that
$$\aligned
c_\mu(\tp&, \tq, \tr)=\a_\tp(\xi(q, r))\xi(p, qr)\{\xi(p, q)\xi(pq, r)\}^*.
\endaligned
$$
Therefore, we get
$$\aligned
\a_\tp(f(&q, r))f(p, qr)\{f(p, q)f(pq, r)\}^*\\
&=\a_p\Big(d_\mu(s; \dfs(q))\a_q(d_\mu(s; \dfs(r)))
d_\mu(s; \dfs(qr))^*\Big)\\
&\hskip.2in\times \a_\tp(\xi(q, r))\xi(p, qr)\{\xi(p, q)\xi(pq, r)\}^*,
\endaligned
$$
equivalently
$$\aligned
\a_\tp((\xi^*f)(&q, r))(\xi^*f)(p, qr)\{(\xi^*f)(p, q)(\xi^*f)(pq, r)\}^*\\
&=\a_p\Big(d_\mu(s; \dfs(q))\a_q(d_\mu(s; \dfs(r)))
d_\mu(s; \dfs(qr))^*\Big).
\endaligned
$$
Setting $s=0$, we obtain $\partq (\xi^*f|_Q)=1$. With $p=1$, we get
$$\aligned
\th_s((\xi^*f)(&q, r))(\xi^*f)(q, r))^*
\\&=d_\mu(s; \dfs(q))\a_q(d_\mu(s; \dfs(r)))
d_\mu(s; \dfs(qr))^*.
\endaligned\tag2.16
$$
We now use the formula (2.9), which states that $d_\mu$ gives rise to an element 
$[d_\mu]\in \tZ_\a^1(H, \thth1)$. The assumption $\nu_\chi=1$ entails that the cocycle
$[d_\mu]$ factors through $Q$, i.e., there exists a map $a\colon\ (m, h)\in 
L\times H\mapsto a(m, h)\in A$ such that
$$\aligned
d_\mu(s; mh)&=\th_s(a(m, h))a(m, h)^*d_\mu(s; h), \quad m\in L, h\in H.\\
\endaligned\tag2.17
$$
We write $H$ in term of the cross-section $\dfs$ and the cocycle $\fnl$: $H=L\rtimes_\fnl Q$.
Writing $g=\txm_L(g)\dfs(\dpi(g)), h\in H,$  with 
$$
b(g)=a(\txm_L(g), \dfs(\dpi(g))\in A,\tag2.18
$$
we obtain
$$\aligned
d_\mu(s&; g)=\th_s(b(g))b(g)^*d_\mu(s; \dfs(\dpi(g))), \quad g\in H.
\endaligned\tag2.19
$$
Then the right hand side of the formula (2.9) becomes:
$$\aligned
d_\mu(s&; g)\a_g(d_\mu(s; h))d_\mu(s; gh)^*\\
&=\th_s(b(g))b(g)^*d_\mu(s; \dfs(\dpi(g)))\a_g\Big(\th_s(b(h))b(h)^*d_\mu(s;
\dfs(\dpi(h)))\Big)\\ &\hskip .5in\times
\Big(\th_s(b(gh))b(gh)^*d_\mu(s; \dfs(\dpi(gh)))\Big)^*\\
&=\th_s(b(g))b(g)^*\a_g\Big(\th_s(b(h))b(h)^*\Big)  
\Big(\th_s(b(gh)^*)b(gh) \Big) \\
&\hskip.5in\times 
\th_s\Big((\xi^*f)(\dpi(g), \dpi(h))\Big)(\xi^*f)(\dpi(g), \dpi(h)))^*\\
&=\th_s\Big(b(g)\a_g(b(h))(\xi^*f)(\dpi(g), \dpi(h))b(gh)^*\Big)\\
&\hskip.5in \times
\Big(b(g)\a_g(b(h))(\xi^*f)(\dpi(g), \dpi(h))b(gh)^*\Big)^*.
\endaligned
$$
Equating this to the left hand side of (2.9), we get
$$\left.\aligned
\th_s(\muh(g&, h))\muh(g, h)^*\\
&=\th_s\Big(b(g)\a_g(b(h))(\xi^*f)(\dpi(g),
\dpi(h))b(gh)^*\Big)\\ &\hskip.5in \times
\Big(b(g)\a_g(b(h))(\xi^*f)(\dpi(g), \dpi(h))b(gh)^*\Big)^*
\endaligned\right\}\quad g, h \in H.
$$
Hence $\mu_0=\dpi^*(\xi f^*)(\part_H b^*)\muh\in \tztw(H, \T)$.  Finally we compare 
$$
\Res(\mu_0)=(\la_{\mu_0}, \mu_0|_L)
$$
and  $(\la_\mu, \mu)$. First, we compare the
$\mu$-components of the characteristic cocycle and obtain
$$
\mu_0(m, n)=b(m)^*b(n)^*b(mn)\muh(m, n), \quad m, n\in L,
$$
since $(\xi f^*)(\dpi(m), \dpi(n))=1$. Second, we also get
$$\aligned
\la_{\mu_0}&(m; g, s)=\a_g(d_{\mu_0}(s; g\inv m g))\mu_0(g; g\inv mg)\mu_0(m; g)^*\\
&=\mu_0(g; g\inv mg)\mu_0(m; g)^*\\
&=(\part_H b^*)(g; g\inv mg)\muh(g; g\inv mg)(\part_H b^*)(m; g)^*\muh(m; g)^*\\
&=b(g)^*\a_g(b(g\inv mg)^*)b(mg)b(m)b(g)b(mg)^*\muh(g; g\inv mg)\muh(m; g)^*\\
&=\a_g(b(g\inv mg)^*)b(m)\muh(g; g\inv mg)\muh(m; g)^*\\
&=\a_g(b(g\inv mg)^*)b(m)\la_\mu(m; g, s).
\endaligned
$$
Therefore we conclude  that 
$$
\Res([\mu_0])=[\la_\mu, \mu]=\chi\in \La_\a(\wtH, L, M, A).
$$
This completes the proof of $\Ker(\d)\subset \text{\rm Im}(\Res)$ and so $\Ker(\d)=
\text{\rm Im}(\Res)$.

\proclaim{Lemma 2.11} There is a natural commutative diagram of exact sequences{\rm:}
$$\eightpoint\CD
\La_\ta(\wtH, L, M, A)@>\d>>\thasth(\wtQ,
A)*_\fs\Hom_G(N, \thth1(\R, A)))@>\Inf=\inf\scirc\part>>\thath(H, \T)\\
@VV i^*_{L, M}V@V\part VV@|\\
\La(H, M, \T)@>\dhjr>>\ththr(G, \T)@>\inf>>\ththr(H, \T)
\endCD
$$
\endproclaim
\demo{Proof} {\bf Map $\pmb \part\!:$}
Fix a cross-section $\sz\!: \thth1(\R, A)\mapsto \tzth1(\R,
A)$ and set
$$
\z_\nu(s; n)=(\sz(\nu(n))_s\in A, \quad s\in\R,n\in N,\ \nu\in
\Hom_G(N, \thth1(\R, A)).
$$
Choose $([c], \nu)\in \thasth(\wtQ, A)
*_\fs\Hom_G(N, \thth1(\R, A)))$ so that
$$
\part_2[c]=\nu\cup \fn_\fs,
$$
i.e.,
$$
[d(\ \cdot\ ; q, r)]=\nu(\fn_\fs(q, r))\qquad\text{in}\quad \thth1(\R, A),
\quad q, r\in Q.
$$
Hence there exists $f\in \tcatw(Q, A)$ such that
$$
d_c(s; q, r)=\th_s(f(q, r))f(q, r)^*\z_\nu(s; \fn(q, r)),\quad q, r \in Q, s\in \R, \tag2.20
$$
and therefore,
$$\aligned
c(\tp, \tq, \tr)&=c_Q(p, q, r)\a_p\Big(\th_s(f(q, r))f(q, r)^*\z_\nu(s; \fn(q, r))\Big),\\
&\hskip 1.5in \tp=(p, s), \tq, \tr\in \wtQ.
\endaligned\tag2.21
$$
The necessary condition for $c\in \tzasth(\wtQ, A)$ in (2.12) gives the
following:
$$\aligned
\th_s(c_Q&(p, q, r))c_Q(p, q, r)^*\\
&=\a_p\Big(\th_s(f(q, r))f(q, r)^*\z_\nu(s; \fn(q, r))\Big)\\
&\hskip .5in\times
\th_s(f(p, qr))f(p, qr)^*\z_\nu(s; \fn(p, qr))\\
&\hskip .5in\times \{\th_s(f(p, q))f(p, q)^*\z_\nu(s; \fn(p, q))\\
&\hskip .5in\times
\th_s(f(pq, r))f(pq, r)^*\z_\nu(s; \fn(pq, r))\}^*
\endaligned\tag2.22
$$
for each $p, q, r\in Q$ and $s\in \R$.

Now we are going to consider the pull back $\tpi^*(c)\in \tzath(\wtG,
A)$  with 
$$
\tpi(g, s)=(\pi(g), s),\quad \tilde g=(g, s)\in \wtG=G\times \R,
$$
But we first check the pull back $\tpi^*(\nu\cup \fnn)\in \tzatw(G,
\thth1)$. To this end, with $\text m_N(g)=g\fs(\pi(g))\inv\in N$ and
$\txn_N(g)=\fs(\pi(g))g\inv\in N$, we observe first
$$\aligned
\fnn(\pi(g)&, \pi(h))=\fs(\pi(g))\fs(\pi(h))\fs(\pi(gh))\inv, \quad g, h\in
G,\\
&=\txn_N(g) g \txn_N(h) h \{\txn_N(gh) gh\}\inv\\
&=\txn_N(g) g\txn_N(h) g\inv \txn_N(gh)\inv;
\endaligned
$$
and that
$$\aligned
\nu(\fnn(&\pi(g),\pi(h)))\\
&=\nu(\txn_N(g))\a_g(\nu(\txn_N(h))
\nu(\txn_N(gh))\inv \ \text{in}\ \thth1, 
\endaligned\quad g, h \in G.
$$
Hence we can choose $a(g, h)\in A, g, h\in G,$ such that
$$\aligned
\z_\nu(s; \fnn(\pi(g), \pi(h)))&=\th_s(a(g, h))a(g, h)^*\z_\nu(s; \txn_N(g))\\
&\hskip .5in\times
\a_g(\z_\nu(s; \txn_N(h)))\z_\nu(s; \txn_N(gh))^*
\endaligned\tag2.23
$$
We apply now this to the pull back of the above (2.22) and 
obtain for each $g, h, k\in G$:
$$\aligned
\th_s&(c_Q(\pi(g), \pi(h), \pi(k)))c_Q(\pi(g), \pi(h), \pi(k))^*\\
&=\a_g\Big(\th_s( \pi^*(f)(h, k)) \pi^*(f)(h, k)^*\th_s(a(h, k))a(h, k)^*\z_\nu(s; \txn_N(h))\\
&\hskip .5in\times
\a_h(\z_\nu(s; \txn_N(k)))\z_\nu(s; \txn_N(hk))^*\Big)\\
&\hskip .5in\times \th_s(\pi^*(f)(g, hk))\pi^*(f)(g, hk)^*\th_s(a(g, hk))
a(g, hk)^*\\
\endaligned
$$
$$\aligned
&\hskip .5in \times\z_\nu(s; \txn_N(g))
\a_g(\z_\nu(s; \txn_N(hk)))\z_\nu(s; \txn_N(ghk))^*\\
&\hskip .5in\times \Big\{\th_s(\pi^*(f)(g, h))\pi^*(f)(g, h)^*
\th_s(a(g, h))a(g, h)^*\\
&\hskip .5in \times\z_\nu(s; \txn_N(g))
\a_g(\z_\nu(s; \txn_N(h)))\z_\nu(s; \txn_N(gh))^*\\
&\hskip .5in \times
\th_s(\pi^*(f)(gh, k))\pi^*(f)(gh, k)^*\th_s(a(gh, k))a(gh, k)^*\\
&\hskip .5in\times\z_\nu(s; \txn_N(gh))
\a_{  gh }(\z_\nu(s; \txn_N(k)))\z_\nu(s; \txn_N(ghk))^*\Big\}^*\\
&=\a_g\Big(\th_s(\pi^*(f)(h, k))\pi^*(f)(h, k)^*\th_s(a(h, k))a(h, k)^*  \Big)\\
 &\hskip .5in\times \th_s(\pi^*(f)(g, hk))\pi^*(f)(g, hk)^*\th_s(a(g, hk))a(g,
hk)^*)\\ 
&\hskip .5in\times
\Big\{\th_s(\pi^*(f)(g, h))\pi^*(f)(g, h)^*\th_s(a(g, h))a(g, h)^*)\\
&\hskip.5in\times
\th_s\Big(\pi^*(f)(gh, k)\Big)  \pi^*(f)(gh, k)^*\th_s(a(gh, k))a(gh, k)^*\Big\}^*
\endaligned
$$
and hence, for each each $g, h, k \in G$,
$$\aligned
\th_s\Big(c_Q(&\pi(g), \pi(h), \pi(k))\part_G(\pi^*(f)^*a^*)(g, h, k)\Big)\\
&=c_Q(\pi(g), \pi(h), \pi(k))\part_G(\pi^*(f)^*a^*)(g, h, k).
\endaligned
$$
The ergodicity of the flow $\th$ yields that 
$$
\pi^*(c_Q)\part_G(\pi^*(f)a)^*\in \tzthr(G, \T).
$$

Now we change the cocycle $c$ to $c'$ within the cohomology class,
i.e., $c'=(\part_\wtQ b)c$ with $b\in \tcatw(Q, A)$ which gives:
$$\aligned
d'(s; q, r)=\th_s(b(q, r))b(q, r)^*d(s; q, r), \quad s\in \R, q, r \in Q.
\endaligned
$$
We also change  the cross-section $\sz:\thth1\mapsto \tzth1$ to $\fs_\tZ':
\thth1\mapsto \tzth1$. Then there exists a map $n\in N\mapsto e(n)\in A$
such that
$$
\z'_\nu(s; n)=\fs_\tZ'(\nu)_s=\th_s(e(n))e(n)^*\z_\nu(s; n), \quad n\in N.
$$
Thus we obtain, for each $s\in \R, g, h \in G$,
$$\aligned
d'&(s; q, r)=\th_s\Big((f(q, r)e(\fnn(q, r))^*b(q, r)\Big)\\
&\hskip .5in 
f(q, r)^*b(q, r)^* e(\fnn(q, r))\z_\nu'(s; \fnn(q, r));\\
\endaligned
$$
$$\aligned
\z_\nu'&(s; \fnn(\pi(g), \pi(h)))\\
&=\th_s(e(\fnn(\pi(g), \pi(h))))e(\fnn(\pi(g), \pi(h)))^*
\z_\nu(s; \fnn(\pi(g), \pi(h)))\\
& =\th_s\Big(e(\fnn(\pi(g), \pi(h)))a(g,
h)e(\txn_N(g))^*\a_g(\text e(m(h))^*)e(\txn_N(gh))\Big)\\
&\hskip .5in\times e(\fnn(\pi(g), \pi(h)))^*a(g, h)^*e(\txn_N(g)^*)\\
&\hskip .5in\times
\a_g(e(\txn_N(h))^*) e(\txn_N(gh))\\
&\hskip .5in\times\z_\nu'(s; \txn_N(g))
\a_g(\z_\nu'(s; \txn_N(h)))\z_\nu'(s; \txn_N(gh))^*;\\
c_Q'&(p, q, r)=(\part_\wtQ b)(p, q, r)c_Q(p, q, r)\\
&=(\part_Qb)(p, q, r)c_Q(p, q, r), \quad p, q, r\in Q.
\endaligned
$$
Therefore, the cochains $f$ and $a$ are transformed to the following
$f'$ and $a'$:
$$\aligned
f'(p, q)&=f(p, q)e(\fnn(p, q)  )^*b(p, q), \quad p, q \in Q;\\
a'(g, h)&=e(\fnn(\pi(g), \pi(h)))a(g, h)e(\txn_N(g))^*\\
&\hskip .5in\times
\a_g(e(\txn_N(h))^*) e(\txn_N(gh)), \quad g, h \in G.\\
\endaligned
$$
Thus we get
$$\aligned
\pi^*(c_Q')\part_G(\pi^*(f')a')^*
&= \pi^*(c_Q)\pi^*(\part_Q b)\part_G\Big(\pi^*(fb)\pi^*(e\scirc \fnn )a\Big)^* \\
&\hskip 1in\times
\part_G(\pi^*(e\scirc \fnn))
\part_G^2(e\scirc \txn_N)^*)\\
&=\pi^*(c_Q)\part_G(\pi^*(f)a)^*.
\endaligned
$$
Finally, in the choice of $f$ and $a$ we have pecisely the ambiguity of
$\tctw(Q, \T)$ and $\tctw(G, \T)$ which result the change on
$\pi^*(c_Q)\part_G(\pi^*(f)a)^*\in  \tzthr (G, \T)$ by $\tbthr(G, \T)$. Thus we
have a well-defined homomorphism
$$
\part: \thasth(\wtQ, A)*_\fs\Hom_G(N, \thth1(\R, A))  \mapsto \ththr(G, \T).
$$
which depends on the choice of the section $\fs \,\coron Q\longmapsto G. $

Now fix $\chi=[\la, \mu]\in \La_\a(\wtH, L, M, A)$ and set 
$$
([c^{\la, \mu}], \nu_\chi)=\d(\chi)\in \thasth(\wtQ, A)*_\fs\Hom_G(N, \thth1(\R, A))  .
$$
Associated with $(\la, \mu)$ is an $\wtH$-equivariant exact square:
$$\CD
@.1@.1@.1\\
@.@VVV@VVV@VVV\\
1@>>>\T@>>>A@>>>\tB@>>>1\\
@.@VVV@Vi VV@VVV\\
1@>>>V@>>>E@>>>F@>>>1\\
@.@VVV@Vj V\big\uparrow \sj V@VVV\\
1@>>>M@>>>L@>>>N@>>>1\\
@.@VVV@VVV@VVV\\
@.1@.1@.1
\endCD
$$
The far left column exact sequence corresponds to the restriction 
$$
(\la, \mu)|_{M\times H}=i^*_{L, M}(\la, \mu)\in \La(H, M, \T).
$$
To go further, we need the following:

\proclaim{Sublemma 2.12} In the above context, we have
$$
\dhjr(i^*_{L, M}(\chi))=\part (\d(\chi))\in \ththr(G, \T), \quad \chi\in \La_\a(\wtH, L, M, A).
$$
\endproclaim
\demo{Proof} First we arrange the cross-sections  in the following way:
$$\aligned
&\hskip .5in \sh(n\fs(p))=\sh(n)\sh(\fs(p)), \quad n\in N, p\in Q;\\
& \sj(m\sh(n))=\sj(m)\sj(\sh(n)),\quad m\in M, n\in N. 
\endaligned
$$
We further arrange the cross-sections $\sh$ on $N$, $\sdpart$ and $\sz$ on $\tH$,
so that they satisfy the following composition rules:
$$\CD
E@>\part_\th>>\tZ\\
@Vj V\big\uparrow \sj V@VV\big \uparrow \sz V\\
L@>\dot \part>\underset{\sdpart}\to\longleftarrow>\tH\\
@V\pig V \big\uparrow \sh V@A\nu A\big\downarrow \fs_\nu A\\
N@=N
\endCD\hskip .3in 
\aligned
&\sz=\part_\th\scirc\sj\scirc \sdpart:\tH\mapsto \tZ;\\
&\part_\th\scirc\sj\scirc\sh=\part_\th\scirc\sj\scirc\sdpart\scirc \nu=\sz\scirc
\nu: N\mapsto \tZ;
\\ &\partth\scirc \sj=\sz\scirc \nu\scirc\pig: L\mapsto \tZ.
\endaligned
$$
As $\sh\neq\fs_{\dot\part}\scirc \nu$ if $M\neq K(\la, \mu)$,  i.e., if $\nu$
is not injective, the second composition rule needs to be justified. For each
$n\in N$, set $m=\fs_{\dot \part}\scirc \nu(n)\sh(n)\inv
\in M$ so that $m\sh(n)=\fs_{\dot\part}(\nu(n))$. Then we have
$$\aligned
\sj(\fs_{\dot\part}&(\nu(n)))=\sj(m)\sj(\sh(n));\\
\sz(\nu(n))&=\partth(\sj(\fs_{\dot\part}(\nu(n)))=\partth(\sj(m)\sj(\sh(n)))\\
&=\partth(\sj(\sh(n))),
\endaligned
$$
which justfies the second composition rule.

For each $\ell\in L$, we write 
$$
\ell=\txmm(\ell)\sdpart(\dot\part(\ell))=\txmm(\ell)\sdpart(\nu(\pig(\ell)))
$$
and obtain
$$\aligned
\sj(\ell)&=\sj(\txmm(\ell)) \sj\scirc\sdpart(\nu(\pig(\ell))) , \quad \ell\in L;\\
\part_\th(\sj(\ell))&=\part_\th (\sj\scirc\sdpart(\nu(\pig(\ell))) =\sz(\nu(\pig(\ell))).
\endaligned
$$
Each $g\in G$ is uniquely written in the form:
$$
g=\txmn(g)\fs(\pi(g)), \quad g\in G,
$$
with $\txmn(g)\in N$. Therefore we have
$$
\sh(g)=\sh(\txmn(g))\dfs(\pi(g)), g\in G, \sh(\txmn(g))\in L.
$$
Then the product $gh$ of each pair $g, h\in G$ gives:
$$\aligned
\txmn(gh)&\fs(\pi(gh))= gh\\
&=\txmn(g)\fs(\pi(g))\txmn(h)\fs(\pi(h))\\
&=\txmn(g)\fs(\pi(g))\txmn(h)\fs(\pi(g))\inv\fs(\pi(g))\fs(\pi(h))\\
&=\txmn(g)\fs(\pi(g))\txmn(h)\fs(\pi(g))\inv\fnn(\pi(g), \pi(h))\fs(\pi(gh));\\
1&=\txmn(g)\fs(\pi(g))\txmn(h)\fs(\pi(g))\inv \fnn(\pi(g), \pi(h))\txmn(gh)\inv.
\endaligned
$$
We observe the following relation between the cocycles $\fnl$ and $\fnn$.
$$\aligned
\pig(\fnl(p&, q))=\pig(\dfs(p)\dfs(q)\dfs(pq))\inv)\\
&=\fs(p)\fs(q)\fs(pq)\inv=\fnn(p, q), \quad p, q \in Q.
\endaligned
$$
We then further compute:
$$\aligned
\fnm(g, h)&=\sh(g)\sh(h)\sh(gh)\inv, \quad g, h \in G,\\
&=\sh(\txmn(g)\fs(\pi((g)))\sh(\txmn(h)\fs(\pi((h)))\\
&\hskip .5in\times
\sh(\txmn(gh)\fs(\pi((gh)))\inv\\
&=\sh(\txmn(g))\sh(\fs(\pi((g))))\sh(\txmn(h))\sh(\fs(\pi(h))\\
&\hskip .5in\times
\{\sh(\txmn(gh))\sh(\fs(\pi(gh)))\}\inv\\
&=\sh(\txmn(g))\dfs(\pi(g)))\sh(\txmn(h))\dfs(\pi(g))\inv 
\dfs(\pi(  g ))\dfs(\pi(h))\\
&\hskip .5in\times
\{\sh(\txmn(gh))\dfs(gh)))\}\inv\\
&=\sh(\txmn(g))\dfs(\pi(g)))\sh(\txmn(h))\dfs(\pi(g))\inv 
\fnl(\pi(  g )),\pi(h))\\
&\hskip .5in\times
\{\sh(\txmn(gh))  \}\inv.
\endaligned
$$
We now take the cross-section $\sj$ and choose $b(g, h)\in A$ so that the
following computation is valid:
$$\aligned
\sj(\fnm(g&, h))=\sj\Big(\sh(g)\sh(\txmn(h))\sh(g)\inv\\
&\hskip .5in\times
 \sh(\txmn(g))\fnl(\pi(g), \pi(h))
\sh(\txmn(gh))\inv\Big)\\
\endaligned
$$
$$\aligned
&=b(g, h)\a_{\sh(g)}(\sj(\sh(\txmn(h))))\sj(\sh(\txmn(g)))\\
&\hskip 1in\times\sj(\fnl(\pi(g),\pi(h))) \sj(\sh(\txmn(gh)))\inv\\
&=b(g, h)\a_{\sh(\txmn(g))\dfs(\pi(g))}(\sj(\sh(\txmn(h))))\sj(\sh(\txmn(g)))\\
&\hskip 1in\times\sj(\fnl(\pi(g),\pi(h))) \sj(\sh(\txmn(gh)))\inv\\
&=b(g, h)\sj(\sh(\txmn(g)))\a_{\dfs(\pi(g))}(\sj(\sh(\txmn(h))))\\
&\hskip 1in\times\sj(\fnl(\pi(g),\pi(h))) \sj(\sh(\txmn(gh)))\inv.
\endaligned
$$
We summerlize this here for later use:
$$\aligned
\sj(\fnm(g&, h))=b(g, h)\sj(\sh(\txmn(g)))\a_{\dfs(\pi(g))}(\sj(\sh(\txmn(h))))\\
&\hskip 1in\times\sj(\fnl(\pi(g),\pi(h))) \sj(\sh(\txmn(gh)))\inv
\endaligned\tag2.24
$$
We then apply the coboundary operator $\partth$ to the both side to obtain:
$$\aligned
1&=\partth(b(g, h))\partth(\sj(\txmn(g))\partth(\a_{\dfs(\pi(g))}(\sj(\sh(\txmn(h))))\\
&\hskip .5in\times \partth(\sj(\fnl(\pi(g),\pi(h))))\partth( \sj(\sh(\txmn(gh)))\inv)
\endaligned
$$
and use the compostion rules among cross-sections to drive:
$$\aligned
\sz(\nu(&\fnn(\pi(g),\pi(h)))=\partth(b(g,
h)\inv)\sz(\nu(\txmn(g)\inv)\\
&\hskip.5in\times\a_g(\sz(\nu(\txmn(h))\inv))
\sz(\nu(\txmn(gh)))).
\endaligned
$$
Since $\txnn(g)=\txmn(g)\inv, g\in G$, we have
$$\aligned
\sz(\nu(&\fnn(\pi(g),\pi(h)))=\partth(b(g,
h)\inv)\sz(\nu(\txnn(g))\\ 
&\hskip1in\times
\a_g(\sz(\nu(\txnn(h))))\sz(\nu(\txnn(gh))\inv)),
\endaligned
$$
equivalently
$$\aligned
\z_\nu(s; \fnn(&\pi(g), \pi(h))=\th_s(b(g, h)\inv)b(g, h)\z_\nu(s; \txnn(g))\\
&\hskip .5in\times 
\a_g(\z_\nu(s; \txnn(h)))\z_\nu(s; \txnn(gh))^*, \quad g, h \in G.
\endaligned
$$
Therefore the elements $b(g, h)\inv\in A$ serves as $a(g, h)$ of (2.23) in
the construction of  $\part (\d(\chi))$.

With $u(g)=\sj(\sh(\txmn(g)))\in E, g\in G,$ and $w(g, h)=\sj(\fnl(\pi(g), \pi(h))),$ 
we apply the coboundary operation $\partg$ to the both side of (2.24) relative to the
outer action $\a_{\sh}$ of $G$ on $E$ to obtain:
$$\aligned
c^{\la, \mu}_G(g&, h, k)= \a_{\sh(g)}\big(\sj(\fnm(h,
k))\big)\sj(\fnm(g, hk))\\
&\hskip .5in\times
\{\sj(\fnm(g, h))\sj(\fnm(gh, k))\}\inv \\
&=(\partg a)(g, h, k)\inv\a_{\sh(g)}\Big(u(h)\a_{\dfs(\pi(h))}(u(k))w(h,
k)u(hk)\inv\Big)\\ &\hskip.5in\times
u(g)\a_{\dfs( \pi(g) )}(u(hk))w(g, hk)u(ghk)\inv\\
&\hskip.5in\times
\{u(g)\a_{\dfs( \pi(g) )}(u(h))w(g, h)u(gh)\inv u(gh)\\
&\hskip.5in\times\a_{\dfs( \pi(gh) )}(u(k))w(gh,
k)u(ghk)\inv\}\inv\\
&=(\partg a)(g, h, k)\inv\\
&\hskip .5in\times
\Ad(u(g))\scirc\a_{\dfs(\pi(g))}\Big(u(h)\a_{\dfs(\pi(h))}(u(k))w(h, k)u(hk)\inv\Big)\\
&\hskip.5in\times u(g)\a_{\dfs( \pi(g) )}(u(hk))w(g, hk)\\ &\hskip.5in\times
\{u(g)\a_{\dfs( \pi(g) )}(u(h))w(g, h)u(gh)\inv u(gh)\\
&\hskip.5in\times\a_{\dfs(\pi(gh))}(u(k))w(gh, k)\}\inv\\
&=(\partg a)(g, h, k)\inv\a_{\dfs(\pi(g))}\Big(\a_{\dfs(\pi(h))}(u(k))w(h,
k)\Big)\\ &\hskip.5in\times
w(g, hk) \{w(g, h)\a_{\dfs(\pi(gh))}(u(k))w(gh, k)\}\inv\\
&=(\partg a)(g, h, k)\inv w(g, h)\a_{\dfs(\pi(gh))}(u(k)))\\
&\hskip .5in\times
w(g, h)\inv \a_{\dfs(\pi(g) ) }(w(h, k))\\ &\hskip.5in\times
w(g, hk) \{w(g, h)\a_{\dfs(\pi(gh))}(u(k))w(gh, k)\}\inv\\
&=(\partg a)(g, h, k)\inv
\a_{\dfs(\pi(g) ) }(w(h, k))w(g, hk) \{w(g, h)w(gh, k)\}\inv\\
&=(\partg a)(g, h, k)\inv c^{\la,\mu}_Q(\pi(g), \pi(h), \pi(k)).
\endaligned
$$
Since the cochain $f\in \tcatw(Q, A)$ of (2.20) is taken to be 1 in our case, the
above  computation shows that the third cohomology class:
$$
[\partg (a\inv) \pi^*(c^{\la, \mu}_Q)]\in \ththr(G, \T)
$$
is indeed the 3-cohomology class associated with the far
left column $H$-equivariant exact sequence of the exact square before the lemma:
$$\CD
1@>>>\T@>>>V@>>>M@>>>1
\endCD
$$ 
which corresponds to the characteristice invariant $i^*_{L, M}(\chi)\in 
\La(H, M, \T)$.
\QED
\enddemo
\enddemo

\subhead{$\pmb{\text{\rm Im}(\d)\subset \Ker(\inf\scirc \part)}$}
\endsubhead
As seen above, we have $\part(\d(\chi))=\dhjr(i^*_{L, M}(\chi))\in \ththr(G, \T)$. Hence we 
conclude
$$
1=\inf\scirc \dhjr(i^*_{L, M}(\chi))=\inf\scirc \part (\d(\chi)),\quad  \chi\in \La_\a(H, L, M, A).
$$

\subhead{$\pmb{\text{\rm Im}(\d)\supset \Ker(\inf\scirc \part)}$}
\endsubhead
First, we compare our sequence with {\hjr}:
$$\eightpoint\CD
\La_\a(\wtH, L, A)@>\dhjr>>\ththra(\wtQ, A)@>\inf>>\ththra(\wtH, A)\\
@AAA@AAA@A {i_{A, \T}}_*AA\\
\La_\a(\wtH, L, M, A)@>\d>>\thasth(\wtQ, A)*_\fs\Hom_G(N, \thth1(\R, A))
@>\Inf>>\ththr(H, \T)
\endCD
$$
Now suppose that $\Inf([c], \nu)=1$ in $\ththr(H, \T)$. The 3-cocycle $c\in \tzasth(\wtQ, A)$ is
naturally an element of $\tzath(\wtQ, A)$. We denote this element by $\tilde c$ and its cohomology
class $[\tilde c]\in \thath(\wtQ, A)$. First, the image $\part([c], \nu)\in \ththr(G, \T)$ is obtained
as the class of $\partg(\pi^*(f)a\inv)\pi^*(c_Q)$ where $f\in \tcatw(Q, A)$ and $a\in \tcatw(G,
A)$ are obtained subject to the following conditions:
$$\aligned
d_c(s &; q, r)=\th_s(f(q, r))f(q, r)^*\z_\nu(s; \fnn(q, r)),\quad s\in \R, q, r \in Q;\\
\z_\nu(s&; \fn(\pi(g), \pi(h))=\th_s(a(g, h))a(g, h)^*\z_\nu(s; \txnn(g)), \quad g, h \in G,\\
&\hskip 1in\times\a_g(\z_\nu(s; \txnn(h)))\z_\nu(s; \txnn(gh))^*,
\endaligned\tag2.25
$$
where $\z_\nu(s; n)=\sz(\nu(n))_s, n\in N, s\in \R$. The image $\Inf([c], \nu)$ is obtained as
the class of 
$$
\pi_H^*(\partg(\pi^*(f)a)\inv\pi^*(c_Q))=\part_H(\dpi^*(f)\pi_H^*(a))\inv\dpi^*(c_Q)\in \tzthr(H,
\T).
$$
The assumption that $\Inf([c], \nu)=1$ means that
$\part_H(\dpi^*(f)\pi_H^*(a))\inv\dpi^*(c_Q)\in \tbthr(H, \T)$, i.e., there exists 
$b\in \tctw(H, \T)$ such that
$$\aligned
\part_H(\dpi^*(f)\pi_H^*(a))\inv\dpi^*(c_Q)=\parth b
\endaligned
$$
Hence for each triple $g, h, k\in H$ we have
$$\aligned
c_Q(\dpi(g)&, \dpi(h), \dpi(k))=\a_g\Big(b(h, k)f(\dpi(h), \dpi(k))a(\pig(h),
\pig(k))\Big)\\
&\hskip.2in\times
b(g, hk)f(\dpi(g), \dpi(hk))a(\pig(g), \pig(hk))\\
&\hskip.2in\times \Big\{b(g, h)f(\dpi(g), \dpi(k))a(\pig(g),\pig(k))\\
&\hskip.2in\times b(gh, k)
f(\dpi(gh), \dpi(k))a(\pig(gh),\pig(k))\Big\}^*. 
\endaligned
$$
With $u(g, h)=a(\pig(g), \pig(h))b(g, h)f(\dpi(g), \dpi(h))$, we get
$$
\dpi^*(c_Q)=\parth u,
$$
and
$$\aligned
c(\dpi(\tilde g)&, \dpi(\tilde h), \dpi(\tilde k))=\a_g(d_c(s; \dpi(h), \dpi(k))) c_Q(\dpi(g),
\dpi(h), \dpi(k)) \\
&=\a_g(d_c(s; \dpi(h), \dpi(k)))\a_g\Big(u(h, k)\Big)u(g, hk)\{u(g, h) u(gh, k)\}^*.
\endaligned
$$
The identities (2.25) yields the following computations, for each $g, h, k\in H$,
$$\aligned
d_c(s &; \dpi(h), \dpi(k))=\th_s(f(\dpi(h), \dpi(k)))f(\dpi(h), \dpi(k))^*
\z_\nu(s; \fn(\dpi(h), \dpi(k)));\\
\z_\nu(s&; \fnn(\dpi(g), \dpi(h))=\th_s(a(\pig(g), \pig(h)))a(\pig(g), \pig(h))^*\\
&\hskip.5in\times z_\nu(s; \txnn(\pig(g)))\a_g(\z_\nu(s; \txnn(\pig(h))))\z_\nu(s;
\txnn(\pig(gh)))^*;\\
d_c(s &; \dpi(h), \dpi(k))=\th_s(f(\dpi(h), \dpi(k)))f(\dpi(h), \dpi(k))^*\\
&\hskip.2in\times
\th_s(a(\pig(h), \pig(k)))a(\pig(h), \pig(k))^*\\
&\hskip.5in\times z_\nu(s; \txnn(\pig(h)))\a_h(\z_\nu(s; \txnn(\pig(k))))\z_\nu(s;
\txnn(\pig(hk)))^*.
\endaligned
$$
With $v(s; g)=\z_\nu(s; \txnn(\pig(g)))$, we get
$$\aligned
d_c(s &; \dpi(h), \dpi(k))=\th_s(u(h, k))u(h, k)^*v(s; g)\a_g(v(s; h))v(s; gh)^*\\
&=(\partth (u(h, k)))_s(\parth v)(s; h, k)
\endaligned
$$
Substituting this to the above computation of $\dpi^*(c)$ and setting
$$
w(g, s; h, t)=u(g, h)\a_g(v(s; h)^*),
$$ we obtain:
$$\aligned
(\dpi^* c)(\tilde g&, \tilde h, \tilde k)= \a_g\Big(\th_s(u(h, k))u(h, k)^*v(s; h)\a_h(v(s;
k))v(s; hk)^*\Big)\\
&\hskip .5in\times \a_g(u(h, k))u(g, hk)\{u(g, h) u(gh, k)\}^*; \
\endaligned
$$
$$\aligned
(\part_{\wtH}w)(\tilde g&, \tilde h, \tilde r)=\a_{\tilde g}(w(\tilde h; \tilde k))
w(\tilde h; \tilde h\tilde k)\{w(\tilde g; \tilde h)w(\tilde g\tilde h; \tilde k)\}^*\\
&=\a_g\scirc\th_s\Big(u(h, k)\a_h(v(t; k))^*))\Big)\\
&\hskip .5in\times u(g, hk)\a_g(v(s; hk)^*)\\
&\hskip .5in\times
\{u(g, h)\a_g(v(s; h))^*
u(gh, k)\a_{gh}(v(s+t; k)^*)\}^*\\
&=\a_g\scirc\th_s\Big(u(h, k)\a_h(v(t; k))^*))\Big)\\
&\hskip .5in\times u(g, hk)\a_g(v(s; hk)^*)
\Big\{u(g, h)\a_g(v(s; h)^*) u(gh, k)\\
&\hskip .5in\times\a_{gh}\big(  v(s; k)^*\th_s(v(t; k)^*) \big)\Big\}^*\\
&=\a_g\scirc\Big(\th_s(u(h, k))u(h, k)^*v(s, h)\a_h(v(s, k))v(s; hk)^*\Big)\\
&\hskip .5in\times \a_g(u(h, k))u(g, hk)
\{u(g, h)u(gh, k)\}^*\\
&=(\dpi^*c)(\tilde g, \tilde h, \tilde k).\\
\endaligned
$$
Therefore, we conclude
$$
\dpi^*(c)=\part_\wtH w.\tag2.26
$$
Hence the element $(\la, \mu)$ given by
$$\aligned
\la(\ell&; g, s)=w(g, s; g\inv \ell g)w(\ell; g, s)^*,\quad \ell\in L, g\in H, s\in \R;\\
&\mu(m, n)=w(m, n), \quad m, n \in L,
\endaligned
$$
is a characteristice cocycle in $\tZ_\ta(\wtH, L, A)$. In terms of the original $a, b$ and
$f$, we get
$$\aligned
\la(m&; g, s)=a(\pig(g), \pig(g\inv m g))b(g, g\inv m g)f(\dpi(g), \dpi(g\inv m g))\\
&\hskip.5in\times
\a_g(\z_\nu(s; \txnn(\pig(g\inv m g)))^*)a(\pig(m), \pig(g))^*\\
&\hskip.5in\times
b(m, g)^*f(\dpi(m),\dpi(g))^*;\\
\mu(m&; n)=a(\pig(m), \pig(n))b(m, n).
\endaligned\tag2.27
$$
Now observe that if $m, n\in M$, then both $\la$ and $\mu$ takes values in $\T$, so
that $\chi=[\la, \mu]\in \La_\a(\wtH, L, M, A)$. 

We are now going to compare the new cocycle $c^{\la, \mu}$ and the original $c$ in the next
lemma to complete the proof of Lemma 2.12 and therefore Theorem 2.7:
\proclaim{Lemma 2.13} The cochain $W\in \tcatw(\wtQ, A)$  defined by
$$
W(\tp, \tq)=w(\dfs(\tp), \dfs(\tq))w(\fnl(p, q), \dfs(\tp\tq))^*, \quad \tp=(p, s), \tq=(q, t)\in \wtQ,
$$
falls in $\tcatw(Q, A)$ and its coboundary $\part_\wtQ W$ bridges the difference 
between
$(c^{\la, \mu}, \nu_{\chi})$ and the original $(c, \nu)$, i.e., $([c], \nu)=\d(\chi)$. Therefore
$$
\Ker(\Inf)\subset \text{\rm Im}(\d).
$$
\endproclaim
\demo{Proof} First we observe that for any pair $\tp=(p, s), \tq=(q, t)\in \wtQ$
$$\aligned
W(\tp, \tq)&=
w(\dfs(\tp), \dfs(\tq))w(\fnl(p, q), \dfs(\tp\tq))^*\\
&=u(\dfs(\tp), \dfs(\tq))\a_p(v(s;
\dfs(q)) u(\fnl(p, q), \dfs(pq))^*\\
&=a(\fs(p), \fs(q))b(\dfs(p), \dfs(q))f(p, q)\a_p(\z_\nu(s; \txnn(\fs(q)))\\
&\hskip.5in\times
a(\fn(p, q), \fs(pq))^*b(\fnl(p, q), \dfs(pq))^*f(p, q)^*\\
&=a(\fs(p), \fs(q))b(\dfs(p), \dfs(q))f(p, q)\quad(\text{as } \txnn(\fs(q))=1)\\
&\hskip.5in\times
a(\fn(p, q), \fs(pq))^*b(\fnl(p, q), \dfs(pq))^*f(p, q)^*\\
&=W(p, q).
\endaligned
$$
Thus $W$ is constant on $\R$-variables, so that it  belongs to $\tcatw(Q, A)$.

By Lemma 2.1, we have
$$
c=(\part_\wtQ W)c^{\la, \mu}
$$
and therefore
$$
c^{\la, \mu}\equiv c\quad \mod\ \tbasth(\wtQ, A),\ \text{i.e.,}\quad [c^{\la, \mu}]=[c]\quad
\text{in}\ \thasth(\wtQ, A).
$$

Setting $g=1$ in (2.27), we obtain for each $m\in L$
$$\aligned
\la(m; s)&=a(1, \pig( m ))b(1,  m )f(1, \dpi(m))\\
&\hskip.5in\times
\z_\nu(s; \txnn(\pig(m)))^* a(\pig(m), 1)^*\\
&\hskip.5in\times
b(m, 1)^*f(\dpi(m),1)^*\\
&=\z_\nu(s; \txnn(\pig(m)))^*=\z_\nu(s; \pig(m)\inv))^*\\
\endaligned
$$
since the cochains $a, b$ and $f$ can be chosen such a way that whenever $1$ appears in the
arguments they take value $1$. As 
$$
\z_\nu(\ \cdot; \pig(m)\inv)^*\equiv \z_\nu(\ \cdot; \pig(m))\quad \mod\ \tbth1,
$$
i.e., 
$\nu(\pig(m)\inv)\inv =\nu(\pig(m)), m \in L,$ we conclude that $\nu=\nu_{[\la, \mu]}$.
Therefore we conclude that $([c], \nu)=\d([\la, \mu])$. This completes the proof of the inclustion
$\Ker(\Inf)\subset \text{\rm Im}(\d)$.
\QED
\enddemo

\proclaim{Lemma 2.14} Let $A$ denote the unitary group $\sU(\sC)$ of an abelian separable {\vna} $\sC$ or
the torus group $\T$. Let $\a$ be an action of a countable discrete group $G$ on $\sC$.
 To each $c\in \tZ_\a^3(G, A)$, there corresponds a countable group 
$H=H(c)$ and a normal subgroup $M=M(c)$ such that{\rm:}
\roster
\item"i)"  the group $G$ is identified with the quotient group $H/M${\rm;}
\item"ii)" there corresponds a characteristic cocycle 
$$
(\la, \mu)=(\la_c, \mu_c)\in \tZ_\a(H, M, A)
$$ 
such that 
$$
[c]=\dhjr([\la, \mu]) 
$$
in the {\hjr} relative to $\{H, M, A\}${\rm;}
\item"iii)" the group $M$ is abelian.
\endroster
\endproclaim
\demo{Proof} First extend the coefficient group $A$ to the unitary group 
$$B=\sU(\sC\botimes \ell^\infty(G))
$$
on which $G$ acts by $\a\otimes \rho$ with $\rho$ the right translation action of
$G$ on $\ell^\infty(G)$, which will be denoted by $\a$ again whenever it will not cause any
confusion, and obtain an exact sequence:
$$\CD
1@>>>A@>i>>B@>j>\underset \sj\to \longleftarrow>C@>>>1,
\endCD
$$
where $i(a)=a\otimes 1\in B, a\in A,$ and $\sj$ a cross-section which can be fixed without reference to the
cocycle $c$. Then set
$$ 
u(x, g, h)=u_c(x, g, h)=\a_x^{-1}(c(x, g, h))\in A, \quad x, g, h\in G,
$$
and view $u(\cdot, g, h)$ as an element of $B=\sU(\ell^\infty(G)\botimes \sC))=\Map(G, A)$. The cocycle
identity gives that $c=\part_{\a}^G u\in \tB_{\a}^2(G, C)$. Since
$j(A)=\{1\}$ in $C$, $\mu=\mu_c=j_*(u)$ is in $\tZ_{\a}^2(G, C)$. Let $M$ be the subgroup of
$C$ generated by the saturation $\{\a_g(\mu(h, k)): g, h, k\in G\}$ of the range of $\mu$, so that 
$\mu$ belongs to $\tZ_\a^2(G, M)$. Now consider the twisted semi-direct product:
$$
H=H(c)=M\rtimes_{\a, \mu} G
$$
and obtain an exact sequence: 
$$\CD
1@>>> M@>>> H@>\pig>> G@>>> 1.
\endCD
$$
Set $E=j\inv(M)$ to obtain a crossed extension $E\in \X_{\a}(H, M, A)$. With $\sh$ the
cross-section given by
$$
\sh(g)=(1, g)\in H, \quad g \in G,
$$
we obtain
$$
\mu(g, h)=\sh(g)\sh(h)\sh(gh)\inv, \quad g, h \in G.
$$
Thus  $\mu\in \tZ_{\a}^2(G, M)$. Now observe that 
$$
f(g, h)=\sj(\mu(g, h))\inv u(g, h)\in i(A)
$$
and that
$$
(\part_G f)(g, h, k)=\part_G(\sj\scirc \mu)(g, h, k)\inv c(g, h, k), \quad g, h, k\in G.
$$
Thus we conclude that $\dhjr(\chi_E)=[\part_G(\sj\scirc \mu)]=[c]\in \tH_\a^3(G, A)$.
\QED
\enddemo

\noindent
{\bf Remark 2.15.} The last lemma shows that if $G$ is a countable discrete amenable group, then so is 
$H$ because $M$ is abelian and countable, and the quotient $G=H/M$ is amenable. Another important fact is
that the groups $H$ and $M$ depend heavily on the choice of the cocycle $c$. Two cohomologous 
cocycles $c, c'\in \tZ_\a^3(G, A)$  need not produce isomorphic $H(c)$ and $H(c')$. In fact, the subgroups
$M(c)$ and
$M(c')$ are not isomorphic. We will address this incovenience later. If we use the entire $C$ in place of $M$,
then the resulting groups are isomorphic in a natural way. But in this way,
we will lose
the countability of $H$.

\vskip.1in
\noindent
{\bf Definition 2.16.} We call the group $H(c)$ the {\it resolution group of the cocycle} $c\in \tzath(G, A)$
and the characteristic coycle $(\la_c, \mu_c)\in \tZ_\a(H, M, A)$ a {\it resolution} of
the cocycle $c$. We also call the map $\pig:H(c)\mapsto G$ {\it resolution map} and the pair
$\{H(c), \pig\}$ a {\it resolution system}.

\proclaim{Corollary 2.17} Let $\{\sC, \R, \th\}$ be an ergodic flow and $G$ a discrete countable group acting
on the flow $\{\sC, \R, \th\}$ via $\a$. Let $N$ be a normal subgroup of $G$ such that $N\subset \Ker(\a)$.
Then with $Q=G/N$ the quotient group of $G$ by $N$ and $\fs\!: Q\mapsto G$ a cross-section of the quotient
map $\pi\!: G\mapsto Q$, for any pair
$$
([c],\nu)\in \thasth(Q\times \R, \sU(\sC))*_\fs\Hom_G(N, \thth1(\R, \sU(\sC)))
$$
there exist a countable discrete group $H$ and a surjective homomorphism $\pig\!: H\mapsto
G$ and $\chi\in
\La_{\pi_G^*(\a)}(H\times \R, L, M, A)$ such that 
$$
([c], \nu)=\d(\chi)
$$
where $L=\pi_G^{-1}(N)$, $M=\Ker(\pig)$ and $\d$ is the modified {\rm HJR}-map in 
{\rm Lemma 2.11}
associated with the exact sequence{\rm:}
$$\CD
1@>>>M@>>>L@>\pig>>G@>\pi>\underset{\fs}\to\longleftarrow>Q@>>>1.
\endCD
$$
Moreover, the kernel $M=\Ker(\pig)$ is chosen to be abelian. Hence
if $G$ is amenable in addition, then $H$ is amenable.
\endproclaim
\demo{Proof} Let $\part$ be the map in Lemma 2.11:
$$
\part: \thasth(Q\times \R, \sU(\sC))*_\fs\Hom_G(N, \thth1(\R, \sU(\sC)))\mapsto \ththr(G, \T).
$$
Set $[c_G]=\part([c], \nu)\in \ththr(G, \T)$ and choose a cocycle $c_G\in \tzthr(G, \T)$ which represents the
cohomology class $[c_G]$. Let $H=H(c_G)$ be the resolution group of $c_G$ in
Lemma 2.14, i.e., the group $H$ is equipped with a surjective homomorphism $\pig\!:
H\mapsto G$ such that $\pi_G^*([c_G])=1$ in $\ththr(H, \T)$. Thus with 
$L=\pi_G^{-1}(N)\triangleleft H$  and $M=\Ker(\pig)\triangleleft H$,
we have an exact sequence:
$$\CD
1@>>>M@>>>L@>\pig>>G@>\pi>\underset{\fs}\to\longleftarrow>Q@>>>1,
\endCD
$$
with specified cross-section $\fs$ of $\pi\!: G\mapsto Q$. This generates the associated
modified {\hjr} of (2.13) in Theorem 2.7. As
$$
\Inf([c], \nu)=\pi_G^*(\part([c], \nu)=\pi_G^*([c_G])=1,
$$
there exists $\chi\in \La_\a(\wtQ, L, M, A)$ such  that $([c], \nu)=\d(\chi)$ where
$\wtQ=Q\times \R$ and $A=\sU(\sC)$ of course. 
If $G$ is amenable, then the group $H$ in Lemma 2.14 must be amenable as the subgroup
$M$ of $H$ in Lemma 2.14 is abelian and the quotient group $H/M=G$ is amenable.
\QED
\enddemo

\subhead\nofrills{Change on the Cross-Section $\pmb{\fs\!: Q\mapsto G}$:}
\endsubhead\
As mentioned repeatedly, the group $\thasth(Q\times G, A)*_\fs\Hom_G(N, \thth1)$ depends
heavily on the cross-section $\fs\!: Q\mapsto G$. So we are going to examine what change
occurs if we change the cross-section from $\fs\!: Q\mapsto G$ to another $\fs'\!: Q\mapsto G$.
The change  does not alter the groups $\thasth(\wtQ, A)$ not $\Hom_G(N, \thth1)$, but the
fiber product $\thasth(Q\times G, A)*_\fs\Hom_G(N, \thth1)$ changes to $\thasth(Q\times G,
A)*_\fsp\Hom_G(N, \thth1)$.

\proclaim{Proposition 2.18} In the setting as above, if $\fsp\!: Q\mapsto G$ is another
cross-section of the homomorphism $\pi\!: G\mapsto Q=G/N$, then there is a natural
isomorphism 
$$
\sig_{\fsp, \fs}\!:\ \thasout(G\times\R, N, A)\mapsto \thasth(\wtQ,
A)*_\fsp\Hom_G(N, \thth1),\tag2.28
$$
where $\wtQ=Q\times \R$ as before.
Furthermore, if $\frak s\dprime\!: Q\mapsto G$ is the third cross-section of $\pi$, then the
ismomorphisms satisfy the following chain rule{\rm:}
$$
\sig_{\frak s\dprime, \fs}=\sig_{\frak s\dprime, \fsp}\scirc\sig_{\fsp, \fs}\tag2.29
$$
\endproclaim
\demo{Proof} The new  cross-section $\fsp\!:Q\mapsto G$ generates a new  $N$-valued 
2-cocycle:
$$
\fnnp(p, q)=\fsp(p)\fsp(q)\fsp(pq)', \quad p, q \in Q.
$$
Set
$$
\txnsps(p)=\fsp(p)\fs(p)\inv\in N, \quad p\in  Q.
$$
Then the 2-cocycle $\fnnp$ is written in terms of $\fnn$ and $\txnsps$ as follows:
$$\left.\aligned
\fnnp(p, q)&=\txnsps(p)\fs(p)\txnsps(q)\fs(q)\{\txnsps(pq)\sp(pq)\}\inv\\
&=\txnsps(p)\fs(p)\txnsps(q)\fs(p)\inv \fnn(p, q)\txnsps(pq)\inv
\endaligned\right\}\quad p, q \in Q.
$$
Hence for each $\nu\in \Hom_G(N, \thth1)$ we have
$$
\nu(\fnnp(p, q))=\nu(\txnsps(p))\a_p(\nu(\txnsps(q)))\nu(\fnn(p, q))\nu(\txnsps(pq))\inv.
$$
For each $([c], \nu)\in \thasout(G\times\R, N, A)$, we set
$$\aligned
d_{c'}(s; q, r)&=d_c(s; q, r)\z_\nu(s; \txnsps(q))
\a_q(\z_\nu(s; \txnsps(r)) \z_\nu(s; \txnsps(qr))^*;\\
&\hskip .5in c_Q'(p, q, r)=c_Q(p, q, r), \quad s\in \R, p, q, r\in Q,
\endaligned
$$
where $\z_\nu(s; n)=\sz(n)_s, n\in N, s\in \R$. As  $\partq d_{c'}=\partq d_c$, the new pair
$(d_{c'}, c_Q')$ gives a standard 3-cocycle
$c'\in  \tzasth(\wtQ, A)$ which is not  congruent to $c=(d_c, c_Q)\in \tzasth(\wtQ, A)$
modulo $\tbasth(\wtQ, A)$ in general although they are congruent modulo $\tbath(\wtQ, A)$.
Now we define the map $\sigsps$ in the  following way:
$$\aligned
\sigsps([c], \nu)=([c'], \nu),\quad ([c], \nu)\in \thasout(G\times\R, N, A).
\endaligned
$$
Then as 
$$
[d_{c'}(\ \cdot, q, r)]=\nu(\fnnp(q, r))\in \thth1, \quad q, r\in Q,
$$
the pair $([c'], \nu)$ belongs to $\thasth(\wtQ, A)*_{\fsp}\Hom_G(N, \thth1)$.

To check the multiplicativity of $\sigsps$, for each pair $h, k\in \thth1$ we choose
$a(h, k)\in A$ such that
$$
\sz(h)\sz(k)=\partth (a(h, k))\sz(hk).
$$
Then for each pair $([c], \nu), ([\bar c], \bar \nu)\in \thasout(G\times\R, N, A)$,
we have
$$\aligned
d_{c'}(\ \cdot\ ; q, r)&=d_c(\ \cdot\ ; q, r)
\sz(\nu(\txnsps(q)))\a_q(\sz(\nu(\txnsps(r))))\\
&\hskip .5in\times\sz(\nu(\txnsps(qr)))\inv;\\
d_{\bar c'}(\ \cdot\ ; q, r)&=d_{\bar c}(\ \cdot\ ; q, r)
\sz(\bar\nu(\txnsps(q)))\a_q((\bar\sz(\nu(\txnsps(r))))\\
&\hskip .5in\times
\sz((\bar\nu(\txnsps(qr)))\inv;\\
\endaligned
$$
$$\aligned
(d_{c'}&d_{\bar c'})(\ \cdot\ ; q, r)=(d_c d_{\bar c})(\ \cdot\ ; q,
r)\sz(\nu(\txnsps(q)))\\ 
&\hskip 1.5in\times
\a_q(\sz(\nu(\txnsps(r))))\sz(\nu(\txnsps(qr)))\inv\\
&\hskip.5in\times \sz(\bar\nu(\txnsps(q)))\a_q(\sz(\bar\nu(\txnsps(r))))
\sz(\bar\nu(\txnsps(qr)))\inv;\\
&=(d_c d_{\bar c})(\ \cdot\ ; q, r)\sz(\nu\bar\nu(\txnsps(q)))
\a_q(\sz(\nu\bar\nu(\txnsps(r))))
\sz(\nu\bar\nu(\txnsps(qr)))\inv\\
&\hskip .5in\times
\partth(a(\nu(\txnsps(q)), \bar\nu(\txnsps(q))))\partth(\a_q(a(\nu(\txnsps(r)),
\bar\nu(\txnsps(r))))\\
&\hskip .5in  \times
\partth(a(\nu(\txnsps(qr), \bar \nu(\txnsps(qr))))\inv\\
&=d_{(c\bar c)'}(\ \cdot\ ; q, r)(\partth\part_Q b))(q, r),
\endaligned
$$
where $b\in \tC_\a^1(Q, A)$ is given by
$$
b(q)=a(\nu(\txnsps(q)), \bar \nu(\txnsps(q)))\in A.
$$
Also we have
$$
(c\bar c)_Q'=(c\bar c)_Q=c_Q\bar c_Q(\part_Q\part_Q b)=c_Q'\bar c_Q'.
$$
Therefore, we get $[c'][\bar c']=[(c\bar c)']$ in $\thasth(\wtQ, A)$ by Lemma 2.5, i.e., $\sigsps$ is
multiplicative. 

The chain rule follows from the definition of $\sigsps$. We leave it to the reader. The chain
rule also gives that the map $\sigsps$ is an isomorphism.
\QED
\enddemo

The chain rule (2.29) allows us to define the cohomology group independent of the cross-section in
the following way: first let  $S$ be the set of all cross sections $\fs\!: Q\mapsto G$ of the
homomorphism $\pi$ and set:
$$\aligned
\thaout(G, N, A)
&=\Big\{\{([c], \nu)_\fs\!: \fs \in S\}\in\prod_{\fs\in S}\thasout(G\times\R, N, A):\\
&\hskip .5in ([c], \nu)_\fsp=\sigsps(([c], \nu)_\fs),\quad \fsp, \fs \in S\Big\}.
\endaligned\tag2.30
$$

\vskip.1in
\noindent
{\bf Definition 2.19.} The group 
$\tH_{\a}^\out(G, N, A)$ will be called {\it the modular obstruction group}. Each pair
$(c, \z)\in \tzasth(\wtQ, A)\times \Map(N, \tzth1)$ which gives rise to an element $([c],
[\z])\in \thaout(G, N, A)$ will be called a {\it modular
obstruction cocycle}.

\head{\bf \S3. Outer Actions of a Discrete Group on a Factor}
\endhead
Let $\sM$ be a separable factor. Associated with $\sM$ is the characteristic square:
$$\CD
@.1@.1@.1\\
@.@VVV@VVV@VVV\\
1@>>>\T@>>>A@>\partth>>\tbth1(\R, A)@>>>1\\
@.@VVV@VVV@VVV\\
1@>>>\sU(\sM)@>>>\tsU(\sM)@>\partth>>\tzth1(\R, A)@>>>1\\
@.@V\Ad VV@V\tAd VV@VVV\\
1@>>>\Int(\sM)@>>>\cntr(\sM)@>\dpartth>>\thth1(\R, A)@>>>1\\
@.@VVV@VVV@VVV\\
@.1@.1@.1\\
\endCD\tag3.1
$$
where $A=\sU(\sC)$ is the unitary group of the flow $\{\sC, \R, \th\}$ of weights on $\sM$,
which is $\Aut(\sM)\times \R$-equivariant. Applying the previous section to the groups
$$\aligned
&H=\Aut(\sM), M=\Int(\sM), G=\Out(\sM), N= \cntr(\sM) \\
&Q=\Outt(\tM)=\Out(\sM)/\thth1(\R, A)=\Aut(\sM)/\cntr(\sM),
\endaligned
$$
we obtain the intrinsic invariant and the intrinsic modular obstruction:
$$\aligned
&\Theta(\sM)\in \La_{\mod\times\th}(\Aut(\sM)\times \R, \cntr(\sM), A);\\
&\Obm(\sM)\in \tH_{\mod\times \th, \txs}^\out(\Out(\sM), \thth1(\R, A), A).
\endaligned
$$
Choosing a cross-section: $g\in \Out(\sM)\mapsto \a_g\in \Aut(\sM)$, we obtain an
outer action of $\Out(\sM)$ on $\sM$, i.e., 
$$\aligned
\a_g\scirc \a_h&\equiv \a_{gh} \quad \mod\ \Int(\sM), g, h\in \Out(\sM);\\
\a_{\id}&=\id; \quad \a_g\not\in \Int(\sM) \quad \text{if}\ g\neq \id.
\endaligned
$$
Choosing $\{u(g, h)\in \sU(\sM): g, h \in \Out(\sM)\}$ so that
$$
\a_g\scirc \a_h=\Ad(u(g, h))\a_{gh},\quad g, h \in \Out(\sM),
$$
we obtain a 3-cocycle $c\in \tzthr(\Out(\sM), \T)$:
$$
c(g, h, k)=\a_g(u(h, k))u(g, hk)\{u(g, h)u(gh, k)\}^*, \quad g, h, k\in \Out(\sM).
$$
Its cohomology class $[c]\in \ththr(\Out(\sM), \T)$ does not depend on the choice of
the cross-section $\a:g\in \Out(\sM)\mapsto \a_g\in \Aut(\sM)$ nor on the choice of
$\{u(g, h)\}$. The {\it intrinsic obstruction} $\Ob(\sM)=[c]$ of $\sM$ is, by
definition, the cohomology class $[c]\in \ththr(\Out(\sM), \T)$.

\proclaim{Proposition 3.1} The intrinsic obstruction $\Ob(\sM)$ of the factor
$\sM$ is the image $\part(\Obm(\sM))$ of the instrinsic modular obstruction $\Obm(\sM)$
of $\sM$ under the map 
$$
\part: \tH_{\mod\times \th, \txs}^\out(\Out(\sM), \thth1, A)\mapsto \ththr(\Out(\sM), \T)
$$
given by {\rm Lemma 2.11.}
\endproclaim
\demo{Proof} In the notations in the last section, we take $\Aut(\sM)$ for $H$,
$\Out(\sM)$ for $G$, $\Outt(\tM)$ for $Q$, $\cntr(\sM)$ for $L$, $\Int(\sM)$ for $M$ and $N$
for $\thth1(\R, A)$. Then with 
$$
\tilde\chi=\Theta(\sM)\in
\La_{\mod\times \th}(\Aut(\sM)\times\R, \cntr(\sM), A)
$$
the characteristic square (3.1) gives that $M=K(\tilde\chi)$.
The characteristic invariant $\chi \in \La(\Aut(\sM), \Int(\sM), \T)$ associated with the
$\Aut(\sM)$-equivariant exact sequence:
$$\CD
1@>>>\T@>>>\sU(\sM)@>>>\Int(\sM)@>>>1
\endCD
$$
is precisely $\chi=i_{\cntr(\sM), \Int(\sM)}^*(\tilde \chi)$ the pull back in Lemma 2.11. Then
the obstruction $\Ob(\sM)=\dhjr(\chi)$ is $\part(\Obm(\sM))=\part(\d(\tilde \chi))$ by Lemma
2.11.
\QED
\enddemo

Therefore, the modular obstruction $\Obm(\sM)$ contains the information carried by the
obstruction $\Ob(\sM)$.

 Let $G$ be a countable group. Fix a free outer action $\a$ of $G$ on $\sM$.
If $\dot \a_g$ is the class of $\a_g$ in $\Out(\sM)$, then the map $\dot \a\!: g\in G\mapsto \dot \a_g\in
\Out(\sM)$ is an injective homomorphism. With $N=\dot \a\inv (\thth1(\R, A))\triangleleft G$, the
quotient map
$\pi\!: g\in G\mapsto \pi(g)=gN\in Q=G/N$
and a cross-section $\fs\!: Q\mapsto G$ of $\pi$, we get the modular obstruction 
$$
\Obm(\a)\in \thasout(G\times\R, N,  A).
$$
Two outer actions $\a$ and $\b$ of $G$ on the same factor $\sM$ are, by definition, {\it outer
conjugate} if there exists an automorphism $\sig\in \Aut(\sM)$ such that
$$
\dot \b_g=\dot\sig\dot\a_g\dot \sig\inv, \quad g \in G,
$$
where $\dot \sig\in \Out(\sM)$ is the class of $\sig$ in $\Out(\sM)$, i.e., $\dot \sig=\sig \Int(\sM)
\in \Out(\sM)$.

\proclaim{Theorem 3.2} Let $G$ be a countable discrete group and $\sM$ a separable
infinite factor with flow of weights $\{\sC, \R, \th\}$.
Suppose that $\a\!: g\in G\mapsto \a_g\in \Aut(\sM)$ is a free outer action of $G$ on
$\sM$. Set
$N=N(\a)=\a\inv(\cntr(\sM))$, $Q=G/N$ and fix a cross-section $\fs\!: Q\mapsto G$ of the quotient
map $\pi\!: G\mapsto Q$.

{\rm i)} The modular obstruction{\rm:} 
$$
\Obm(\a)\in \tH_{\mod(\a)\times\th\!,\ \txs}^3(Q\times \R,
\sU(\sC))*_\fs\Hom_G(N, \thth1(\R,
\sU(\sC)))
$$ 
is an invariant for the outer conjugacy class of $\a$.

{\rm ii)} If $\sM$ is an approximately finite dimensional factor and $G$ is amenable, then  
the triplet $(N(\a), \mod(\a), \Obm(\a) ) $ is a complete invariant of the outer conjugacy class of $\a$
in the sense that if $\b: G\mapsto \Aut(\sM)$ is another outer action of $G$ on $\sM$ such that
$N(\a)=N(\b)$, and there exists an automorphism $\sig\in \Aut_\th(\sC)$ such that
$$
\sig\scirc\mod(\a_g)\scirc\sig\inv=\mod(\b_g), \quad g\in G;\quad \sig_*(\Obm(\a))=\Obm(\b),
$$
then the automorphism $\sig$ of $\sC$ can be extended to an automorphism denoted by $\sig$ again
to the non-commutative flow of weights $\{\tM, \R, \th, \tau\}$ such that
$$
\sig\scirc \a_g\scirc \sig\inv \equiv \b_g\quad \mod \ \Int(\sM), \quad g \in G.
$$
\endproclaim
\demo{Proof} We continue to denote the unitary group $\sU(\sC)$ by $A$ for short. Let $[c^\a]\in
\ththr(G, \T)$ be the obstruction $\Ob(\a)$ and $c^\a\in \tzthr(G, \T)$ represent
$[c]$ which is obtained  by fixing a family $\{u(g, h)\in \sU(\sM)\!:\ g, h \in G\}$
such that
$$
\a_g\scirc \a_h=\Ad(u(g, h))\scirc \a_{gh}, \quad g, h \in G,
$$
and by setting
$$
c^\a(g, h, k)= \a_g(u(h, k))u(g, hk)\{u(g, h)u(gh, k)\}^*\in \T,\quad g, h, k\in G.
$$
Let $\pig\!: H=H(c^\a)\mapsto G$ be the resolution group of the cocycle $c^\a\in
\tzthr(G, \T)$ and the resolution map, i.e., $\pi_G^*(c^\a)\in \tbthr(H, \T)$.
Choose $b: h\in H\mapsto b(h)\in\T$ such  that
$$
c^\a(\pig(g), \pig(h), \pig(k))=b(h, k)b(g, hk)\{b(g, h)b(gh, k)\}^*, \quad g, h, k\in
H.
$$
Setting
$$
\bar u_H(g,  h)=b(g, h)^*u(\pig(g),\pig(h)),  \quad g, h \in H,
$$
we obtain
$$
\a_{\pig(g)}(\bar u_H(h, k))\bar u_H(g, hk)\{\bar u_H(g, h)\bar u_H(gh, k)\}^*=1.
$$
Hence $\{\a_\pig, \bar u_H\}$ is a cocycle twisted action of $H$. Then by \cite{ST1:
Theorem 4.13, page 156}, there exits a family $\{v_H(g)\in \sU(\sM):\ g\in H\} $
such that
$$
\bar u_H(g, h)=\a_{\pig(g)}(v_H(h)^*)v_H(g)^*v_H(gh), \quad g, h \in H,
$$
so that the map 
$$
g\in H\mapsto \b_g=\Ad(v_H(g))\scirc \a_{\pig(g)}\in \Aut(\sM)
$$
is an action of $H$ on $\sM$ as seen below:
$$\aligned
\b_g\scirc \b_h&=\Ad(v_H(g))\scirc \a_{\pig(g)}\scirc \Ad(v_H(h))\scirc \a_{\pig(h)}\\
&=\Ad(v_H(g)\a_{\pig(g)}(v_H(h)))\scirc \a_{\pig(g)}\scirc \a_{\pig(h)}\\
&=\Ad(v_H(g)\a_{\pig(g)}(v_H(h)))\scirc \Ad(\bar u_H(g, h))\scirc \a_{\pig(gh)}\\
&=\Ad(v_H(gh))\scirc \a_{\pig(gh)}=\b_{gh}, \quad g, h \in H.
\endaligned
$$
Therefore, the outer action $\a_\pig$ is perturbed to an action $\b$ by inner
automorphisms. With $\sh$ a cross-section of $\pig$, the map $\dot \b\!:\ g\in
G\mapsto \b_{\sh(g)}\in\Aut(\sM)$ is an outer action of $G$ on $\sM$ which is an
inner perturbation of the original outer action $\a$ because
$$\aligned
\dot \b_g&=\b_{\sh(g)}=\Ad(v_H(\sh(g))\scirc \a_{\pig(\sh(g))}\\
&=\Ad(v_H(\sh(g))\scirc \a_g, \quad g\in G.
\endaligned
$$
Hence we may and do replace the outer action $\a$ by $\dot \b$. Then the outer action
$\a$ is given by an action $\b$ of $H$ in the following way:
$$
\a_g=\b_{\sh(g)},\qquad g\in G.
$$
The action $\b$ of $H$ gives rise to the characteristic invariant $\chi(\b)\in \La(H, M,
\T)$ with $M=\Ker(\pig)=\a\inv(\Int(\sM))$, so that  the obstraction $\Ob(\a)$ becomes
the image $\dhjr(\chi(\b))$ of $\chi(\b)$ under the HJR-map $\dhjr$.

 i) We only need to prove that the modular obstruction is unchanged under the
perturbation by inner automorphisms. Choose $\{w(p, q): p, q \in Q\}\i \tsU(\sM)$
so that
$$
\a_p\scirc \a_q=\tAd(w(p, q))\scirc \a_{pq}, \quad p, q \in Q,
$$
where $\a_p$ means $\a_{\fs(p)}$ for short. We write $\a_\tp\in \Aut(\tM)$ for $\a_p\scirc \th_s,
\tp=(p, s)\in \wtQ=Q\times \R$. Then for each triple $\tp=(p, s), \tq=(q, t), \tr=(r, u)\in
\wtQ$, the cocycle $c=c^\a$ representing $\Obm(\a)$ is given by:
$$\aligned
c^\a(\tp, \tq, \tr)&=\a_{\tp}(w(q, r))w(p, qr)\{w(p, q)w(pq, r)\}^*\\
&=\a_p(\th_s(w(q, r))w(q, r)^*)\a_p(w(q, r))w(p, qr)\{w(p, q)w(pq, r)\}^*\\
&=\a_p(d(s; q, r))c_Q(p, q, r),
\endaligned
$$
where
$$\aligned
&d(s; q, r)=\th_s(w(q, r))w(q, r)^*;\\
c_Q(p, q, r)&=\a_p(w(q, r))w(p, qr)\{w(p, q)w(pq, r)\}^*.
\endaligned
$$
The $G$-equivariant homomorphism $\nu: N\mapsto \thth1(\R, A)$ is given by
$\nu_\a(m)=\dot\partth(\a_m)\in \thth1(\R, A), m \in N.$ Let $\{v(g): g\in G\}\i \sU(\sM)$ and
set 
$$
\b_g=\Ad(v(g))\scirc \a_g, \quad g\in G.
$$
Then we have, with $\b_p=\b_{\fs(p)}, p\in Q,$
$$\aligned
\b_p\scirc \b_q&=\Ad(v(\fs(p)))\scirc \a_p\scirc \Ad(v(\fs(q)))\scirc \a_q\\
&=\Ad(v(\fs(p))\a_p(v(\fs(q))))\scirc \a_p\scirc\a_q\\
&=\tAd(v(\fs(p))\a_p(v(\fs(q)))w(p, q))\scirc\a_{pq}\\
&=\tAd(v(\fs(p))\a_p(v(\fs(q)))w(p,
q)v(\fs(pq))^*)\scirc\Ad(v(\fs(pq)))\scirc\a_{pq}\\
&=\Ad(v(\fs(p))\a_p(v(\fs(q)))w(p, q)v(\fs(pq))^*)\scirc \b_{pq}, \quad p, q \in Q.
\endaligned
$$
Therefore, we get 
$$\aligned
c^\b(\tp&, \tq, \tr)=\b_\tp\Big(v(\fs(q))\a_q(v(\fs(r)))w(q, r)v(\fs(qr))^*\Big)\\
&\hskip.5in\times
v(\fs(p))\a_p(v(\fs(qr))w(p, qr)v(\fs(pqr))^*\\
&\hskip.5in\times
\Big(v(\fs(p))\a_p(v(\fs(q))w(p, q)v(\fs(pq))^*\\
&\hskip.5in\times
v(\fs(pq))\a_{pq}(v(\fs(r))w(pq, r)v(\fs(pqr))^*\Big)^*\\
&=\b_p\Big\{\th_s\Big(v(\fs(q))\a_q(v(\fs(r)))w(q, r)v(\fs(qr))^*\Big)\\
&\hskip.5in\times\Big(v(\fs(q))\a_q(v(\fs(r)))w(q,
r)v(\fs(qr))^*\Big)^*\Big\}\\
&\hskip.5in\times
\b_p\Big(v(\fs(q))\a_q(v(\fs(r)))w(q, r)v(\fs(qr))^*\Big)\\
&\hskip.5in\times
v(\fs(p))\a_p(v(\fs(qr))w(p, qr)v(\fs(pqr))^*\\
&\hskip.5in\times
\Big(v(\fs(p))\a_p(v(\fs(q))w(p, q)v(\fs(pq))^*\\
&\hskip.5in\times
v(\fs(pq))\a_{pq}(v(\fs(r))w(pq, r)v(\fs(pqr))^*\Big)^*\\
&=v(\fs(p))\a_p\Big(\th_s(w(q, r))w(q, r)^*\Big)v(\fs(p))^*\\
&\hskip.5in\times
v(\fs(p))\a_p\Big(v(\fs(q))\a_q(v(\fs(r)))w(q,
r)v(\fs(qr))^*\Big)v(\fs(p))^*\\
&\hskip.5in\times
v(\fs(p))\a_p(v(\fs(qr))w(p, qr)v(\fs(pqr))^*\\
&\hskip.5in\times
\Big(v(\fs(p))\a_p(v(\fs(q))w(p, q)v(\fs(pq))^*\\
&\hskip.5in\times
v(\fs(pq))\a_{pq}(v(\fs(r))w(pq, r)v(\fs(pqr))^*\Big)^*\\
&=\a_p\Big(\th_s(w(q, r))w(q, r)^*\Big)
\a_p\Big(\a_q(v(\fs(r)))w(q,
r)\Big)\\
&\hskip.5in\times
w(p, qr)\Big(w(p, q)
\a_{pq}(v(\fs(r))w(pq, r)\Big)^*\\
&=\a_p\Big(\th_s(w(q, r))w(q, r)^*\Big)
w(p, q)\a_{pq}(v(\fs(r)))w(p, q)^*\\
&\hskip.5in\times
\a_p(w(q, r))
w(p, qr)\Big(w(p, q)\a_{pq}(v(\fs(r))w(pq, r)\Big)^*\\
\endaligned
$$
$$\aligned
&=\a_p\Big(\th_s(w(q, r))w(q, r)^*\Big)
\a_p(w(q, r))
w(p, qr)\Big(w(p, q)w(pq, r)\Big)^*\\
&=c^\a(\tp, \tq, \tr).
\endaligned
$$
Therefore, the inner perturbation $\b$ of the outer action $\a$ of $G$ does not
change the modular obstruction cocycle $c^\a$, i.e., $c^\a=c^\b$ as seen above.
Hence $\Obm(\a)=\Obm(\b)$. This proves the assertion (i).

ii) Assume that $\sM$ is an approximately finite dimensional factor with
non-commutative flow $\{\tM, \R, \th, \tau\}$ of weights and the flow $\{\sC, \R, \th\}$
of weights on $\sM$, and suppose that $G$ is a countable discrete amenable group. Let
$\dot\a$ and $\dot\b$ be outer actions of $G$ on $\sM$ such that 
\roster
\item"a)" $N=N(\dot\a)=\dot \a\inv(\cnt(\sM))=N(\dot \b)=\dot\b\inv(\cnt(\sM))$;
\item"b)" $\mod(\dot\a_g)=\mod(\dot\b_g), \quad g\in G;$
\item"c)" with $Q=G/N$ and $A=\sU(\sC)$
$$
(\Obm(\dot\a),  \nu_{\dot\a})=(\Obm(\dot\b), \nu_{\dot\b})\in \thasth(Q\times
\R, A)*\Hom_G(N, \thth1(\R, A)).
$$
\endroster
We want to conclude from this data that the outer actions $\dot \a$ and $\dot \b$ of
$G$ are outer conjugate. The assumption (c) implies that
$$
\Ob(\dot\a)=\part((\Obm(\dot\a),  \nu_{\dot\a}))=\part((\Obm(\dot\b),  \nu_{\dot\b}))
=\Ob(\dot \b)\in \ththr(G, \T).
$$
Therefore, we may and do choose the same obstruction cocycle
$c=c^{\dot\a}=c^{\dot\b}$, which allows us to pick up the common resolution system
$\pig: H=H(c)\mapsto G$ and actions $\a$ and $\b$ of $H$ on $\sM$ which give
$\dot\a$ and $\dot \b$ respectively:
$$
\dot\a_g=\a_{\sh(g)}, \quad \dot\b_g=\b_{\sh(g)}, \quad g\in G.
$$
First, the resolution group $H$ is amenable by  Lemma 2.14 and 
the actions $\a$ and $\b$ of $H$ give rise to the following invariants:
$$\aligned
&L=\pi_G^{-1}(N)=\a\inv(\cnt(\sM))=\b\inv(\cnt(\sM)),\\
& \chim(\a),\  \chim(\b)\in \La_{\mod(\a)\times  \th}(H\times  \R, L, A).
\endaligned
$$
Since $M=\Ker(\pig)=\a\inv(\Int(\sM))=\b\inv(\Int(\sM))$, we have
$$
M=K(\chim(\a))=K(\chim(\b)).
$$
Therefore, the modular characteristic invariant
$\chim(\a)$ and $\chim(\b)$ are both members of $\La_{\a\times\th}(\wtH, L, M, A)$
with
$\wtH=H\times\R$, where we are now going to use $\a$ for $\mod(\a)=\mod(\b)$. The
resolution system $\{H, \pig\}$ generates the following modified {\hjr}:
$$\CD
\cdots@>>>\thtw(H,\T)@>\Res>>\La_{\a\times\th}(\wtH, L, M, A)\\
@>\d>>\thasth(\wtQ, A)@>\part>>\ththr(H, \T),
\endCD
$$
such that 
$$
\d(\chim(\a))=\Obm(\dot\a)=\Obm(\dot\b)=\d(\chim(\b)).
$$
With  this, our assetion (ii) follows immediately from the next theorem
\QED
\enddemo

\proclaim{Theorem 3.3} Let $\a$ and $\b$ be two actions of a countable discrete group
$H$ on an infinite factor $\sM$ with $L=\a\inv(\cntr(\sM))=\b\inv(\cntr(\sM))$ and
$M=\a\inv(\Int(\sM))=\b\inv(\Int(\sM))$. Let $G=H/M$ and $\pig\!: H\mapsto G$ be
the quotient map. Suppose that $\sh\!: G\mapsto H$ is a cross-section and set
$$
\dot\a_g=\a_{\sh(g)},\quad \dot \b_g=\b_{\sh(g)}, \quad g \in G,
$$
to obtain outer actions $\dot\a$ and $\dot \b$ of $G$ on $\sM$.

{\rm i)} The two outer actions $\da$ and $\db$ of $G$ are outer conjugate if and only
if the two original actions $\a$ and $\b$ of $H$ are outer conjugate.

{\rm ii)} If the two actions $\a$ and $\b$ of $H$ on $\sM$ are outer conjugate, then
there exists an automorphism $\sig\in \Aut_\th(\sC)$ such that
\roster
\item"a)" $\mod(\a)=\sig\scirc\mod(\b)\scirc\sig\inv;$
\item"b)" their modular characteristic invariants $\chim(\a)\in \La_{\a\times\th}(\wtH,
L, M, A)$ and $\sig_*(\chim(\b))\in \La_{\a\times\th}(\wtH, L, M, A)$, where
$\wtH=H\times\R$, have the same image{\rm:}
$$\aligned
\Obm(\da)&=\d(\chim(\a))=\d(\sig_*(\chim(\b)))=\Obm(\sig\db\sig\inv)\\
\endaligned
$$
in the group $\thasth(\wtQ, A)*\Hom_G(N, \thth1(\R, A))$
under the modified {\rm HJR}-map $\d$, where $N=L/M$, $Q=H/L=G/N$ and
$\wtQ=Q\times\R.$
\endroster

{\rm iii)} If $\sM$ is an approximately finite dimensional infinite factor and $H$ is
amenable in addition, then the existence of an automorphism $\sig\in
\Aut_\th(\sC)$ such that{\rm:}
\roster
\item"a)" $\mod(\a)=\sig\scirc\mod(\b)\scirc\sig\inv;$
\item"b)" their modular characteristic invariants $\chim(\a)\in \La_{\a\times\th}(\wtH,
L, M, A)$ and $\sig_*(\chim(\b))\in \La_{\a\times\th}(\wtH, L, M, A)$, where
$\wtH=H\times\R$, have the same image{\rm:}
$$
\d(\chim(\a))=\d(\sig_*(\chim(\b)))\in \thasth(\wtQ, A)*\Hom_G(N,
\thth1(\R, A))
$$
\endroster
is sufficient for $\a$ and $\b$ to be outer conjugate.
\endproclaim
\demo{Proof} i) It is obvious that the outer conjugacy of the outer actions $\da$ and 
$\db$ of $G$ follows from that of the original actions $\a$ and $\b$ of $H$. So
suppose that the outer actions $\da$ and  $\db$ of $G$ are outer conjugate, which
means the existence of an automorphism $\sig\in \Aut(\sM)$ and a family $\{u(g):
g\in G\}\subset \sU(\sM)$ of unitaries such that
$$
\Ad(u(g))\scirc \da_g=\sig\scirc \db_g\scirc \sig\inv, \quad g \in G.
$$
Writing each $h\in H$ in the form:
$$
h=\txmm(h)\sh(\pig(h)), \quad h\in H,\ \txmm(h)\in M,
$$
we have
$$\aligned
&\a_h=\a_{\txmm(h)}\scirc \da_{\pig(h)};\\
&\b_h=\b_{\txmm(h)}\scirc \db_{\pig(h)},
\endaligned \quad h\in H.
$$
As $\a_m$ and $\b_m$ are inner for each $m\in M$, they are in the following form:
$$\aligned
\a_m=\Ad(v(m)),\quad \b_m=\Ad(w(m)), \quad m\in M,
\endaligned
$$
for some $v(m), w(m)\in \sU(\sM)$. Therefore, we have, for each $h\in H$,
$$\aligned
\sig\scirc\b_h\scirc \sig\inv&=\sig\scirc \b_{\txmm(h)\sh(\pig(h))}\scirc\sig\inv
=\sig\scirc\b_{\txmm(h)}\scirc\db_{\pig(h)}\scirc\sig\inv\\
&=\sig\scirc\Ad(w(\txmm(h)))\scirc\sig\inv\sig\scirc \db_{\pig(h)}\scirc\sig\inv\\
&=\Ad(\sig(w(\txmm(h))))\scirc\Ad(u(\pig(h))\scirc \da_{\pig(h)}\\
&=\Ad(\sig(w(\txmm(h))))\scirc\Ad(u(\pig(h))\scirc
\a_{\txmm(h)}^{-1}\scirc\a_{h}\\
&=\Ad(\sig(w(\txmm(h))))\scirc\Ad(u(\pig(h))\scirc
\Ad(v(\txmm(h))^*)\scirc\a_{h}\\
&=\Ad(u(h))\scirc\a_h,
\endaligned
$$
where $u(h)=\sig(w(\txmm(h)))u(\pig(h))v(\txmm(h))^*$. Hence the actions $\a$
and $\b$ of $H$ are outer conjugate.

ii) Assume that the two actions $\a$ and $\b$ of $H$ on $\sM$ are outer conjugate.
Then there exist $\sig\in \Aut(\sM)$ and a family $\{u(h): h\in H\}\subset \sU(\sM)$
such that $u(1)=1$ and
$$
\sig\scirc\b_h\scirc \sig\inv=\Ad(u(h))\scirc \a_h, \quad h\in H.
$$
Since $\Int(\sM)$ acts on the flow of weights trivially, we have
$$
\mod(\sig)\scirc\mod(\b_h)\scirc\mod(\sig)\inv=\mod(\a_h), \quad h\in H,
$$
and conclude that $\mod(\sig)$ conjugates $\mod(\a)$ and $\mod(\b)$, i.e., the
assertion (a). Replacing $\b_g, g\in H,$ by $\sig\scirc\b_g\scirc \sig\inv, g\in H,$ we
may and do assume from now on for short that $\mod(\a)=\mod(\b)$ and
$$
\b_g=\Ad(u(h))\scirc \a_h, \quad h\in H.
$$
As $\a$ and $\b$ are both actions, we have
$$\aligned
\Ad(u(gh))&\scirc \a_{gh}=
\b_{gh}=\b_g\scirc \b_h =\Ad(u(g))\scirc
\a_g\scirc\Ad(u(h))\scirc \a_h\\
&=\Ad(u(g)\a_g(u(h)))\scirc\a_g\scirc\a_h\\
&=\Ad(u(g)\a_g(u(h)))\scirc\a_{gh}, \quad g, h \in H.
\endaligned
$$
Thus we get
$$
\mu(g, h)=u(g)\a_g(u(h))u(gh)^*\in \T, \quad g, h\in H,
$$
and $\mu\in \tztw(H, \T)$. Each $\a_m, m\in L,$ falls in $\cntr(\sM)$, so that it
is of the form:
$$
\a_m=\tAd(v(m)), \quad v(m)\in \tsU(\sM).
$$
As $\b_m=\Ad(u(m))\scirc \a_m$, we may choose $w(m)=u(m)v(m), m\in L$. The
unitary families $\{v(m): m\in L\}$ and $\{w(m): m \in L\}$ generate the
corresponding modular characteristic cocycles:
$$\alignat2
\a_{\tilde g}(v(g\inv m g))&=\la_\a(m; g, s)v(m);\quad \b_{\tilde g}(w(g\inv m
g))&=\la_\b(m; g, s)w(m);\\
v(m)v(n)&=\mu_\a(m, n)v(mn);\qquad\quad w(m)w(n)&=\mu_\b(m, n)w(mn),
\endalignat
$$
with $\tilde g=(g, s)\in \wtH$ and $m, n \in L$. Now we take a closer look at
$(\la_\b, \mu_\b)\in \tZ_\a(\wtH, L, A)$:
$$\aligned
\mu_\b&(m, n)=w(m)w(n)w(mn)^*, \quad m, n \in L,\\
&=(u(m)v(m))(u(n)v(n))(u(mn)v(mn))^*\\
&=u(m)\a_m(u(n))v(m)v(n)v(mn)^*u(mn)^*\\
&=u(m)\a_m(u(n))\mu_\a(m, n)u(mn)^*\\
&=\mu_\a(m, n)u(m)\a_m(u(n))u(mn)^*\\
&=\mu_\a(m, n)\mu(m, n),
\endaligned
$$
and for $(m, \tilde g)=(m, g, s)\in L\times \wtH$ 
$$\aligned
\la_\b&(m; \tilde g)u(m)v(m)=\la_\b(m; \tilde g)w(m)\\
&=\b_{\tilde g}(w(g\inv m g))\\
&=\Ad(u(g))\scirc \a_{\tilde g}\Big(u(g\inv mg)v(g\inv mg)\Big)\\
& =u(g)\a_g(u(g\inv m g))\la_\a(m; \tilde g)v(m)u(g)^*\\
&=\mu(g, g\inv mg)u(mg)\la_\a(m; \tilde g)v(m)u(g)^*\\
&=\la_\a(m; \tilde g)\mu(g, g\inv mg)\mu(m, g)^*u(m)\a_m(u(g))v(m)u(g)^*\\
&=\la_\a(m; \tilde g)\mu(g, g\inv mg)\mu(m, g)^*u(m)v(m),
\endaligned
$$
Therefore the characteristic cocycles $(\la_\b, \mu_\b)\in \tZ_{\a\times\th}(\wtH, L,
A)$ is of the form:
$$\aligned
\mu_\b(m, n)&=\mu(m, n)\mu_\a(m, n), \quad m, n \in L;\\
\la_\b(m; \tilde g)&=\mu(g, g\inv mg)\overline{\mu(m, g)}\la_\a(m; \tilde g), \quad 
\tilde g=(g, s)\in \wtH.
\endaligned
$$
Thus we conclude that $\chim(\chi(\b))=\Res([\mu])\chim(\chi(\a))$ in
$\La_{\a\times\th}(\wtH, L, M, A)$. In virtue of Theorem 2.7, this is equivalent to
the fact that 
$$
\Obm(\dot\a)=\Obm(\dot \b)\in \thasth(\wtQ, A)*_{ \fs }\Hom_G(N, \thth1).
$$

iii) Suppose that $\sM$ is an infinite AFD factor and $H$ is amenable. The
automorphism $\sig\in \Aut_\th(\sC)$ can be extended to an automorphism in
$\Autt(\tM)$ by \cite{ST3} which will be denoted again by $\sig$. Replacing 
$\{\b_g: g\in H\}$ by $\{\sig\scirc\b_g\scirc\sig\inv: g\in H\}$, we may and
do assume that $\mod(\a)=\mod(\b)$ and $\d(\chim(\a))=\d(\chim(\b))$ in the
invariant group $\thasth(\wtQ, A)*_{ \fs }\Hom_G(N, \thth1)$. The 
modified {\hjr} of Theorem 2.7 yields the existence of a cohomology class
$[\mu]\in \thtw(H, \T)$ such that
$$
\chim(\b)=\Res([\mu])\chim(\a)\in \La_{\a\times\th}(\wtH, L, M, A).
$$
A cocycle perturbation of $\a$, denoted by $\a$ again, leaves a subfactor $\sB$ of
type {\oneinf} pointwise invariant. Let $u: g\in H\mapsto u(g)\in \sU(\sB)$ be a
projective unitary representation of $H$ in $\sB$ with the multiplier
$\mu\in\tztw(H, \T)$ representating $[\mu]$ such that
$$
u(g)u(h)=\mu(g, h)u(gh), \quad g, h \in H.
$$
Set ${}_u\a_g=\Ad(u(g))\scirc \a_g, g\in H$. Then it is a straightforward calculation
to show that $\chim({}_u\a)=\Res([\mu])\chim(\a)$. Therefore, the characteristic
invariant $\chim({}_u\a)$ is precisely  $\chim(\b)$ of $\b$. Hence the cocycle
conjugacy classification theorem, \cite{KtST1}, guarantees the concycle conjgacy of
${}_u\a$ and $\b$. Therefore, the original actions $\a$ and $\b$ are outer
conjugate.
\QED

\enddemo

\head{\bf \S4. Model Construction}
\endhead
As laid down in \cite{KtST1}, the construction of a model from a set of invariants is
an integral part of the classifiction theory. It is particulary important here because the
invariants associated with outer actions do not form a standard Borel space. For
example, the classification functor cannot be  Borel in the case of type \threez. 
So we have to begin with a desingularization of the space of invariants. 

We fix an ergodic flow $\{\sC, \R, \th\}$ to begin with. An action $\a$ of a group $G$ on the
flow $\{\sC, \R, \th\}$ means a homomorphism $g\in G\mapsto \a_g\in \Aut_\th(\sC)$, where
$\Aut_\th(\sC)=\{\sig\in \Aut(\sC): \sig\scirc\th_s=\th_s\scirc \sig, s\in \R\}$.

As before, we denote the
unitary group $\sU(\sC)$ of $\sC$ by $A$ for short.  The first cohomology group
$\thth1=\thth1(\R, A)$ can not be a standard Borel space if the flow $\th$ is properly
ergodic. So we have to consider the first cocycle group $\tzth1=\tzth1(\R, A)$ instead
together with the coboundary subgroup $\tbth1=\tbth1(\R, A)$. Next we fix a
countable discrete amenable group $G$ and an exact sequence:
$$\CD
1@>>>N@>>>G@>\pi>\underset{\fs}\to\longleftarrow>Q@>>>1
\endCD
$$
together with a cross-section $\fs$ which will be  fixed throughout as in the previous
section and therefore the $N$-valued cocyle $\fnn:$
$$
\fnn(p, q)=\fs(p)\fs(q)\fs(pq)\inv, \quad p, q \in Q,
$$
is also fixed. Let $\Hom_\R(Q, \Aut(\sC))$ be the set of all
homomorphisms $\a\!: p\in Q\mapsto \a_p\in \Aut_\th(\sC)$ from $Q$ into the group of
all automorphisms of $\sC$ commuting with the flow $\th$. It is easily seen to be a
Polish space. Each $\a\in \homr(Q, \Aut(\sC))$ can be identified with an action of $G$ 
whose kernel contains $N$. So we view $\a$ as an action of $G$ on $\sC$ freely
whenever necessary. We also use the notations $\wtQ=Q\times \R$ and
$\wtG=G\times \R$ freely. We fix the action $\a$ of $Q$ and consequently of $G$ on the flow
$\{\sC, \R, \th\}$ throughout this section and the joint action $\a\times \th$ will be denoted by
the single character $\a$ for short. With these data, we have the group of modular
obstructions: $\thasout(G\times\R, N, A)$ with $A=\sU(\sC)$ which will be fixed throughout
this section. The group $\thasout(G\times\R, N, A)$ is not a standard Borel group in general,
in particular it will never be standard except trivial cases if  the flow $\th$ is properly
ergodic. Also there is no way to construct a model directly from an element 
$([c],\nu)\in \thasout(G\times\R, N, A)$ either. We must desingularize the group of invariants
first. To this end, we first consider the group $\tzasout(G, N, A)$ of modular
obstruction cocycles $(c, \z)$. However, being an obstruction  cocycle, $(c, \z)$ does
not allow us to construct an outer action of $G$ directly. We recall Corollary 2.17
to find a resolution system:
$$\CD
1@>>>M@>>>H@>\pig>>G@>>>1
\endCD
$$
with $H$ a countable discrete group such that
\roster
\item"i)" the normal subgroup $M$ is abelian{\rm;}
\item"ii)" relative to the modified {\hjr}{\rm:}
$$\left.\CD
\cdots@>>>\thtw(H, \T)@>\Res>>\La_\a(\wtH, L, M, A)\\
@.\\
@>\d>>\thasout(G\times\R, N, A)@>\Inf>>\ththr(H, \T)\\
@.@V\part VV@A\pi_G^*AA\\
@. \ththr(G, \T)@=\ththr(G, \T)
\endCD\right\}
$$
where $L=\pi_G^{-1}(N)$
the $\Inf$-image of $([c], \nu)$ vanishes, i.e., 
$$
\Inf([c], \nu)=1\in \ththr(H, \T).
$$
Hence there exists a modular characteristic invariant 
$$
\chi\in \La_\a(\wtH, L, M, A)
$$
such that $\d(\chi)=([c], \nu)$.
\endroster
We also recall that in this resolution proceedure we need several extra data. For
example the map $\part: \thasout(G\times\R, N, A)\mapsto \tzthr(G, \T)$ requires a choice of
$(a, f)\in \tcatw(G, A)\times \tztw(Q, A)$ so that
$c_G=\pi^*(c_Q)\partg(\pi^*(f)a)^*\in \tzthr(G, \T)$. But in any case we do have
a resolution system $\{H, \pig,  L, M\}$ of $([c], \nu)$. So instead of going through
all steps of desingularizations starting from the cocycle $(c, \z)\in \tzasout(G, N, A)$,
we move directly to $\{H, \pig,  L, M\}$ and call $(\la, \mu)\in \La_\a(\wtH, L, M,
A)$ a {\it resolution} of the modular obstruction $([c], \nu)\in \thasout(G\times\R, N, A)$ if
$$
\d([\la, \mu])=([c], \nu).
$$
If we begin with $(\la, \mu)\in \tZ_\a(\wtH, L, M, A)$, it is easy to see the
corresponding obtruction cocycle $(c, \z)\in \tzasout(G, N, A)$: 
\roster
\item"a)" First we  fix a cross-section $\sh\!: G\mapsto H$ of the map $\pig$;
\item"b)" With $\dfs(p)=\sh\scirc \fs(p), p\in Q,$ and $\fnl(p,
q)=\dfs(p)\dfs(q)\dfs(pq)\inv\in L$ we have
$$\aligned
\z(s;  n)&=\la(\sh(n); s), \quad n\in N;\\
d_c(s; q, r)&=\la(\fnl(q, r); s), \quad q, r \in Q, s \in \R;\\
c_Q(p, q, r)&=\la(\dfs(p)\fnl(q, r)\dfs(p)\inv; \dfs(p))\mu(\dfs(p)\fnl(q, r)\dfs(p)\inv,
\fnl(p, qr))\\
&\hskip 1in\times
\{\mu(\fnl(p, q), \fnl(pq, r)\}^*, \quad p, q, r\in Q.
\endaligned
$$
\endroster
We will  write $(c, \z)=\part_\dfs(\la, \mu)$. Let $\rsn(H, \pig, ([c],
\nu))$ be the set of all $(\la, \mu)\in \tZ_\a(\wtH, L, M, A)$ such
that $\d([\la, \mu])=([c], \nu)$. On the space $\tzasout(G,\allowmathbreak N, A)$ of
modular obstruction cocycles, the group $\tcatw(Q,
A)$ acts in the following way:
$$
(c, \z)\mapsto ((\partwtq b)c, \z),b\in \tcatw(Q, A),
$$
which does not change the cohomology class of $(c, \z)$. Also the
group 
$$
\tztw(H, \T)\times \tC^1(N, A)
$$ 
acts on $\tZ_\a(\wtH, L,
M, A)$ without changing the cohomology class of $\d_\fs(\la,
\mu)$:
$$
(\la, \mu) \mapsto ((\part_1 a)\la_\xi\la, (\part_2 a)\xi_L\mu), \quad
(\xi, a)\in \tztw(H, \T)\times \tC^1(N, A),
$$
where $(\la_\xi, \xi_L)$ is the characteristic cocycle given by
(2.10):
$$\aligned
\la_\xi(m; g, s)&=\xi(g, g\inv m g)\overline{\xi(m, g)}, \quad m\in
L, (g, s)\in \wtH;\\
\xi_L&=\text{\rm the restriction of }\xi\ \text{to }L\times L.
\endaligned
$$

Now as soon as we have a characteristic cocycle $(\la, \mu)$, we
have a covariant cocycle $\{\sM, H, \a^{\la, \mu}\}$ equipped with
a map
$u: m\in L\mapsto u(m)\in\tsU(\sM)$ such that
$$\aligned
u(m)u(n)&=\mu(m, n)u(mn), \quad m, n \in L;\\
\a_m^{\la, \mu}&=\tAd(u(m)), ;\\
\a_g^{\la, \mu}\scirc\th_s(u(g\inv m g))&=\la(m; g, s)u(m), \quad
(g, s)\in
\wtH;\\
\endaligned
$$
which therefore gives:
$$
\da_g^{\la, \mu}=\a_{\sh(g)}^{\la, \mu}, \quad g\in G,
$$
whose modular obstruction cocycle is precisely 
$ \d_\sh(\la, \mu) $.
The action of $a\in \tC^1(N, A)$ on $(\la, \mu)$ does not change
the action $\a^{\la, \mu}$ itself, but on the unitary family
$\{u(m): m\in L\}$ which is perturbed to $\{(au)(m)=a(m)u(m):
m\in L\}$. So this does not cause any interesting change. The
perturbation by
$\xi\in \tztw(H, \T)$ gives somewhat non trivial change on
$\a^{\la, \mu}$. Namely, what we need is to consider the left
regular $\xi$-projective representation, say $v^\xi\!: g\in H\mapsto
v^\xi(g)\in\sU(\ell^2(H))$, so that
$$
v^\xi(g)v^\xi(h)=\xi(g, h)v^\xi(gh), \quad g, h \in H.
$$
Now the new action $g\in H\mapsto \a_g^{\la, \mu}\otimes
\Ad(v^\xi(g))\in \Aut(\sM\botimes \sL(\ell^2(H))$ has the modular
characteristic cocycle $(\la_\xi \la,  \xi_L  \mu)$, which is of course
does not change the outer conjugacy class of the outer action
$\da^{\la, \mu}$ of $G$. The change caused by the action of $b\in
\tcatw(Q, A)$ is again absorbed by changing the unitary family
$\{u(\fnl(p, q)): p, q \in Q\}$ to $\{b(p, q)u(\fnl(p, q)): p,  q \in
Q\}$, which does not change the outer action $\da^{\la, \mu}$
itself. Therefore the scheme of model constructions looks like:
$$\CD
(\la, \mu)\in\tZ_\a(\wtH, L, M, A)@>>>\d_\sh(\la,
\mu)\in\tzasout(G, N, A)\\
@VVV@AAA\\
\a^{\la, \mu}\in \act(H, \sM)@>>>\da^{\la, \mu}=\a_\sh\in \oct(G, \sM)
\endCD
$$
where $\act(G, \sM)$ and $\oct(G, \sM)$ are respectively the
spaces of actions and outer actions of $G$ on $\sM$. Summerizing the
discussion, we get the following:
\proclaim{Theorem 4.1} Let $G$ be a {\cdag} and $N$ a normal
subgroup. Let $\{\sC, \R, \th\}$ be an ergodic flow and $\a$ an action of
$G$ on the flow $\{\sC, \R, \th\}$ with $\Ker(\a)\supset N$, i.e., $\a$ is
a homomorphism of $G$ into the group $\Aut_\th(\sC)$ of
automorphisms commuting with the flow $\th$ whith $\a_m=\id, m \in
N.$ Let $A$ denote the unitary group $\sU(\sC)$.

 For every modular obstruction cocycle $(c, \z)\in \tzasout(G, N,
A)$, there exists an amenable resolution system $\{H, L, M, \pig, \la, \mu\}$ with
$(\la, \mu)\in \tZ_\a(\wtH, L, M,\allowmathbreak A)$ and a cross-section $\sh\!:
G\mapsto H$ of the map $\pig$ such that 
$$
\d_\sh(\la, \mu)\equiv (c, \z)\quad \mod\ \tbasout(G, N, A).
$$
Consequently, the action $\a^{\la, \mu}$ associated with the characteristic cocycle
$(\la, \mu)$ gives an outer action $\da^{\la, \mu}=\a_\sh$ of $G$ on the
approximately finite dimensional factor $\sM$ with flow of weights $\{\sC, \R,
\th\}$ such that
$$
\Obm(\da^{\la, \mu})=([c], [\z])\in \thasout(G\times\R, N, A).
$$
The homomorphism $\nu=[\z]\in \Hom(N, \thth1)$ is injective if and only if $\da$ is
free.
\endproclaim

\head{\bf \S 5. Non-Triviality of the Exact Sequence:
$$
\pmb{\CD1@>>>\thth1@>>>\Out(\sM)@>>>\Outt(\tM)@>>>1
\endCD}
$$}
\endhead
\proclaim{Theorem  5.1} Let $\a$ be an outer action of a countable
discrete group $G$
on a separable factor $\sM$ with $N=\a\inv(\cntr(\sM))$ with modular
obstruction
$$
\Obm(\a)=([c], \nu)\in \thaout(G, N, \sU(\sC)).
$$
Let $Q$ be the  quotient group $Q=G/N$ and $\fs$ be a cross-section
of the quotient
map $\pi\!: G\mapsto Q$. Then
the map $\a_\fs\!: p\in Q\mapsto \a_{\fs(p)}\in \Aut(\sM)$ can be
perturbed by $\cntr(\sM)$ to an action of $Q$ if and only if the
modular obstruction
$$
\Obm(\a)=([c], \nu)\in \thasout(G, N, A)=\tH_{\a, \txs}^3(\wtQ,
A)*_\fs \Hom_G(N,
\thth1)
$$
has  trivial 
$$
[c\cdot\a_p\left(\partial_Q(b)\right)]=
[c(\tp,\tq,\tr)\a_p(\partial_Q (b)(s;q,r))]=1
$$
for some $ b(\cdot,\ q)\in \tZ_\th^1(\Bbb R, A)$, which implies
$\nu\cup\fnn\in \tbatw(Q,
\thth1). $ 
\endproclaim

\demo{Proof} Suppose $ [c\cdot\a_p\left(\partial_Q(b)\right)]=1$
for some $b(\cdot, q)\in \tZ_\th^1(\Bbb R, A). $ Choose
$\{u(p, q)\in \tsU(\sM)\!:\ p, q \in Q\}$ so that
$$
\a_\tp\scirc \a_\tq=\Ad(u(p, q))\scirc \a_{\tp\tq}, \quad \tp, \tq\in \wtQ.
$$
The associated modular obstruction cocycle $(c, \nu)\in
\tzasout(\wtG, N, A)$ is given
by:
$$\aligned
c(\tp, \tq, \tr)&=\a_\tp(u(q, r))u(p, qr)\{u(p, q)u(pq, r)\}^*;\\
\nu(m)&=[\a_m]\in \thth1(\R, A), \quad m\in N.
\endaligned
$$
The triviality of $ [c\cdot\a_p\left(\partial_Q(b)\right)] $ means the existence of $f\in \tctw(Q, A)$ such that
$  c\cdot\a_p\left(\partial_Q(b)\right)=\part_\wtQ f $.
Setting $$  v(p, q)= f(p, q)^* w(p)\a_p(w(q))u(p, q)w(pq)^*, $$
where $  w(p)\in \tM$ with $\ w(p)^*\th_t(w(p))=b(t,p) ,$
we get
$$
  \Ad w(p)\scirc \a_\tp\scirc \Ad w(q)\scirc \a_\tq=
\Ad(v(p, q))\scirc\Ad w(pq)\scirc \a_{\tp\tq}\quad\text{and}\quad 
\part_\wtQ v=1.
$$
Since $v(q,r)^*\th_t(v(q,r))=1$ for $t\in \Bbb R$,  
the unitaries $v(q,r)$ are elements of $\sM.$ 
Setting 
$$
{}_w\a_p=\tAd w(p)\scirc\a_p,
$$ obtain
a cocycle crossed action $\{{}_w\a, v\}$ of $Q$ on $\sM $. As
the fixed point algebra
$ \sM^{{}_w\a} $ can be assumed to be properly infinite without loss of
generality, we can
find a family $ \{a(p)\in \sU(\sM): p\in Q\} $ such that
$$\aligned
 1&=a(p)\a_p(a(q))v(p, q)a(pq)^*\\
&=f(p, q)^*a(p){}_{w}\a_p(a(q))w(p)\a_p(w(q))u(p, q)w(pq)^*a(pq)^*\\
&=f(p, q)^*a(p)w(p)\a_p(a(q)w(q))u(p,q)\left(a(pq)w(pq)\right)^*;\\
f(p, q)&=a(p)w(p)\a_p(a(q)w(q))u(p,q)\left(a(pq)w(pq)\right)^*. 
\endaligned
$$
Hence $ \b={}_{a\cdot w}\a: \tp\in \wtQ\mapsto 
{}_{a\cdot w}\a_\tp=\Ad(a(p)w(p))\scirc
\a_\tp\in \Aut(\tM) $ is
an action of $\wtQ$ on $\tM$. The restriction of $\b$ to $\sM$ is precisely a
$\cntr(\sM)$-perturbation of the original action $\a$ since
$ \tAd(a(p)w(p))\in \cntr(\sM) $.
Now we have
$$\allowdisplaybreaks\aligned
  d_c(s, q, r)&=\th_s(u(q, r))u(q, r)^*\\
&=\th_s\Big(f(q, r)\a_q((a(r)w(r))^*)(a(q)w(q))^*a(pq)w(qr)\Big)\\
&\hskip.2in\times
\Big(f(q, r)\a_q((a(r)w(r))^*)(a(q)w(q))^*a(pq)w(qr)\Big)^*\\
&=(\partth f(q, r))_s \a_q(b(s, r)^*)b(s, q)^*b(s, pq) 
\endaligned
$$
where $ \th_s(a(p)w(p))(a(p)w(p))^*=\th_s(w(p))(w(p))^*=b(s, p)\in A. $  Hence  we have
$$\allowdisplaybreaks\aligned
\nu(\fnn(q, r))=[d_c(\cdot, q, r)]&=[\part_\wtQ(b(\cdot, \cdot)^*)(q, r)]\quad  \text{in}\
\thth1.
\endaligned
$$
Thus we conclude that $\nu\cup \fnn\in 
\tbatw(\wtQ, \thth1)$.

Conversely, suppose that $\a_\fs$ is perturbed to an action of $Q$ by
$\cntr(\sM)$.
Choose $\{  w(p) \in \tsU(\sM): p\in  Q\}$ so that
$$
\tAd(  w(p))\scirc \a_p\scirc \tAd(w(p))\scirc \a_q=\tAd(w(pq))\scirc
\a_{pq} , \quad p, q \in
Q.
$$
Let $\{u(p, q)\in \tsU(\sM): p, q \in Q\}$ be a family such that
$$
\a_p\scirc \a_q=\Ad(u(p, q))\scirc \a_{pq}.
$$
Then we have
$$
 f(p, q)=  w(p)\a_p(w(q))u(p, q)w(pq)^*\in A,
$$
and that
$$\aligned
(\part_\wtQ &f)(\tp, \tq, \tr)=\a_\tp\Big(w(q)\a_q(w(r))u(q,
r)w(qr)^*\Big)\\
&\hskip .5in\times
\Big(w(p)\a_p(w(qr))u(p, qr)w(pqr)^*\Big)\\
&\hskip .5in\times
\Big\{\Big(w(p)\a_p(w(q))u(p, q)w(pq)^*\Big)\\
&\hskip .5in\times
\Big(w(p)\a_{pq}(w(r))u(pq, r)w(pqr)^*\Big)\Big\}^*\\
&=c(\tp, \tq, \tr)\a_p\left(b(s,q)\a_q(b(s,r)b(s,qr)^*\right), 
\endaligned
$$
$  $ where $b(s,p)=w(p)^*\th_s(w(p))\in \tZ^1_\th(\Bbb R, A). $
Thus we conclude 
$$[c(\tp,\tq,\tr)\a_p(\partial_Q (b)(s;q,r))]=1. $$
\QED
\enddemo
This characterization has an immediate consequence:

\proclaim{Theorem 5.2} If $\sM$ is an approximately finite
dimensional factor of type
{\threee} with flow of weights $\{\sC, \R, \th\}$, then
the exact sequence{\rm:}
$$\CD
1@>>>\thth1(\R, \sU(\sC))@>>>\Out(\sM)@>>>\Outt(\tM)@>>>1
\endCD
$$
does not split.
\endproclaim

\demo{Proof} Let $G$ be the discrete Heisenberg group:
$$
G=\left\{\pmatrix 1&a&c\\0&1&b\\0&0&1\endpmatrix:  a, b, c\in \Z \right\}
$$
and
$$
N=\left\{\pmatrix 1&0&c\\0&1&0\\0&0&1\endpmatrix: c\in \Z \right\}
$$
be the center  of $G$ as in Example 7.1 on \cite{KtST1}. We write an
element of $G$ as
$(a, b, c)\in \Z$ with the multiplication rule:
$$
(a, b, c)(a', b', c')=(a+a', b+b', c+c'+ab').
$$
We then form the quotient group $Q=G/N$ and obtain an
exact sequence:
$$\CD
1@>>>N@>>>G@>\piq>>Q@>>>0.
\endCD
$$
The quotient group $Q$ is isomorphic to $\Z^2$. Define a
cross-section $\fs$ of $\piq$ in the
following way:
$$
\fs(a, b)=\pmatrix 1&a&0\\0&1&b\\0&0&1\endpmatrix, \quad (a, b)\in \Z^2=Q,
$$
and compute
$$\aligned
 \fnn (a, b&; a', b')
=ab'\in \Z
\endaligned\tag{5.1}
$$
where $N$ is identified with $\Z$.
Choose $(\la, \mu)\in \tZ(G, N, \T)$ to be trivial, i.e.,
$\mu(m, n)=1, m, n\in N$ and $\la(m, g)=1, m\in N, g\in G$. But
choose $\nu=e^{iT}\in
\Hom(N, \T)=\widehat N=\T$  with $T>0$ to be determined, so that the
characteristic cocycle
$(\la, \mu)\in \tZ(G\times \R, N, \T)$ is given by:
$$\allowdisplaybreaks\aligned
&\mu(m, n)=1, \quad m, n \in N;\\
&\la(m, (g, s))=\exp({\txti T'ms}), \quad  s\in \R, g \in G,
\endaligned
$$
where $T'={2\pi}/T.$
Let $\sM$ be an AFD factor of type {\threee} with flow of weights
$\{\sC, \R, \th\}$. Viewing
the torus $\T$ as the subgroup $\sU(\sC^\th)$ of the unitary group
$A=\sU(\sC)$, we view
the cocycle $(\la, \mu)$ as an element of $\tZ_\a(\wtG, N, A)$.
Now choose $T>0$ such
that $\sig_{nT}^\f\not\in \Int(\sM)$ for every $n\in \Z, n\neq 0,$ with $\f$ a
preassigned faithful semi-finite normal weight on $\sM$. Such a $T\in
\R$ exists because $\{t\in \R, \sigft\in\Int(\sM)\}$ must be a meager
subgroup of
$\R$. Let $\a=\a^{\la, \mu}$ be the action of $G$ on $\sM$ associated
with the cocycle
$(\la, \mu)$ and $\mod (\a_g)=\id$. The construction yields that the 
action $\a$ is free
and it enjoys the following
property:
$$
\a_m=\sig_{mT}^\f,\quad m\in N,
$$
with $\f$ a dominant weight on $\sM$.  We can assume the invariance
$\f=\f\scirc\a_g,
g\in G$ for $\a$. The freeness of $\a$ shows that the map:
$\da: g\in G\mapsto \da_g=[\a_g]\in \Out(\sM)$ is an injective
homomorphism such
that $\da_m=\sig_{mT}\in  \thth1(\R, A) \subset \Out(\sM)$. We are now
going to compute
the modular obstruction cocycle $(c, \nu)\in\tzasout(\wtQ, A)*_\fs
\Hom(N, \thth1).$ Since
$$\aligned
\a_\tp\scirc \a_\tq&=\a_{\fs(p), s}\scirc\a_{\fs(q),
t}=\a_{\fs(p)\fs(q)}\scirc\th_{s+t}\\
&=\a_{\fnn(p, q)\fs(pq)}\scirc\th_{s+t}\\
&=\Ad(\f^{\txti T'\fnn(p, q)})\scirc \a_{\tp\tq},
\endaligned
$$
with $u(p, q)=\f^{\txti T'\fnn(p, q)}$ we get
$$\allowdisplaybreaks\aligned
c(\tp, \tq, \tr)&=\a_\tp(u(q, r))u(p, qr)\{u(p, q)u(pq, r)\}^*\\
&=\th_s(u(q, r))u(q, r)^*\a_p(u(q, r))u(p, qr)\{u(p, q)u(pq, r)\}^*\\
&=\exp(-\txti T's\fnn(q, r))\f^{\txti T'(\fnn(q, r)+\fnn(p,
qr)-\fnn(p, q)-\fnn(pq, r))}\\
&=\exp(-\txti T's\fnn(q, r))
\endaligned
$$
and $c_Q=1$.
In order for $  c\cdot \a_p(\partial_Q(b))$ with some
$b(\cdot,q)\in \tZ_\th^1(\Bbb R, A) $ to be trivial, it is necessary 
and sufficient that
there exists
$f\in \tctw(Q, A)$ such that $\part_\wtQ f=
  c\cdot \a_p(\partial_Q(b))$.
 The function $f$
satisfies the equations:
$$
\exp(-\txti T's\fnn(q, r))= \a_q(b(s,r)^*)b(s,q)^*b(s,qr)f(q, r)^*\th_s(f(q,r)) 
$$
which means that $[\exp(-\txti T's\fnn(q, r))]\in 
 \tB^2(Q, \tH^1_\th). $ As
$\mod(\a_p)=\id, p \in  Q,$ and $ Q$ is a free abelian group, the second
cohomology group $\thtw(Q, \tH^1_\th)$ is isomorphic to the group
$X(Q^2, \tH^1_\th)$ of all $\tH^1_\th$-valued 
skew symmetric bihomomorphisms.
We have
$$\aligned
  \exp\Big(&-\txti T's\fnn(q, r)\Big)\exp\Big(\txti T's\fnn(r, q)\Big) \\
=&\Big(\a_q(b(s,r)^*)b(s,q)^*b(s,qr)f(q, r)^*\th_s(f(q,r))\Big)\\
&\hskip.2in \times
\Big(\a_r(b(s,q)^*)b(s,r)^*b(s,rq)f(r, q)^*\th_s(f(r,q))\Big)^*\\
=&f(r, q)f(q, r)^*\th_s(f(r, q)^*f(q, r)).
\endaligned\tag{5.2}
$$
By (5.1), we have
$$\exp(-\txti T's\fnn(q, r)+\txti T's\fnn(r, q))
=\exp(-\txti T's(ab'-a'b))
$$
where $q=(a, b), r=(a', b')$ Thus it follows from (5.2) that
the modular automorphisms $\sig_{T(ab'-a'b)}^\f$ are
inner, which contradicts to the choice of $T$. Therefore
$[c(\tp,\tq,\tr) \a_p(\partial_Q (b)(s;q,r))]\neq 1$
in $\tH_{\a,
\txs}^3(\wtQ, A)$. Theorem 5.1 says that $\a_\fs$ cannot be perturbed
into an action of
$Q$ by $\cntr(\sM)$.

Since we have a commutative diagram of exact sequences:
$$\CD
1@>>>N@>>>G@>\piq>\underset{\fs}\to\longleftarrow>Q@>>>1\\
@.@V\nu VV@V\da VV@V\ta VV\\
1@>>>\thth1@>>>\Out(\sM)@>\pi>\underset{\sp}\to\longleftarrow>\Outt(\tM)@>>>1,
\endCD
$$
if the second sequence splits via cross-section $\sp$, then the
associated injection $\ta$ of
$Q$ into $\Outt(\tM)$ is composed with the cross-section $\sp$ to
be an outer action of
$Q$, say $\b$. But $\ththr(Q, \T)=1$,  so that $\b$ can be perturbed
into an action of
$Q$, denoted by $\b$ again.
Then we have  $\b_p\equiv \a_{\fs(p)}\ \mod\ \cntr(\sM)$. Therefore
the $\a_\fs$ is
perturbed to an action by $\cntr(\sM)$, which contradicts to the fact
$[c(\tp,\tq,\tr) \a_p(\partial_Q (b)(s;q,r))]\neq 1$ for any
$  b(\cdot, p) \in\tZ^1(\Bbb R, \tH^1_\th)$ as seen
above.
\QED
\enddemo

\end{document}